\documentclass[preprint,12pt]{elsarticle}

\usepackage[numbers]{natbib}
\usepackage{color}
\usepackage{amssymb}
\usepackage{mathrsfs}
\usepackage{bm}
\usepackage{psfrag}
\usepackage{amsthm}
\usepackage{appendix}
\usepackage{algorithm}
\usepackage{algorithmic}
\usepackage{epsf,amsfonts,amsmath}
\usepackage{amssymb,amsmath}
\usepackage{graphicx}

\newtheorem{assump1}{Assumption}
\newtheorem{remark1}{Remark}
\newtheorem{prop1}{Proposition}

\journal{Elsevier}

\begin{document}

\begin{frontmatter}

\title{Partitioned treatment of uncertainty in coupled domain problems: A separated representation approach}

\author[UCB1]{Mohammad Hadigol}
\ead{mohammad.hadigol@colorado.edu}

\author[UCB1]{Alireza Doostan\corref{cor1}}
\ead{alireza.doostan@colorado.edu}

\author[TUB2]{Hermann G. Matthies}
\ead{h.matthies@tu-braunschweig.de}

\author[TUB2]{Rainer Niekamp}
\ead{r.niekamp@tu-bs.de}

\cortext[cor1]{Corresponding Author: Alireza Doostan}

\address[UCB1]{Aerospace Engineering Sciences Department, University of Colorado, Boulder, CO 80309, USA}

\address[TUB2]{Institute of Scientific Computing, Technische Universit\"at Braunschweig, Braunschweig, Germany}

\begin{abstract}
\label{Abstract}

This work is concerned with the propagation of uncertainty across coupled domain problems with high-dimensional random inputs. A stochastic model reduction approach based on low-rank separated representations is proposed for the partitioned treatment of the uncertainty space. The construction of the coupled domain solution is achieved though a sequence of approximations with respect to the dimensionality of the random inputs associated with each individual sub-domain and not the combined dimensionality, hence drastically reducing the overall computational cost. The coupling between the sub-domain solutions is done via the classical Finite Element Tearing and Interconnecting (FETI) method, thus providing a well suited framework for parallel computing. Two high-dimensional stochastic problems, a 2D elliptic PDE with random diffusion coefficient and a  stochastic linear elasticity problem, have been considered to study the performance and accuracy of the proposed stochastic coupling approach.

\end{abstract}


\end{frontmatter}


%
\section{Introduction}
\label{sec:introduction}
Simulation-based prediction of most physical systems is subject to either lack of knowledge about the governing physical laws or incomplete/limited information about model parameters such as material properties,  initial, or boundary conditions. In order to obtain realistic predictions of these systems one, therefore, needs to characterize such uncertainties and quantify their impact on Quantities of interest (QoI). Uncertainty quantification (UQ), an emerging field in computational engineering and science, is concerned with the development of rigorous and efficient solutions to this exercise. 

A major class of UQ approaches are probabilistic where uncertain parameters are represented by random variables or processes. Several techniques have been developed to study the propagation of such uncertainties in engineering systems, e.g., see \cite{Ghanem91a,HGMatth08,Najm09,Xiu09a,Xiu10a,LeMaitre10} and the references therein. Among these techniques, stochastic spectral methods \cite{Ghanem91a,Ghanem99c,Xiu03,Xiu10a,LeMaitre10} based on polynomial chaos (PC) expansions \cite{Wiener38,Cameron47} have received special attention due to their advantages over traditional UQ techniques such as perturbation-based and Monte Carlo sampling (MCS) methods. In particular, under certain regularity conditions, these schemes converge faster than MCS methods \citep{LeMaitre10} and, unlike perturbation methods, are not restricted to problems with small uncertainty levels \citep{Ghanem91a}. Stochastic spectral methods are based on expanding the solution of interest in PC bases. The coefficients of these expansions are then computed, for instance, via Galerkin projection \citep{Ghanem91a}, referred to as stochastic Galerkin (SG), or pseudo-spectral collocation \cite{Xiu05a,Mathelin03}, named stochastic collocation (SC).

Although PC-based techniques benefit from elegant mathematical analyses, e.g.,  formal convergence studies, {\color{black}they may suffer from the so-called \textit{curse-of-dimensionality} if executed carelessly: the computational cost may grow exponentially as a function of the number of independent random inputs \citep{Xiu10a,LeMaitre10,Doostan11a}, and as a function of the total number of unknowns in the system. This may be a common situation when one is dealing with UQ of engineering problems involving, for instance, coupled domains or separated scales with independent and separate sources of uncertainty. In the past few years, several alternative methods have been proposed to introduce some form of {\it sparsity} in an effort to counter the curse-of-dimensionality \cite{Mathelin03,Xiu05a,Todor07a,Bieri09a,Bieri09c,Nouy08,Nobile08b,Doostan09,Ma09a,Foo10,Gao10b,Nouy_2010,Bieri11,Zhang11,Doostan11a,Khoromskij11,Chinesta11,Falco11,Gao11a,Hackbusch12,Doostan13}. The ideas used here are based on sparse grids, element-like partitioning, ANOVA expansions, separated representations, and more general tensor decompositions. Although these techniques have been found efficient in alleviating curse-of-dimensionality, the high computing cost is still a major bottleneck in PC expansion of high-dimensional random solutions. }

Of particular interest in the present study is the problem of uncertainty propagation across coupled domain problems where independent, high-dimensional random inputs are present in each sub-domain. In such cases, the solution depends on all random inputs; hence, a direct application of PC expansion may not be feasible {\color{black} or at least not desirable}. While integration of PC expansions with standard domain decomposition (DD) techniques may partially reduce the overall computational complexity by partitioning the physical space, e.g., see \cite{Ghosh09,Cottereau11,Subber12}, expansions (or sampling) with respect to the combined set of random inputs is still required. Instead, we propose an approach that additionally enables a {\it partitioned treatment} of the stochastic space; that is, the solution is computed through a sequence of approximations with respect to the random inputs associated with each individual sub-domain, and not the combined set of random inputs. To this end, we adopt a stochastic expansion based on the so-called {\it separated representations} and demonstrate how it can be obtained in conjunction with a DD approach. {\color{black}As one of the simplest cases of such coupling, we consider a linear problem on two non-overlapping sub-domains with a common interface. Although more elaborate coupling techniques are possible, we restrict ourselves to a finite element tearing and interconnecting (FETI) approach \cite{Farhat91} for the sake of simplicity.}

Model reduction techniques based on separated representation, a.k.a. canonical decomposition, of high-dimensional stochastic functions have been recently proposed in \citep{Nouy07,Doostan07,Doostan09,Nouy08,Nouy_2010,Khoromskij11,Hackbusch12}  and proven effective in reducing the issue of curse-of-dimensionality. We here adopt a special form of separated representations for the stochastic computation of coupled domain problems. 

Let $\left(\Omega,\mathcal{T},\mathcal{P} \right)$ be a complete probability space where $\Omega$ is the sample set and $\mathcal{P}$ is a probability measure on the $\sigma-$field $\mathcal{T}$. Also assume {\color{black}that the input uncertainty has been discretized and approximated by random variables, such that the vector $\bm{\xi} = \left(\xi_1,\cdots, \xi_d\right):\Omega \rightarrow \mathbb{R}^{d}$, $d \in \mathbb{N}$, represents the set of independent random inputs associated with a PDE defined on a domain $\mathcal{D}\subset\mathbb{R}^D$, $D\in\{1,2,3\}$, composed of two non-overlapping sub-domains $\mathcal{D}_1$ and $\mathcal{D}_2$}. We further assume that the random vector $\bm{\xi} = (\bm{\xi}_1, \bm{\xi}_2)$ is such that $\bm{\xi}_1:\Omega \rightarrow \mathbb{R}^{d_1}$ and $\bm{\xi}_2:\Omega \rightarrow \mathbb{R}^{d_2}$ denote random inputs corresponding to $\mathcal{D}_1$ and $\mathcal{D}_2$, respectively. Here, $d_1,d_2 \in \mathbb{N}$ are the sizes of $\bm \xi_1$ and $\bm \xi_2$, respectively, and $d=d_1+d_2$. Assuming that $u_1(\bm x, \bm{\xi}):\bar{\mathcal{D}}_1\times\Omega\rightarrow \mathbb{R}$ and $u_2(\bm x,\bm{\xi}):\bar{\mathcal{D}}_2\times\Omega\rightarrow \mathbb{R}$ are the sub-domain solutions, we consider the separated representation of the form
\begin{equation}
\label{eqn:sr_general}
u_i(\bm x,\bm{\xi}) = \sum_{l=1}^r {u}_{0,i}^l(\bm x) \phi_1^l(\bm{\xi}_1) \phi_{2}^l(\bm{\xi}_{2}) + \mathcal{O}(\epsilon),\qquad i=1,2.
\end{equation}
Here $u_{0,i}^l(\bm{x}):\bar{\mathcal{D}}_i\rightarrow\mathbb{R}$ and $\phi_i^l(\bm{\xi}_i):\Omega\rightarrow\mathbb{R}$, $l=1,\dots,r$, are, respectively, deterministic and stochastic functions -- or {\it factors} --  to be determined along with the {\it separation rank} $r$. These quantities are not fixed \textit{a priori} and are computed through an optimization scheme such that a prescribed target accuracy $\epsilon$ is reached for a minimum $r$. 

{\color{black}Notice that, by definition, finding the expansion (\ref{eqn:sr_general}) is in fact a non-linear problem, even for a linear problem}. However, it can be computed through a sequence of alternating linear problems, where $\{u_{0,i}^l(\bm{x})\}_{l=1}^{r}$, $i=1,2$, $\{\phi_1^l(\bm{\xi}_1)\}_{l=1}^{r}$, or $\{\phi_2^l(\bm{\xi}_2)\}_{l=1}^{r}$ are solved for one at a time while others are fixed at their recent values. This alternating construction, together with the separated form of (\ref{eqn:sr_general}), enables computing the stochastic functions $\{\phi_1^l(\bm{\xi}_1)\}_{l=1}^{r}$ and $\{\phi_2^l(\bm{\xi}_2)\}_{l=1}^{r}$ with {\color{black}computational} complexities that depend on $d_1$ and $d_2$ but not $d=d_1+d_2$. Additionally, it allows for a natural extension of DD techniques for computing $\{u_{0,i}^l(\bm{x})\}_{l=1}^{r}$. In the present work, we employ the standard FETI approach \citep{Farhat91} for the latter purpose. Moreover, for situations where the separation rank $r$ is small, (\ref{eqn:sr_general}) provides a reduced order approximation and representation of the coupled solution. We will describe the details of computing (\ref{eqn:sr_general}) in Section \ref{sec:SR}.

Among the limited earlier effort on separating random inputs for the solution of coupled problems, we particularly mention the recent work of Arnst et al. \citep{Arnst12} that is based on the so-called {\it reduced chaos} expansions \citep{Soize09}. In particular, in each component, e.g., physics or sub-domain, the solution is expanded in a PC basis that is generated based on random inputs associated with that component. The PC coefficients are then considered to be functions of random inputs corresponding to the other component. {\color{black}This may be seen as a special case of (\ref{eqn:sr_general}).} As in (\ref{eqn:sr_general}) we do not prescribe a stochastic basis {\it a priori}, our approach is different from that of \citep{Arnst12}. The construction of (\ref{eqn:sr_general}) is also different from that of the reduced chaos expansions. More importantly, when the coupled solution admits a low separation rank $r$, the computational complexity and storage requirement of (\ref{eqn:sr_general}) may be significantly smaller than those of the reduced chaos expansion, especially when $d$ is large. 

{\color{black}Separated representations have also been utilized in the context of multi-scale domain decomposition with high-dimensional local uncertainties \cite{Chevreuil12}. While our coupled domain formulation is similar to that of \cite{Chevreuil12}, our particular choice of separated representation (\ref{eqn:sr_general}) and it's numerical construction, both at the physical and stochastic levels, are considerably different from those in \cite{Chevreuil12}.}

The remainder of this paper is organized as follows. Section \ref{sec:pf} summarizes the problem formulation. There, we start with a general linear stochastic PDE as an abstract problem and present its coupled formulation via FETI. {\color{black}As already remarked, this is one of the simplest problems which has been chosen for the sake of simplicity. The techniques presented here are nonetheless applicable also to non-linear problems and more elaborate coupling schemes, as well as to more than two coupled components.} The choice of spatial discretization is given in Section \ref{sec:Det_Dis} and following that, in Section \ref{sec:SG_Scheme}, we review the application of the SG scheme to the original problem as well as its coupled formulation. In Section \ref{sec:SR}, we present our approach based on separated representations along with an extension of the FETI algorithm to solve the resulting coupled problems at the spatial level. Finally, two numerical examples; a 2D linear diffusion problem on an L-shaped domain with random diffusivity coefficient and a 2D linear elasticity problem describing the deflection of a cantilever beam with random Young's modulus are considered in Section \ref{sec:Numerical_Examples}.  The accuracy and performance of the proposed coupling approach are demonstrated using these two examples. 

\section{Problem Formulation}
\label{sec:pf}
%

\subsection{A linear PDE with random inputs}
\label{sec:original_problem}

\begin{figure}[htb]
  \begin{center}
	\includegraphics [width=7cm]{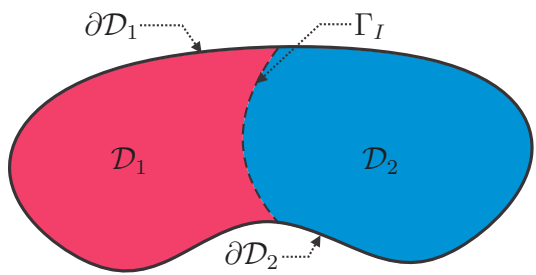}
    \put(-160,40){$\mathcal{D}_1$}
    \put(-65,40){$\mathcal{D}_2$}
    \put(-170,91){$\partial \mathcal{D}_1$}
    \put(-117,4){$\partial \mathcal{D}_2$}
    \put(-68,91){$\Gamma_I$}
    \caption{Geometry of the original problem and partitioning of $\mathcal{D}$ into two non-overlapping sub-domains $\mathcal{D}_1$ and $\mathcal{D}_2$. $\mathcal{D}=\mathcal{D}_1\cup\mathcal{D}_2$, $\mathcal{D}_1\cap\mathcal{D}_2=\emptyset$, and $\partial\mathcal{D}_1\cap\partial\mathcal{D}_2=\Gamma_I$.}
    \label{fig:DD-schematic}
   \end{center}
\end{figure}

Let $\left(\Omega,\mathcal{T},\mathcal{P} \right)$ be a suitable probability space. We consider computing the solution $u(\bm{x},\bm{\xi}):\bar{\mathcal{D}} \times \Omega \rightarrow \mathbb{R}$ satisfying the linear PDE 
\begin{equation} 
\label{eqn:full_SPDE}
\begin{aligned}
\mathcal{L}(\bm{x},\bm{\xi})(u(\bm{x},\bm{\xi}))=f(\bm{x})\ \ \ \ \ \bm{x} \in \mathcal{D}, \\
\mathcal{B}(\bm{x},\bm{\xi})(u(\bm{x},\bm{\xi}))=g(\bm{x})\ \ \ \ \bm{x} \in \partial \mathcal{D},\\
\end{aligned}
\end{equation}
$\mathcal{P}-$ almost surely in $\Omega$. Here, $\mathcal{D} \subset \mathbb{R}^D$, $D$=1,2,3, is the spatial extent of the problem with boundary $\partial \mathcal{D}$ as displayed in Fig. \ref{fig:DD-schematic}, $\mathcal{L}$ is a linear differential operator, and $\mathcal{B}$ is a boundary operator taking possibly various forms on different boundary segments. {\color{black}We assume (\ref{eqn:full_SPDE}) is the variational equation associated with a minimization problem. Additionally, as mentioned before, the uncertainty in (\ref{eqn:full_SPDE}) is assumed to have been approximated in such a way it may be represented by the random vector $\bm{\xi} = \left(\xi_1,\cdots, \xi_d\right):\Omega \rightarrow \mathbb{R}^{d}, d \in \mathbb{N}$, consisting of independent identically distributed (i.i.d.) random variables $\xi_k$ with probability density function $\rho(\xi_k):\mathbb{R}\rightarrow\mathbb{R}_{+}$}. Without loss of generality, the functions $f(\bm{x})$ and $g(\bm{x})$ in (\ref{eqn:full_SPDE}) are considered to be deterministic. We hereafter refer to (\ref{eqn:full_SPDE}) as the \textit{original problem}. 

We assume that $\mathcal{D}$ is composed of two non-overlapping sub-domains $\mathcal{D}_1$ and $\mathcal{D}_2$ ($\mathcal{D}=\mathcal{D}_1\cup\mathcal{D}_2$ and $\mathcal{D}_1\cap\mathcal{D}_2=\emptyset$) sharing an interface boundary $\Gamma_I=\partial\mathcal{D}_1\cap\partial\mathcal{D}_2$ as shown in Fig. \ref{fig:DD-schematic}. We further assume that the random inputs representing the uncertainty in $\mathcal{D}_1$ and $\mathcal{D}_2$ are independent and denoted by $\bm{\xi}_1 = \left(\xi_{1,1},\cdots, \xi_{1,d_1} \right): {\color{black}\Omega \rightarrow \mathbb{R}^{d_1}}, d_1 \in \mathbb{N}$ and $\bm{\xi}_2 = \left(\xi_{2,1},\cdots, \xi_{2,d_2} \right): {\color{black}\Omega \rightarrow \mathbb{R}^{d_2}}, d_2 \in \mathbb{N}$, respectively. This, therefore, provides a natural partitioning of $\bm \xi$ into $\bm{\xi}=(\bm{\xi}_1,\bm{\xi}_2)$.  

We note that while we consider the case of two coupled sub-domains, the extension of the subsequent algorithms to the case of multiple coupled sub-domains is straight-forward. 

\subsection{The coupled formulation of original problem}
\label{sec:coupled_problem}

We next reformulate the original problem (\ref{eqn:full_SPDE}) in terms of the coupled stochastic PDEs,
\begin{eqnarray} 
\label{eqn:coupled_SPDE}
&&\mathcal{L}(\bm{x},\bm{\xi}_1)(u_1(\bm{x},\bm{\xi}))=f(\bm{x})\ \ \ \ \ \ \bm{x} \in \mathcal{D}_1, \nonumber \\
&&\mathcal{L}(\bm{x},\bm{\xi}_2)(u_2(\bm{x},\bm{\xi}))=f(\bm{x})\ \ \ \ \ \ \bm{x} \in \mathcal{D}_2, \nonumber \\
&&\mathcal{B}(\bm{x},\bm{\xi}_1)(u_1(\bm{x},\bm{\xi}))=g(\bm{x})\ \ \ \ \ \ \bm{x} \in \partial \mathcal{D}_1, \nonumber \\
&&\mathcal{B}(\bm{x},\bm{\xi}_2)(u_2(\bm{x},\bm{\xi}))=g(\bm{x}) \ \ \ \ \ \ \bm{x} \in \partial \mathcal{D}_2, \nonumber \\
&&u_1(\bm{x},\bm{\xi})=u_2(\bm{x},\bm{\xi})\quad\quad\ \ \ \ \ \ \ \ \ \bm{x} \in \Gamma_I, 
\end{eqnarray}
which we henceforth refer to as the \textit{coupled problem}. Here, $u_1(\bm{x},\bm{\xi}):\bar{\mathcal{D}}_1 \times \Omega \rightarrow \mathbb{R}$ and $u_2(\bm{x},\bm{\xi}):\bar{\mathcal{D}}_2 \times \Omega \rightarrow \mathbb{R}$ are the solutions to (\ref{eqn:coupled_SPDE}). The last equation in (\ref{eqn:coupled_SPDE}) enforces the continuity of the solution on the interface $ \Gamma_I$. {\color{black}Again, for the sake of simplicity, we tacitly assume that the last equation in (\ref{eqn:coupled_SPDE}) is sufficient to enforce that the two partial solutions $u_1(\bm{x},\bm{\xi})$ and $u_2(\bm{x},\bm{\xi})$ of (\ref{eqn:coupled_SPDE}) agree with the restriction of the solution $u(\bm{x},\bm{\xi})$ of (\ref{eqn:full_SPDE}) to the respective sub-domains.} Our goal is then to find the sub-domain solutions $u_1(\bm{x},\bm{\xi})$ and $u_2(\bm{x},\bm{\xi})$ assuming independent solvers for each sub-domain. 

In the following section, we will discuss the spatial discretization of the original and coupled problems. 

\section{Spatial Discretization}
\label{sec:Det_Dis}
\subsection{Original problem}
\label{sec:Det_Dis_original}

Let $\mathcal{X}$ be a suitable Hilbert space for the  spatial discretization of the solution to the variational formulation of (\ref{eqn:full_SPDE}). For example, if $\mathcal{L}$ in (\ref{eqn:full_SPDE}) is a linear elliptic differential operator {\color{black}of second order} and $\mathcal{B}$ corresponds to a homogeneous Dirichlet boundary condition, then $ \mathcal{X} := H^1_0(\mathcal{D}) = \lbrace v \in  H^1(\mathcal{D}): v=0$ on $\partial \mathcal{D} \rbrace$. Also let $\mathcal{W} := L_2(\Omega) = \lbrace v:\Omega\rightarrow \mathbb{R}: \mathbb{E}[v ^2] = \int_{\Omega}  v^2 \mathcal{P} (\mathrm{d}{v}) < \infty \rbrace$ be the space of all square integrable random variables defined on $\left(\Omega,\mathcal{T},\mathcal{P} \right)$. Here, $\mathbb{E}[\cdot]$ denotes the mathematical expectation operator. The solution $u(\bm{x},\bm{\xi})$ to (\ref{eqn:full_SPDE}) then lives in the tensor-product space $\mathcal{X} \otimes \mathcal{W}$ \cite{HGMatth08}.

A finite dimensional subspace $\mathcal{X}_h \subset \mathcal{X}$ may be constructed by considering a finite element (FE) discretization consisting of piecewise linear, continuous functions $\lbrace N_m(\bm{x}) \rbrace_{m=1}^{M}$ on a triangulation of $\mathcal{D}$ with a uniform mesh size $h$. The approximate solution $u_h$ to the original problem (\ref{eqn:full_SPDE}) can then be written in the form 
\begin{equation} 
\label{eqn:sdasop}
u_h(\bm{x},\bm{\xi})= \sum_{m=1}^{M} u_m(\bm{\xi}) N_m(\bm{x}).
\end{equation}

We denote by $\bm{u}(\bm{\xi}) = [u_1(\bm{\xi}),\cdots, u_{M}(\bm{\xi})]^T \in \mathbb{R}^M \otimes \mathcal{W}$ the random vector of expansion coefficients in (\ref{eqn:sdasop}) given by the Galerkin projection of (\ref{eqn:full_SPDE}) onto $\mathcal{X}_h$. Specifically, 
\begin{eqnarray} 
\label{eqn:sdop}
& \bm{K}(\bm{\xi}) \bm{u}(\bm{\xi}) = \bm{f},
\end{eqnarray} 
where $\bm{K}(\bm{\xi}) : \Omega \rightarrow \mathbb{R}^{M \times M}$ is the random stiffness matrix -- assumed to be {\color{black}symmetric, positive definite $\mathcal{P}-$ almost surely in $\Omega$ } -- and $\bm{f} \in \mathbb{R}^M$ is the deterministic force vector.  

We next give a brief discussion on a number of domain decomposition techniques as well as our approach for spatial  discretization of the coupled problem (\ref{eqn:coupled_SPDE}).

\subsection{Coupled problem}
\label{sec:Det_Dis_coupled}
Domain decomposition (DD) methods have been found efficient and powerful means for solving large-scale PDEs with parallel computing {\color{black}\citep{Chan90,Quarteroni99,Smith04,Toselli05,Gosselet06}}. The main idea of DD is to decompose $\mathcal{D}$ into a number of smaller overlapping or non-overlapping sub-domains, and to compute the solution on each sub-domain to obtain the solution over the original domain \citep{Xu98,Cai93}. {\color{black}This is a special case of a partitioned solution approach.} 

In the present study we assume that the partitioning of $\mathcal{D}$ results in non-overlapping sub-domains. Non-overlapping DD methods, also known as iterative sub-structuring methods \citep{Toselli05,Gosselet06}, have some advantages over overlapping DD schemes, e.g., in being more efficient in handling elliptic problems with large jumps in the coefficients \cite{Xu98}. These methods can be classified based on the treatment of the interface solution compatibility. Primal methods, such as balancing domain decomposition (BDD) \cite{Tallec91} and balancing domain decomposition by constraints (BDDC) \cite{Dohrmann03}, enforce the continuity constraint on the sub-domains interfaces via computing a unique unknown interface solution. Alternatively, in dual methods, e.g., the FETI algorithm \citep{Farhat91}, this is done by iterating on intermediary unknown variables, i.e., Lagrange multipliers. There are also hybrid methods such as dual--primal finite element tearing and interconnecting (FETI-DP) \cite{Farhat01}, that are based on the combination of these two approaches. In the present work, as mentioned in Section \ref{sec:introduction}, we use the classical FETI algorithm towards computing the coupled solution to (\ref{eqn:coupled_SPDE}). However, we note that our approach may also be integrated with other non-overlapping DD techniques, specially the FETI-DP. The FETI method is based {\color{black}on finding the stationary or saddle point of an energy functional} associated with the sub-domain solutions, which we will discuss next.

Let $\bm{u}_1(\bm{\xi})\in \mathbb{R}^{M_1} \otimes \mathcal{W}$ and $\bm{u}_2(\bm{\xi})\in \mathbb{R}^{M_2} \otimes \mathcal{W}$ be the vector of nodal values corresponding to the FE discretization of $u_1$ and $u_2$ in (\ref{eqn:coupled_SPDE}), respectively. For simplicity, we assume that the triangulations of $\mathcal{D}_1$ and $\mathcal{D}_2$ result in matching nodes at the interface $\Gamma_I$ with $M_I$ degrees of freedom. Let also $\bm{\lambda}(\bm{\xi}) = [\lambda_1(\bm{\xi}),\cdots, \lambda_{M_I}(\bm{\xi})]^T \in \mathbb{R}^{M_I} \otimes \mathcal{W}$ be the random vector of Lagrange multipliers to enforce the solution continuity at $\Gamma_I$. Then the triple $(\bm{u}_1,\bm{u}_2,\bm{\lambda})$ is the stationary point of the energy functional 
\begin{eqnarray}
\label{eqn:energy}
&&\pi(\bm{u}_1,\bm{u}_2,\bm{\lambda}) =  \\
&&\mathbb{E} \left[\frac{1}{2}\bm{u}_1^{T}\bm{K}_{1}\bm{u}_1-\bm{u}_1^{T}\bm{f}_{1}\right]+
\mathbb{E} \left[ \frac{1}{2}\bm{u}_2^{T}\bm{K}_{2}\bm{u}_2-\bm{u}_2^{T}\bm{f}_{2}\right]+ \mathbb{E} \left[ \bm{\lambda}^{T}\left(\bm{C}_{2}^T \bm{u}_2 - \bm{C}_{1}^T \bm{u}_1 \right) \right],\nonumber
\end{eqnarray}
where, for $i=1,2$, $\bm{C}_{i}\in \mathbb{R}^{M_i\times M_I}$ are matrices that extract the interface nodal solutions from $\bm{u}_i$ and $\bm{K}_{i}(\bm{\xi}_i) : \Omega \rightarrow \mathbb{R}^{M_i \times M_i}$ and $\bm{f}_{i}\in \mathbb{R}^{M_i}$, are, respectively, the random stiffness matrices and deterministic force vectors associated with the FE discretization of (\ref{eqn:coupled_SPDE}) on $\mathcal{D}_i$. 

\begin{assump1}
\label{spd_assumption}
In the present study, we assume that $\bm{K}_{1}(\bm{\xi}_1)$ is symmetric positive definite, while $\bm{K}_{2}(\bm{\xi}_2)$ is symmetric positive semi-definite (which includes definite) {\color{black}$\mathcal{P}-$almost surely in $\Omega$}. 
\end{assump1}

Given Assumption \ref{spd_assumption}, the coupled solution $(\bm{u}_1,\bm{u}_2,\bm{\lambda})$ can be obtained by solving the following saddle point problem, established by the first-order optimality condition corresponding to the stationary point of $\pi$ in (\ref{eqn:energy}),
\begin{eqnarray} 
\label{eqn:sdcp}
\left[\begin{array}{ccc}\bm{K}_{1}(\bm{\xi}_1) & \bm{0} & \bm{C}_{1} \\ \bm{0} & \bm{K}_{2}(\bm{\xi}_2) & -\bm{C}_{2} \\ {\bm{C}_{1}}^T & -{\bm{C}_{2}}^T  & \bm{0}\end{array}\right]\left\{\begin{array}{c}\bm{u}_1(\bm{\xi}) \\\bm{u}_2(\bm{\xi})  \\ \bm{\lambda}(\bm{\xi})\end{array}\right\}=\left\{\begin{array}{c}\bm{f}_{1} \\ \bm{f}_{2} \\ \bm{0}\end{array}\right\},\quad \mathrm{or}\quad \bm{K}_c(\bm{\xi}) \bm{u}_c(\bm{\xi}) = \bm{f}_c.
\end{eqnarray} 

{\color{black}The deterministic version of (\ref{eqn:sdcp}), i.e., corresponding to fixed realizations of $\bm\xi_1$ and $\bm\xi_2$, may be solved efficiently using the FETI method, described in Section \ref{sec:FETI}. This is done by eliminating the primal variables $\bm{u}_1$ and $\bm{u}_2$ form (\ref{eqn:sdcp}) and then solving the resulting system of equations for the dual variable $\bm{\lambda}$ using a projected conjugate gradient solver \citep{Farhat91}. }

We next present a stochastic discretization of (\ref{eqn:sdop}) and (\ref{eqn:sdcp}) using the standard PC expansions. Following that, in Section \ref{sec:SR}, we introduce the proposed stochastic discretization which allows us to solve (\ref{eqn:sdcp}) iteratively through approximations in terms of $\bm\xi_1$ or $\bm\xi_2$.

\section{Stochastic Discretization via Polynomial Chaos}
\label{sec:SG_Scheme}

Similar to spatial discretization, the first step in performing a stochastic discretization is to select an expansion basis. Several choices have been suggested for this purpose, e.g., polynomial chaos (PC) \citep{Ghanem91a,Xiu02}, piecewise polynomial \citep{Babuska07a}, and a multi-wavelet \citep{LeMaitre04} basis. In practice, the choice of expansion basis depends on the probability distribution of $\bm {\xi}$ and the regularity of the solution with respect to $\bm {\xi}$. In the present study, {\color{black}for the sake of simplicity,} we use the PC for the stochastic discretization.

Let $\psi_{{i}_k}(\xi_k)$, $k=1,\dots,d$, be a polynomial of degree ${i}_k  \in \mathbb{N}_0 := \mathbb{N}\cup \lbrace 0\rbrace$ orthogonal with respect to the measure $\rho(\xi_k)$, i.e.,
\begin{equation}
\label{eqn:1d_PCE}
\mathbb{E} [\psi_{i_k} \psi_{j_k}]  = \int \psi_{i_k}(\xi_k) \psi_{j_k}(\xi_k) \rho(\xi_k) \mathrm{d}\xi_k= \delta_{i_kj_k} \mathbb{E} [ \psi_{i_k}^2 ],
\end{equation}
where $\delta_{i_kj_k}$ is the Kronecker delta. The polynomial chaos of maximum (total) degree $p \in \mathbb{N}_0$ is then obtained by the tensorization of $\psi_{{i}_k}(\xi_k)$ corresponding to different directions $k$,
\begin{equation} 
\label{eqn:multivariate_polynomials}
\psi_{\bm{i}}(\bm{\xi}) = \prod_{k=1}^{d} \psi_{{i}_k}(\xi_k), \ \ \ \ \bm{i} \in \mathscr{I}_{d,p},
\end{equation}
where ${\bm{i}} = (i_1, \cdots , i_d) \in \mathscr{I}_{d,p}$ and the set of multi-indeces $\mathscr{I}_{d,p}$ is defined by
\begin{equation} 
\mathscr{I}_{d,p} = \lbrace \bm i = (i_1, \cdots , i_d)\in \mathbb{N}_0^d: \Vert\bm{i}\Vert_1 \leqslant p \rbrace.
\end{equation}

The cardinality of  $\mathscr{I}_{d,p}$, hence the number of PC basis functions of total order not larger than $p$ in dimension $d$, is given by 

\begin{equation} 
\label{eqn:P+1}
P=\vert \mathscr{I}_{d,p}\vert= \frac{(p+d)!}{p!d!}.
\end{equation}

{\color{black}Although the total order truncation is used here for simplicity, this gives only very {\it coarse} control of the size of the basis. It is therefore advisable to use some {\it finer} method to control the size of the basis, preferably in an adaptive or at least heuristic fashion, see, e.g., \cite{Bieri09a,Bieri09c}.}

Following the orthogonality of the polynomials $\psi_{{i}_k}(\xi_k)$ and given that the $\xi_k$ are independent, the PC basis functions $\psi_{\bm i}(\bm\xi)$ are also orthogonal, i.e., $\mathbb{E}[\psi_{\bm i}\psi_{\bm j}]=\delta_{\bm i,\bm j}\mathbb{E}[\psi_{\bm i}^2]$.

We define the finite dimensional space 
\begin{equation}
\label{eqn:wp}
\mathcal{W}_p = \mathrm{span} \lbrace \psi_{\bm{i}}(\bm{\xi}) : \bm{i} \in \mathscr{I}_{d,p}\rbrace \subset \mathcal{W}, \nonumber
\end{equation}
within which we seek an approximation to a finite-variance, stochastic function $u(\bm{\xi})\in\mathcal{W}$ by
\begin{equation} 
\label{eqn:RV_PCE}
u_p(\bm{\xi}) = \sum_{\bm{i} \in \mathscr{I}_{d,p}} u_{\bm{i}} \psi_{\bm{i}}(\bm{\xi}).
\end{equation}
The PC coefficients $u_{\bm{i}}$ are given by
\begin{equation} 
\label{eqn:RV_PCE_coefficients}
u_{\bm{i}}  = \frac{\mathbb{E} [u \psi_{\bm{i}} ]}{\mathbb{E} [ \psi_{\bm{i}}^2 ]}
\end{equation}
and the series in (\ref{eqn:RV_PCE}) converges in the mean-square sense as $p\to\infty$.

For the solutions to the original and coupled problems (\ref{eqn:sdop}) and (\ref{eqn:sdcp}), respectively, the direct projection (\ref{eqn:RV_PCE_coefficients}) may not be applied to compute the PC coefficients, unless (\ref{eqn:RV_PCE_coefficients}) is computed using random sampling or quadrature integration. Alternatively, as proposed in \cite{Ghanem91a}, the PC coefficients may be computed via Galerkin projection that we briefly describe next.

\subsection{Galerkin projection}
\label{sec:GP}
Let $\bm u_p(\bm{\xi}) = \sum_{\bm{i} \in \mathscr{I}_{d,p}} \bm u_{\bm{i}} \psi_{\bm{i}}(\bm{\xi})$ be the PC approximation to the exact solution $\bm u(\bm \xi)$ to (\ref{eqn:sdop}). The coefficients $\bm u_i$ can be computed from the Galerkin projection
\begin{eqnarray} 
\label{eqn:fdop}
&\mathbb{E}[\bm{v}^T (\bm{K} {\bm{u}_p}-\bm f)] = 0, \ \ \ \ \ \ \ \forall \bm{v} \in \mathbb{R}^{M} \otimes \mathcal{W}_p, \nonumber
\end{eqnarray} 
which can be equivalently written in the form of the coupled system of equations
\begin{equation} 
\label{eqn:soeop}
 \sum_{\bm{i} \in \mathscr{I}_{d,p}} \mathbb{E}[\psi_{\bm{i}} \psi_{\bm{j}} \bm{K}] \bm{u}_{\bm{i}} = \mathbb{E}[\psi_{\bm{j}} \bm{f}], \ \ \ \ \ \ \ \forall \bm{j} \in \mathscr{I}_{d,p}.
\end{equation} 

Following the same procedure, the Galerkin solution of the coupled problem (\ref{eqn:sdcp}) is computed from
\begin{equation} 
\label{eqn:soecp}
 \sum_{\bm{i} \in \mathscr{I}_{d,p}} \mathbb{E}[\psi_{\bm{i}} \psi_{\bm{j}} \bm{K}_c] {\bm{u}}_{c,\bm{i}} = \mathbb{E}[\psi_{\bm{j}} \bm{f}_c], \ \ \ \ \ \ \ \forall \bm{j} \in \mathscr{I}_{d,p},
\end{equation} 
where $\bm{K}_c$ and $\bm{f}_c$ are given in (\ref{eqn:sdcp}), and ${\bm{u}}_{c,\bm{i}} = \left[\begin{array}{ccc} {\bm{u}}_{1,\bm{i}}^T, {\bm{u}}_{2,\bm{i}}^T, {\bm{\lambda}}_{\bm{i}}^T \end{array}\right]^T$. Here, ${\bm{u}}_{1,\bm{i}}, {\bm{u}}_{2,\bm{i}},$ and ${\bm{\lambda}}_{\bm{i}}$ are the PC coefficients of the sub-domain solutions $\bm{u}_1$ and $\bm{u}_2$,  and the Lagrange multipliers $\bm{\lambda}$ in (\ref{eqn:sdcp}), respectively. 

\subsection{Need for alternative stochastic discretizations}
\label{sec:cost_SG_scheme}

We note that the size of the linear system of equations (\ref{eqn:soeop}) and (\ref{eqn:soecp}) depends on $P$ which is related to $d=d_1+d_2$ via (\ref{eqn:P+1}). In particular, the linear system (\ref{eqn:soeop}), corresponding to the PC approximation of the original problem (\ref{eqn:sdop}), is of size $M P \times MP$. On the other hand, assuming that ${\bm{\lambda}}_{\bm{i}}$ are given, ${\bm{u}}_{1,\bm{i}}$ and ${\bm{u}}_{2,\bm{i}}$ in the coupled formulation (\ref{eqn:soecp}) may be computed from linear systems of equations of sizes $M_1P \times M_1P$ and $M_2P \times M_2P$, respectively. While the latter formulation provides a reduction with respect to the size $M$ of the spatial discretization, it may be computationally prohibitive when $d$, hence $P$, is large. 

To alleviate this difficulty, we seek an alternative stochastic discretization that allows us to approximate the solution of the coupled problem (\ref{eqn:sdcp}) by solving problems in dimensions $d_1$ and $d_2$, instead of $d=d_1+d_2$, thereby partitioning the stochastic space in addition to the physical space. To this end, we propose the approximation of the solution to (\ref{eqn:sdcp}) via separated representations that we describe next. 

\section{Separated Representations}
\label{sec:SR} 

Given a target accuracy $\epsilon$, we consider the {\color{black}following} approximation of the solution to (\ref{eqn:sdcp}),
\begin{equation}
\label{eqn:rank-r-approx}
\left\{\begin{array}{c}\bm{u}_{1}(\bm{\xi}) \\\bm{u}_{2}(\bm{\xi}) \\ \bm{\lambda}(\bm{\xi})\end{array}\right\} = \sum_{l=1}^{r} \left\{\begin{array}{c} \bm{u}_{0,1}^{l}\\ \bm{u}_{0,2}^{l} \\   \bm{\lambda}_{0}^{l}\end{array}\right\} \phi_1^{l}(\bm{\xi}_1) \phi_2^{l}(\bm{\xi}_2) + \mathcal{O}(\epsilon), 
\end{equation}
which we refer to as the separated representation with {\it separation rank} $r$. {\color{black}Observe that the set of all rank $r$ approximations is not a linear subspace \cite{Hackbusch12}.} The deterministic vectors $\bm{u}_{0,1}^{l}\in\mathbb{R}^{M_1}$, $\bm{u}_{0,2}^{l}\in\mathbb{R}^{M_2}$, and $\bm{\lambda}_{0}^{l}\in\mathbb{R}^{M_I}$, the  stochastic functions $\phi_1^{l}(\bm\xi_1):\Omega\rightarrow\mathbb{R}$ and $\phi_2^{l}(\bm\xi_2):\Omega\rightarrow\mathbb{R}$, as well as the separation rank $r$ are not fixed \textit{a priori} and are sought for to achieve the accuracy $\epsilon$. 

Due to the separated construction of the representation (\ref{eqn:rank-r-approx}) with respect to variables along spatial as well as the stochastic directions $\bm\xi_1$ and $\bm\xi_2$, the unknowns in the RHS of (\ref{eqn:rank-r-approx}) may be computed from a sequence of linear approximations arising from an {\it alternating direction} optimization of a suitable cost function. One instance of such a multi-linear approach is described in Section \ref{sec:ARR}. 

The separation rank $r$ plays an important role in the construction of separated representations. In particular, smaller accuracies $\epsilon$ entail larger $r$ values. This, in turn, results in an increase in the computational complexity of the expansion, as we shall see in Section \ref{sec:ARR}. However, when the solution to (\ref{eqn:sdcp}) admits a low separation rank $r$, (\ref{eqn:rank-r-approx}) is essentially a reduced representation -- or compression -- of the solution. Therefore, we ideally desire to keep $r$ as small as possible. 

We note that the separated representation of the form (\ref{eqn:rank-r-approx}) may not be unique \citep{Beylkin02}. However, we seek one such representation that is within an accuracy $\epsilon$ of the solution to (\ref{eqn:sdcp}) which will be discussed in the sequel. 

\begin{remark1}
A more general separated approximation to the solution of (\ref{eqn:sdcp}) may be considered as
\begin{equation}
\label{eqn:rank-r-approx-general}
\left\{\begin{array}{c} \bm{u}_{1}(\bm{\xi}) \\\bm{u}_{2}(\bm{\xi}) \\ \bm{\lambda}(\bm{\xi})\end{array}\right\} = \sum_{l=1}^{r} \left\{\begin{array}{c} \bm{u}_{0,1}^{l}\\ \bm{u}_{0,2}^{l} \\   \bm{\lambda}_{0}^{l}\end{array}\right\} \prod_{i=1}^{d_1} {{\phi^{l}}_{1,i}}({\xi_{1,i}}) \prod_{j=1}^{d_2} {{\phi^{l}}_{2,j}}({\xi_{2,j}}) + \mathcal{O}(\epsilon).
\end{equation}
where ${{\phi^{l}}_{1,i}}$ and ${{\phi^{l}}_{2,j}}$ are univariate functions in $\xi_{1,i}$ and $\xi_{2,j}$, respectively. Using (\ref{eqn:rank-r-approx-general}), one may approximate $\bm{u}_{1}(\bm{\xi})$, $\bm{u}_{2}(\bm{\xi})$, and $\bm{\lambda}(\bm{\xi})$ via a series of one-dimensional approximations with as fast as a linear increase of the computation cost with respect to $d$. See \citep{Beylkin02, Doostan09} for further discussions on such growth in the case of single domain problems. This, however, is achieved by a potential increase in the separation rank $r$ as compared to the representation (\ref{eqn:rank-r-approx}). We leave further comparison of these different separated representations to a later study.
\end{remark1}
\subsection{Alternating Rayleigh-Ritz (ARR) algorithm}
\label{sec:ARR}

To derive the separated representation (\ref{eqn:rank-r-approx}), we use a Rayleigh-Ritz-type approach: we find the unknowns in the RHS of (\ref{eqn:rank-r-approx}) such that they correspond to a stationary point of the energy functional $\pi$ in (\ref{eqn:energy}). To do so, we adopt an alternating direction optimization approach, hence the name Alternating Rayleigh-Ritz (ARR), in order to compute (\ref{eqn:rank-r-approx}) using a sequence of linear, instead of nonlinear, problems {\color{black}of smaller size}. Notice that the stationary point of $\pi$ is a saddle point, i.e, it is a minimum of $\pi$ with respect to $u_1(\bm\xi)$ and $u_2(\bm\xi)$ and simultaneously a maximum of $\pi$ with respect to $\bm\lambda(\bm\xi)$.  

The ARR approach consists of a sequence of stationary point iterations that are meant to converge to the saddle point of $\pi$ as demonstrated in Proposition \ref{pro:saddle_conv}. {\color{black}Given a value for $r$ and approximations (guesses) $\lbrace \phi_1^{l} \rbrace_{l=1}^{r}$ and $\lbrace \phi_2^{l} \rbrace_{l=1}^{r}$, the ARR iteration starts} by freezing the variables along $\bm\xi_1$ and $\bm\xi_2$ directions, i.e., $\lbrace \phi_1^{l} \rbrace_{l=1}^{r}$ and $\lbrace \phi_2^{l} \rbrace_{l=1}^{r}$, and solving for the deterministic vectors $\lbrace \bm{u}_{0,1}^{l} \rbrace_{l=1}^{r}$, $\lbrace \bm{u}_{0,2}^{l} \rbrace_{l=1}^{r}$, and $\lbrace \bm{\lambda}_{0}^{l} \rbrace_{l=1}^{r}$ as the saddle point of $\pi$. We call this step a {\it deterministic update}. We then alternate to stochastic directions $\bm\xi_1$ and update $\lbrace \phi_1^{l} \rbrace_{l=1}^{r}$ by minimizing $\pi$, while $\lbrace \bm{u}_{0,1}^{l} \rbrace_{l=1}^{r}$, $\lbrace \bm{u}_{0,2}^{l} \rbrace_{l=1}^{r}$ and $\lbrace \phi_2^{l} \rbrace_{l=1}^{r}$ are fixed at their latest values. We dub this step a {\it stochastic update}.  Similarly, we alternate to update $\lbrace \phi_2^{l} \rbrace_{l=1}^{r}$ by minimizing $\pi$, while the rest of the variables are fixed at their current values. This completes one full sweep of the ARR algorithm. As we initialized $\lbrace \phi_1^{l} \rbrace_{l=1}^{r}$ and $\lbrace \phi_2^{l} \rbrace_{l=1}^{r}$ arbitrarily, we repeat the ARR sweeps until the value of $\pi$ does not change {\it much}. We then increase the separation rank $r$ and continue the ARR sweeps until $\pi$ does not change beyond some tolerance. 

As we ideally want the separation rank $r$ to be as small as possible, we may start the ARR updated with a small $r$, e.g., $r=1$ and increase $r$ until, for instance, the maximum norm of sub-domain residuals  
\begin{equation}
\label{eqn:residual_error_general}
{\epsilon_{res}^{r}} = \mathrm{max}\left( \frac{\| \bm{f}_1  - \bm{K}_1  \bm{u}_{1}^r\|_{L_2(\mathbb{R}^{M_1}) \otimes \mathcal{W}}}{\| \bm{f}_1\|_{L_2(\mathbb{R}^{M_1}) \otimes \mathcal{W}}}, \frac{\| \bm{f}_2  - \bm{K}_2 \bm{u}_{2}^r \|_{L_2(\mathbb{R}^{M_2}) \otimes \mathcal{W}}}{\| \bm{f}_2\|_{L_2(\mathbb{R}^{M_2}) \otimes \mathcal{W}}} \right)
\end{equation}
is below a prescribed accuracy $\epsilon$. Here, for $i=1,2$, $\bm{u}_{i}^r = \sum_{l=1}^{r} \bm{u}_{0,i}^{l} \phi_1^{l} \phi_2^{l}$, and $\| \cdot \|_{L_2(\mathbb{R}^{M}) \otimes \mathcal{W}}$ is defined by $\| \bm{u} \|^2_{L_2(\mathbb{R}^{M}) \otimes \mathcal{W}} = \mathbb{E}[\bm{u}^T \bm{u}]$.

In the following, we describe the ARR algorithm in more details.

\mbox{ }\\  
\noindent{\textbf{Deterministic updates for $\lbrace \bm{u}_{0,1}^{l} \rbrace_{l=1}^{r}$, $\lbrace \bm{u}_{0,2}^{l} \rbrace_{l=1}^{r}$, and $\lbrace \bm{\lambda}_{0}^{l} \rbrace_{l=1}^{r}$}.} Assuming an initialization of the separation rank $r$, and the stochastic functions $\lbrace \phi_1^{l} \rbrace_{l=1}^{r}$ and $\lbrace \phi_2^{l} \rbrace_{l=1}^{r}$, the first step of the ARR algorithm is to update deterministic vectors $\lbrace \bm{u}_{0,1}^{l} \rbrace_{l=1}^{r}$, $\lbrace \bm{u}_{0,2}^{l} \rbrace_{l=1}^{r}$, and $\lbrace \bm{\lambda}_{0}^{l} \rbrace_{l=1}^{r}$ while freezing all other variables. Plugging (\ref{eqn:rank-r-approx}) in (\ref{eqn:energy}) and enforcing the condition 
\begin{equation}
\label{eqn:mpe_{u}_{0,1}^{l}}
\delta \pi = 0,\qquad \forall\ \delta \bm{u}_{0,1}^{l}, \delta\bm{u}_{0,2}^{l}, \delta\bm{\lambda}_{0}^{l},\nonumber
\end{equation}
for $l=1,\dots,r$ and arbitrary (but consistent) variations $\delta \bm{u}_{0,1}^{l}$, $\delta\bm{u}_{0,2}^{l}$, and $\delta\bm{\lambda}_{0}^{l}$, we arrive at the saddle point problem
\begin{eqnarray} 
\label{eqn:lin_sys_0_r}
\left[\begin{array}{ccc} \bm{\hat{K}}_{1} & \bm{0} & -\bm{\hat{C}}_{1} \\ \bm{0} & \bm{\hat{K}}_{2} & \bm{\hat{C}}_{2} \\ -\bm{\hat{C}}^{T}_{1} & \bm{\hat{C}}^{T}_{2} & \bm{0}\end{array}\right]\left\{\begin{array}{c}\bm{\hat{u}}_{0,1} \\ \bm{\hat{u}}_{0,2}  \\ \bm{\hat{\lambda}}_0\end{array}\right\}= \left\{\begin{array}{c}{\bm{\hat{f}}_{1}} \\ {\bm{\hat{f}}_{2}}  \\ \bm{0}\end{array}\right\}.
\end{eqnarray} 

Each entry of the liner system (\ref{eqn:lin_sys_0_r}) has a block structure given by
\begin{eqnarray}
\label{eqn:lin_sys_0_r_2}
&& \bm{\hat{K}}_{i}(l,l')=\mathbb{E}\left[\phi_{1}^{l}\phi_{2}^{l} \bm{K}_{i}\phi_{1}^{l'}\phi_{2}^{l'}\right],\nonumber\\
&&\bm{\hat{C}}_{i}(l,l')=\mathbb{E}\left[ \phi_1^{l} \phi_2^{l}\phi_{1}^{l'}\phi_{2}^{l'}\right]\bm{C}_{i},\nonumber\\
&& \bm{\hat{u}}_{0,i}(l)=\bm{u}_{0,i}^l,\\
&& \bm{\hat{\lambda}}_0(l)={\bm{\lambda}_{0}^{l}},\nonumber\\
&&{\bm{\hat{f}}_{i}}(l)=\mathbb{E}\left[\phi_1^{l} \phi_2^{l}\right]\bm{f}_{i},\nonumber
\end{eqnarray}
for $l,l'=1,\dots,r$ and $i=1,2$. A detailed derivation of (\ref{eqn:lin_sys_0_r}) and (\ref{eqn:lin_sys_0_r_2}) can be found in \ref{apx:deterministic_updates}.

It should be noted that the form of (\ref{eqn:lin_sys_0_r}) is similar to that of the saddle point system associated with the deterministic version of the coupled problem (\ref{eqn:sdcp}), except that each entry in (\ref{eqn:lin_sys_0_r}) has a block structure. Such similarity, therefore, suggests an extension of the original FETI method \citep{Farhat91} for a parallel and partitioned solution of (\ref{eqn:lin_sys_0_r}). This will be discussed in more details in Section \ref{sec:FETI}.

\mbox{ }\\  
\noindent{\textbf{Stochastic updates for  $\lbrace \phi_1^{l} \rbrace_{l=1}^{r}$}.} We freeze $\lbrace \bm{u}_{0,1}^{l} \rbrace_{l=1}^{r}$, $\lbrace \bm{u}_{0,2}^{l} \rbrace_{l=1}^{r}$, and $\lbrace \bm{\lambda}_{0}^{l} \rbrace_{l=1}^{r}$ to their updated values from (\ref{eqn:lin_sys_0_r}). We also fix the variables $\lbrace \phi_2^{l} \rbrace_{l=1}^{r}$ and solve for $\lbrace \phi_1^{l} \rbrace_{l=1}^{r}$ by minimizing $\pi$, that is, for $l=1,\dots,r$, we require
\begin{equation}
\label{eqn:mpe_{phi}_{1}^{l}}
\delta \pi= 0,\qquad \forall\ \delta \phi_{1}^{l}(\bm\xi_1).
\end{equation}

This leads to the linear system
\begin{equation}
\label{eqn:lin_sys_1_rg}
\bm{A}_{1}(\bm\xi_1)\bm{\phi}_{1}(\bm\xi_1) =\bm{b}_1,
\end{equation}
where $\bm{A}_{1}\in \mathbb{R}^{r\times r}$, $\bm{\phi}_{1}\in \mathbb{R}^{r}$, and  $\bm{b}_1\in \mathbb{R}^{r}$ are given by
\begin{eqnarray}
\label{eqn:short_1_rg_update}
&&\bm{A}_1(l,l')={\bm{u}^{l}}^{T}_{0,1}\bm{K}_{1}(\bm\xi_1)\bm{u}_{0,1}^{l'}(\bm\xi_1)\mathbb{E}_{\bm\xi_2}\left[\phi_{2}^{l}
\phi_{2}^{l'}\right]+{\bm{u}^{l}}^{T}_{0,2}\mathbb{E}_{\bm\xi_2}\left[
\bm{K}_{2}\phi_{2}^{l}\phi_{2}^{l'}\right]
\bm{u}_{0,2}^{l'},\nonumber\\
&& \bm{\phi}_{1}(l) = \phi_1^l(\bm\xi_1), \\
&&\bm{b}_{1}(l) ={\bm{u}^{l}}^{T}_{0,1}\bm{f}_{1}\mathbb{E}_{\bm\xi_2}
\left[\phi_{2}^{l}\right] + {\bm{u}^{l}}^{T}_{0,2}\bm{f}_{2}\mathbb{E}_{\bm\xi_2}
\left[\phi_{2}^{l}\right],\nonumber
\end{eqnarray}
and $l,l'=1,\dots,r$. Here, $\mathbb{E}_{\bm\xi_2}[\cdot]$ is the mathematical expectation operator with respect to the probability density function of $\bm\xi_2$.  \ref{apx:stochastic_updates} provides further details on the derivation of (\ref{eqn:lin_sys_1_rg}) and (\ref{eqn:short_1_rg_update}).

\mbox{ }\\  
\noindent{\textbf{Stochastic updates for $\lbrace \phi_2^{l} \rbrace_{l=1}^{r}$}.} Following the same procedure for $\lbrace \phi_1^{l} \rbrace_{l=1}^{r}$, $\lbrace \phi_2^{l} \rbrace_{l=1}^{r}$ are updated by solving the linear system of size $r$
\begin{equation}
\label{eqn:lin_sys_2_rg}
\bm{A}_{2}(\bm\xi_2) \bm{\phi}_{2}(\bm\xi_2) =\bm{b}_2, 
\end{equation}
where, for $l,l'=1,\dots,r$,
\begin{eqnarray}
\label{eqn:short_2_update}
&&\bm{A}_2(l,l')={\bm{u}^{l}}^{T}_{0,1}
\mathbb{E}_{\bm\xi_1}\left[\bm{K}_{1}
\phi_{1}^{l}\phi_{1}^{l'}\right]\bm{u}_{0,1}^{l'}
+{\bm{u}^{l}}^{T}_{0,2}\bm{K}_{2}(\bm\xi_2)\bm{u}_{0,2}^{l'}\mathbb{E}_{\bm\xi_1}\left[
\phi_{1}^{l}\phi_{1}^{l'}\right],\nonumber\\
&& \bm{\phi}_{2}(l) = \phi_2^l(\bm\xi_2), \\
&&\bm{b}_{2}(l)={\bm{u}^{l}}^{T}_{0,1}\bm{f}_{1}\mathbb{E}_{\bm\xi_1}
\left[\phi_{1}^{l}\right] + {\bm{u}^{l}}^{T}_{0,2}\bm{f}_{2}\mathbb{E}_{\bm\xi_1}
\left[\phi_{1}^{l}\right],\nonumber
\end{eqnarray}
and $\mathbb{E}_{\bm\xi_1}[\cdot]$ is the mathematical expectation operator with respect to the probability density of $\bm\xi_1$.

We note that the update equations (\ref{eqn:lin_sys_1_rg}) and (\ref{eqn:lin_sys_2_rg}) only depend on $\bm\xi_1$ and $\bm\xi_2$, respectively, and not $\bm\xi = (\bm\xi_1,\bm\xi_2)$, thus yielding a partitioning of the stochastic space. Approximations to the solutions of the stochastic equations (\ref{eqn:lin_sys_1_rg}) and (\ref{eqn:lin_sys_2_rg}) may be obtained using, for instance, the stochastic Galerkin (SG) \citep{Ghanem91a} or the stochastic collocation \citep{Mathelin03, Xiu05a} technique. In the present study we employ SG, described in Section \ref{sec:SG_Scheme}, for this purpose. Specifically, for (\ref{eqn:lin_sys_1_rg}) or (\ref{eqn:lin_sys_2_rg}), we consider the PC expansions
\begin{equation}
\label{eqn:phi_1_2_PCE}
\bm{\phi}_{i}(\bm\xi_i) \approx \sum_{\bm{j} \in \mathscr{I}_{d_i,p}} {\bm\phi_{i,\bm{j}}} \psi_{\bm{j}}(\bm{\xi}_i),\qquad i=1,2,\nonumber
\end{equation}
where the coefficient vectors $\lbrace \bm\phi_{1,\bm{j}} \rbrace_{\bm{j} \in \mathscr{I}_{d_1,p}}$ and $\lbrace {\bm\phi_{2,\bm{j}}} \rbrace_{\bm{j} \in \mathscr{I}_{d_2,p}}$ are computed via Galerkin projection. Let $P_1$ and $P_2$ be the size of $\mathscr{I}_{d_1,p}$ and $\mathscr{I}_{d_2,p}$ evaluated from (\ref{eqn:P+1}) by setting $d=d_1$ and $d=d_2$, respectively. Then, the Galerkin systems to be solved for $\lbrace \bm\phi_{1,\bm{j}} \rbrace_{\bm{j} \in \mathscr{I}_{d_1,p}}$ and $\lbrace {\bm\phi_{2,\bm{j}}} \rbrace_{\bm{j} \in \mathscr{I}_{d_2,p}}$  are of size $rP_1 \times rP_1$ and $rP_2\times rP_2 $, respectively. Notice that, for moderate values of separation rank $r$, the sizes of these linear systems may be significantly smaller than those of (\ref{eqn:soeop}), particularly when $d_1$ and $d_2$ are large. Such a reduction is due to the use of separated representation (\ref{eqn:rank-r-approx}) along with the ARR algorithm for its construction. 

We next demonstrate that the sequence of separated approximations of the form (\ref{eqn:rank-r-approx}) generated by the ARR algorithm results in a non-increasing sequence of $\pi$ in $\bm u_1$ and $\bm u_2$, and a non-decreasing sequence of $\pi$ in $\bm \lambda$. In other words, the ARR algorithm iteratively improves the approximation to the saddle point problem (\ref{eqn:sdcp}) unless $\bm u_1$, $\bm u_2$, or $\bm \lambda$ does not change throughout the iterations. In the latter scenario, either the separation rank $r$ needs to be increased or the separated approximation has converged to the solution of (\ref{eqn:sdcp}). 

\begin{prop1}[Iterative improvement of separated representation] Let $(\bm u_1^{(k)},\bm u_2^{(k)},\bm\lambda^{(k)})$ and $(\bm u_1^{(k+1)},\bm u_2^{(k+1)},\bm\lambda^{(k+1)})$ be the rank $r$ separated representation of the solution $(\bm u_1,\bm u_2,\bm\lambda)$ after $k$ and $k+1$, respectively, ARR updates (\ref{eqn:lin_sys_0_r}), (\ref{eqn:lin_sys_1_rg}), and (\ref{eqn:lin_sys_2_rg}). Then,
\begin{equation}
\label{eqn:saddle_conv}
\pi(\bm u_1^{(k+1)},\bm u_2^{(k+1)},\bm\lambda^{(k)})\le \pi(\bm u_1^{(k)},\bm u_2^{(k)},\bm\lambda^{(k)}) \le \pi(\bm u_1^{(k)},\bm u_2^{(k)},\bm\lambda^{(k+1)}).
\end{equation}
\label{pro:saddle_conv}
\end{prop1}

\begin{proof}
We first note that the updates (\ref{eqn:lin_sys_1_rg}) and (\ref{eqn:lin_sys_2_rg}) along $\bm \xi_1$ and $\bm \xi_2$ directions, respectively, lead only to minimizations of $\pi$, see \ref{apx:stochastic_updates}. Additionally, $\hat{\bm u}_{0,1}$ and $\hat{\bm u}_{0,2}$ (or equivalently, $\lbrace \bm{u}_{0,1}^{l} \rbrace_{l=1}^{r}$ and $\lbrace \bm{u}_{0,2}^{l} \rbrace_{l=1}^{r}$) in (\ref{eqn:lin_sys_0_r}) are obtained by minimizing $\pi$, thus implying the left inequality in (\ref{eqn:saddle_conv}). Finally, the right inequality in (\ref{eqn:saddle_conv}) holds as $\pi$ may increase only through the updates $\hat{\bm\lambda}_0$ in (\ref{eqn:lin_sys_0_r}).
\end{proof}

The ARR algorithm is summarized in Algorithm \ref{alg:ARR}.

\begin{algorithm}[H]
\caption{Alternating Rayleigh-Ritz (ARR) Algorithm} 
\label{alg:ARR}
\begin{algorithmic}[1]
\STATE $\bullet$ \textbf{Input}: target accuracy $\epsilon$
\STATE $\bullet$ \textbf{Output}: separation rank $r$, $\lbrace \bm{u}_{0,1}^{l} \rbrace_{l=1}^{r}$, $\lbrace \bm{u}_{0,2}^{l} \rbrace_{l=1}^{r}$, $\lbrace \bm{\lambda}_{0}^{l} \rbrace_{l=1}^{r}$, $\lbrace \phi_1^{l} \rbrace_{l=1}^{r}$, and $\lbrace \phi_2^{l} \rbrace_{l=1}^{r}$
\STATE $\bullet$ Set $r = 1$; (randomly) initialize $\phi_1^{1}$ and $\phi_2^{1} $
\WHILE{${\epsilon_{res}^{r}} > \epsilon $}
\WHILE{$\pi$ {\color{black}changes} more than some tolerance $tol_r$}
\STATE $\bullet$ Update $\lbrace \bm{u}_{0,1}^{l} \rbrace_{l=1}^{r}$, $\lbrace \bm{u}_{0,2}^{l} \rbrace_{l=1}^{r}$, and $\lbrace \bm{\lambda}_{0}^{l} \rbrace_{l=1}^{r}$ by solving (\ref{eqn:lin_sys_0_r})
\STATE $\bullet$ Update $\lbrace \phi_1^{l} \rbrace_{l=1}^{r}$ by solving (\ref{eqn:lin_sys_1_rg})
\STATE $\bullet$ Update $\lbrace \phi_2^{l} \rbrace_{l=1}^{r}$ by solving (\ref{eqn:lin_sys_2_rg})
\ENDWHILE
\STATE $\bullet$ Compute the residual error ${\epsilon_{res}^{r}}$
\STATE $\bullet$ Set ${r=r+1}$; (randomly) initialize $\phi_1^{r}$ and $\phi_2^{r} $
\ENDWHILE
\end{algorithmic}
\end{algorithm}
\begin{remark1}
It is worthwhile noting that the ARR algorithm presented here may be interpreted as a variation of the block Gauss-Seidel or multiplicative Schwarz \cite{Smith04} methods.
\end{remark1}

To complete our discussion on the ARR algorithm we next present an extension of the classical FETI algorithm to solve the saddle point problem (\ref{eqn:lin_sys_0_r}).

\subsection{FETI solution of saddle point problem (\ref{eqn:lin_sys_0_r})}
\label{sec:FETI}

We are interested in parallel solution of (\ref{eqn:lin_sys_0_r}) such that minimal exchange of information between the sub-domain solvers is required. {\color{black}This may be achieved by, for instance, an extension of the classical FETI method for the particular case of problem (\ref{eqn:lin_sys_0_r}). This is the approach we take in the present study.} FETI, one of the widely used DD approaches, was first introduced by Farhat and Roux \citep{Farhat91} for the parallel FE solution of second order, self-adjoint elliptic equations. It has also been successfully extended to many other problems; multiple right-hand sides \citep{Franca98}, transient problems \citep{Farhat95}, plate bending problems \citep{Farhat98}, and non-linear problems \citep{Geradin97}. {\color{black}We stress that, for the sake of simplicity,} we use the original FETI \citep{Farhat91} approach, but with some extensions one may also utilize more recent versions of FETI such as FETI-DP \citep{Farhat01}.

The FETI method is a non-overlapping DD scheme based on partitioning the computational domain into a number of sub-domains and enforcing the continuity of the solution across the sub-domain interfaces via Lagrange multipliers. In practice, after partitioning the original domain $\mathcal{D}$, there might be some sub-domains with insufficient Dirichlet boundary conditions resulting in positive semi-definite sub-domain stiffness matrices. These are also known as floating sub-domains. For a partition consisting of two sub-domains $\mathcal{D}_1$ and $\mathcal{D}_2$, we assume that $\mathcal{D}_2$ has no Dirichlet boundary condition and is a {\it floating} sub-domain with symmetric positive semi-definite stiffness matrix $\bm{K}_2(\bm{\xi}_{2})$, while $\mathcal{D}_1$ is non-floating with symmetric positive definite  stiffness matrix $\bm{K}_1(\bm{\xi}_{1})$. Therefore, the second equation in (\ref{eqn:lin_sys_0_r}), i.e.,
\begin{eqnarray}
\label{eqn:feti_floating_u}
&&\hat{\bm{K}}_2 \hat{\bm{u}}_{0,2} =\hat{\bm{f}}_2 - \hat{\bm{C}}_2 \hat{\bm{\lambda}}_{0},
\end{eqnarray}
is solvable if and only if, \citep{Farhat91}, 
\begin{equation}
\label{eqn:solvability_condition}
\hat{\bm{R}}_2 \left( \hat{\bm{f}}_2 - \hat{\bm{C}}_2 \hat{\bm{\lambda}}_{0} \right)=\bm{0},
\end{equation}
where $\hat{\bm{R}}_2$ spans the null space of $\hat{\bm{K}}_2$, i.e., $\mathrm{range} (\hat{\bm{R}}_2) = \mathrm{ker} (\hat{\bm{K}}_2)$. In Section \ref{sec:RBMs} we will discuss how $\hat{\bm{R}}_2$ is set in our formulation. The solution of (\ref{eqn:lin_sys_0_r}) can then be written as
\begin{eqnarray}
\label{eqn:feti_solution_1}
&& \left( \hat{\bm{C}}_1^T \hat{\bm{K}}_1^{-1} \hat{\bm{C}}_1 + \hat{\bm{C}}_2^T \hat{\bm{K}}_2^{+} \hat{\bm{C}}_2 \right)\hat{\bm{\lambda}}_{0} =  -\hat{\bm{C}}_1^T \hat{\bm{K}}_1^{-1} \hat{\bm{f}}_1 + \hat{\bm{C}}_2^T \left( \hat{\bm{K}}_2^{+} \hat{\bm{f}}_2 + \hat{\bm{R}}_2 \hat{\bm{\alpha}}_{0} \right),\nonumber\\
&& \hat{\bm{u}}_{0,1} = \hat{\bm{K}}_1^{-1} \left( \hat{\bm{f}}_1 - \hat{\bm{C}}_1 \hat{\bm{\lambda}}_{0} \right),\nonumber\\
&& \hat{\bm{u}}_{0,2} = \hat{\bm{K}}_2^{+} \left( \hat{\bm{f}}_2 - \hat{\bm{C}}_2 \hat{\bm{\lambda}}_{0} \right) + \hat{\bm{R}}_2 \hat{\bm{\alpha}}_{0},
\end{eqnarray}
where $\hat{\bm{K}}_1^{-1}$ and $\hat{\bm{K}}_2^{+}$ are the inverse and pseudo-inverse of $\hat{\bm{K}}_1$ and $\hat{\bm{K}}_2$, respectively. Here, $\hat{\bm{\alpha}}_{0}$ is an unknown vector of coefficients. As it can be seen, there are four unknowns in (\ref{eqn:feti_solution_1}) while only three equations are available. The orthogonality condition in (\ref{eqn:solvability_condition}) is the fourth equation which makes (\ref{eqn:feti_solution_1}) determined. We can write the FETI {\it interface problem} by combining (\ref{eqn:solvability_condition}) and the first equation in (\ref{eqn:feti_solution_1}) as
\begin{equation}
\label{eqn:feti_interface}
\left[\begin{array}{cc} \hat{\bm{F}}_I & -\hat{\bm{R}}_2^I \\ -\hat{\bm{R}}_2^{I^T} & \bm{0} \end{array}\right] \left\{\begin{array}{c} \hat{\bm{\lambda}} \\ \hat{\bm{\alpha}}_{0} \end{array}\right\} = \left \{\begin{array}{c} \hat{\bm{d}} \\ -\hat{\bm{e}} \end{array}\right\},
\end{equation}
where
\begin{eqnarray}
&&\hat{\bm{F}}_I =\hat{\bm{C}}_1^T \hat{\bm{K}}_1^{-1} \hat{\bm{C}}_1 + \hat{\bm{C}}_2^T \hat{\bm{K}}_2^{+} \hat{\bm{C}}_2, \nonumber\\
&&\hat{\bm{R}}_2^I = \hat{\bm{C}}_2^T \hat{\bm{R}}_2, \nonumber\\
&&\hat{\bm{d}} = \hat{\bm{C}}_2^T  \hat{\bm{K}}_2^{+} \hat{\bm{f}}_2-\hat{\bm{C}}_1^T \hat{\bm{K}}_1^{-1} \hat{\bm{f}}_1, \nonumber\\
&&\hat{\bm{e}} = \hat{\bm{R}}_2^{T} \hat{\bm{f}}_2.
\end{eqnarray}

The solution to (\ref{eqn:lin_sys_0_r}) can then be obtained by first solving the interface problem (\ref{eqn:feti_interface}) for $\hat{\bm{\lambda}}$ and $\hat{\bm{\alpha}}_{0}$ and then by substituting these into the  last two equations in (\ref{eqn:feti_solution_1}) to compute $\hat{\bm{u}}_{0,1}$ and $\hat{\bm{u}}_{0,2}$. In order to avoid an explicit assembly of $\hat{\bm{F}}_I$ in (\ref{eqn:feti_interface}), the preconditioned conjugate projected gradient (PCPG) algorithm of \citep{Farhat91}, which utilizes matrix-vector products, may be extended to compute $\hat{\bm{\lambda}}$ in (\ref{eqn:feti_interface}). Algorithm \ref{Algorithm:PCPG} summarizes steps of PCPG solver, in which $\epsilon_{PCPG}$ is a stopping criterion, $\hat{\bm{P}} = \bm{I} - \hat{\bm{R}}_2^{I} (\hat{\bm{R}}_2^{I^T} \hat{\bm{R}}_2^{I})^{-1} \hat{\bm{R}}_2^{I^T}$ is an orthogonal projection onto the null space of $\hat{\bm{R}}_2^{I}$, $\bm{I}$ is the identity matrix of size $rM_I \times rM_I$, and $\bar{\bm{F}}_I^{-1}$ is a suitable preconditioner.
For a detailed description of the PCPG solver and different choices of preconditioners we refer the interested reader to \citep{Farhat99,Klawonn01,Toselli99a,Rapetti01,Galvis10,Charmpis02}. In the numerical result of Section \ref{sec:elasticity}, we have adopted the preconditioner $\bar{\bm{F}}_I^{-1}$ proposed in \cite{Klawonn01}. An {\it optimal} selection or design of $\bar{\bm{F}}_I^{-1}$, however, requires further study. 

After computing $\hat{\bm{\lambda}}$, we obtain $\hat{\bm{\alpha}}_{0}$ from
\begin{eqnarray}
\label{eqn:feti_alpha}
&&\hat{\bm{\alpha}}_{0} =(\hat{\bm{R}}_2^{I^T} \hat{\bm{R}}_2^{I})^{-1} \hat{\bm{R}}_2^{I^T} ( \hat{\bm{F}}_I \hat{\bm{\lambda}} - \hat{\bm{d}}).\nonumber
\end{eqnarray}
\begin{algorithm}[H]
\caption{The FETI PCPG Algorithm}
\begin{algorithmic}[1]
\label{Algorithm:PCPG}

\STATE $\bullet$ Initialize $\hat{\bm{\lambda}}_0 = \hat{\bm{R}}_2^{I} (\hat{\bm{R}}_2^{I^T} \hat{\bm{R}}_2^{I})^{-1} \hat{\bm{R}}_2^{T} \hat{\bm{f}}_2$, $\bm{w}_0 = \hat{\bm{P}} \hat{\bm{d}} - \hat{\bm{P}} \hat{\bm{F}}_I \hat{\bm{\lambda}}_0$, $k=1$

\WHILE{$||\bm{w}_{k-1}||_{L_2} / ||\hat{\bm{d}}||_{L_2} \geqslant \epsilon_{PCPG}$ }

\STATE $\bullet$ $\bm{z}_{k-1} = \bar{\bm{F}}_I^{-1} (\hat{\bm{P}}^T \bm{w}_{k-1})$ 
\STATE $\bullet$ $\bm{y}_{k-1} = \hat{\bm{P}} \bm{z}_{k-1}$ 
\STATE $\bullet$ ${s}_{k} = \bm{y}_{k-1}^T  \bm{w}_{k-1}  /  \bm{y}_{k-2}^T  \bm{w}_{k-2} $ \ \ \  \ \ \ \ $({s}_{1} = 0)$
\STATE $\bullet$ $\bm{p}_{k} = \bm{y}_{k-1} + {s}_{k}\bm{p}_{k-1}$ \ \ \ \ \ \ \ \ \ \ \ \ \ \ \ $(\bm{p}_{1} = \bm{y}_{0})$
\STATE $\bullet$ ${\gamma}_{k} = \bm{y}_{k-1}^T  \bm{w}_{k-1} / \bm{p}_{k}^T \hat{\bm{F}}_I \bm{p}_{k} $
\STATE $\bullet$ $\hat{\bm{\lambda}}_{k} = \hat{\bm{\lambda}}_{k-1} + {\gamma}_{k} \bm{p}_{k}$
\STATE $\bullet$ ${\bm{w}}_{k} = {\bm{w}}_{k-1} - {\gamma}_{k} \hat{\bm{P}} \hat{\bm{F}}_I \bm{p}_{k}$
\ENDWHILE
\end{algorithmic}
\end{algorithm}

Next, we propose a method to set $\hat{\bm{R}}_2$ which does not require an explicit computation of the null space of $\hat{\bm{K}}_2$.

\subsubsection{Setting $\hat{\bm{R}}_2$}
\label{sec:RBMs}

Since $\hat{\bm{K}}_2$ is not assembled in practice, we cannot obtain $\hat{\bm{R}}_2$ via a direct null space computation of $\hat{\bm{K}}_2$. For a deterministic problem, rigid body modes of the physical domain span the null space of the stiffness matrix \citep{Farhat91}. For three-dimensional problems, the maximum number of rigid body modes is six, while for  two-dimensional deterministic problems this number reduces to three. Let ${\bm{R}}_2$ denote the matrix associated with the rigid body modes of the sub-domain $\mathcal{D}_2$ computed based on the geometry of $\mathcal{D}_2$, see, e.g. \citep{Park98}. Then $\hat{\bm{R}}_2$ may be set as follows.

\begin{prop1}
$\hat{\bm{R}}_2$ is an $r \times r$ block diagonal matrix whose $(l,l)$ block, $l=1,\cdots,r$, is given by
\begin{equation}
\label{eqn:rbm}
\hat{\bm{R}}_2(l,l) = {\bm{R}}_2.
\end{equation}
\end{prop1}

\begin{proof}
Since columns of ${\bm{R}}_2$ are linearly independent, columns of $\hat{\bm{R}}_2$ are also linearly independent as $\hat{\bm{R}}_2$ is block diagonal. On the other hand, as the columns of ${\bm{R}}_2$ form a basis for the null space of $\bm{K}_2(\bm\xi)$, we have $\hat{\bm{K}}_2(l,l)\bm{R}_2 = \bm 0$, $l=1,\cdots,r$, where $\hat{\bm{K}}_2(l,l)$ is given in (\ref{eqn:lin_sys_0_r_2}). Consequently $\hat{\bm{K}}_2 \hat{\bm{R}}_2 = \bm{0}$, which, together with the linear independency of the columns of $\hat{\bm{R}}_2$, implies $\mathrm{range} (\hat{\bm{R}}_2) = \mathrm{ker} (\hat{\bm{K}}_2)$.
\end{proof}

\subsection{Response statistics}
\label{sec:Response_Statistics}
In this section, the computation of response statistics based on the separated representation (\ref{eqn:rank-r-approx}) is presented. Let us denote the mean and variance of the separated representation $\bm{u}_{i}^r = \sum_{l=1}^{r} \bm{u}_{0,i}^{l} \phi_1^{l} \phi_2^{l}$, $i=1,2$, by $\mathbb{E}[\bm{u}_{i}^r]$ and $\mathrm{var}[\bm{u}_{i}^r]$, respectively. For example, we can approximate the point-wise mean and the second moment of the sub-domain solutions $\bm{u}_{i}$ by
\begin{eqnarray} 
\label{eqn:residual_error}
&&\mathbb{E}[\bm{u}_{i}^r] = \sum_{l=1}^{r} \bm{u}_{0,i}^{l} \mathbb{E}_{\bm{\xi}_1}[\phi_1^{l}] \mathbb{E}_{\bm{\xi}_2}[\phi_2^{l}], \nonumber \\
&&\mathbb{E}[\bm{u}_{i}^r\circ\bm{u}_{i}^r] = \sum_{l=1}^{r} \sum_{l'=1}^{r} ( \bm{u}_{0,i}^{l} \circ \bm{u}_{0,i}^{l'} ) \mathbb{E}_{\bm{\xi}_1}[\phi_1^{l}\phi_1^{l'}] \mathbb{E}_{\bm{\xi}_2}[\phi_2^{l} \phi_2^{l'}],
\end{eqnarray}
respectively, where $\circ$ denotes the Hadamard product of two vectors. The variance of $\bm{u}_{i}^r$ can then be obtained from $\mathrm{var}[\bm{u}_{i}^r]  = \mathbb{E}[\bm{u}_{i}^r\circ\bm{u}_{i}^r] - \mathbb{E}[\bm{u}_{i}^r]\circ \mathbb{E}[\bm{u}_{i}^r]$.

Alternatively, one may generate Monte Carlo estimates of the statistics of $\bm{u}_{i}^r$ by directly sampling the representation (\ref{eqn:rank-r-approx}). We use the latter approach to compute the probability density function of the solution of interest in our numerical experiments. 

\section{Numerical Examples}
\label{sec:Numerical_Examples}

In this section, we present two numerical examples for the verification of the proposed stochastic coupling framework. In the first example, we consider a linear elliptic PDE defined on a two-dimensional L-shaped (spatial) domain with random diffusion coefficient. The second example deals with the problem of deformation of a linear elastic cantilever beam with uncertain Young's modulus. In both cases, the uncertainly is represented by random fields taking statistically independent values on two non-overlapping sub-domains. To study the convergence of the separated representation on the entire domain, the mean and standard deviation {\color{black}(std)} error measures 
\begin{equation} 
\label{eqn:residual_error}
{\epsilon_{\mu}^{r}}  = \frac{\big\Vert \mathbb{E}[\bm{u}_{ref}] -\mathbb{E}[\bm{u}^{r}] \big\Vert_{L_2(\mathbb{R}^M)}}{\big\Vert \mathbb{E}[\bm{u}_{ref}]\big\Vert_{L_2(\mathbb{R}^M)}}\;\; \mathrm{and}\;\;
{\epsilon_{\sigma}^{r}} = \frac{\big\Vert \mathrm{std}[\bm{u}_{ref}] -\mathrm{std}[\bm{u}^{r}] \big\Vert_{L_2(\mathbb{R}^M)}}{\big\Vert \mathrm{std}[\bm{u}_{ref}] \big\Vert_{L_2(\mathbb{R}^M)}}
\end{equation}
are used. Here $\bm{u}^{r}$ is the rank $r$ separated representation of the solution $\bm u$ to the coupled problem (\ref{eqn:sdcp})  and $\bm{u}_{ref}$ is the solution to the original problem (\ref{eqn:sdop}). The reference solution $\bm{u}_{ref}$ is the PC expansion of $\bm u$ constructed with respect to all random inputs $\bm \xi$. The PC coefficients were computed via least-squares regression using sufficiently large numbers of random realizations of $\bm u$. 

\subsection{Example I: 2D elliptic PDE with random diffusion coefficient}
\label{sec:elliptic_SPDE}
\begin{figure}[h]
  \begin{center}
	\includegraphics [width=8cm]{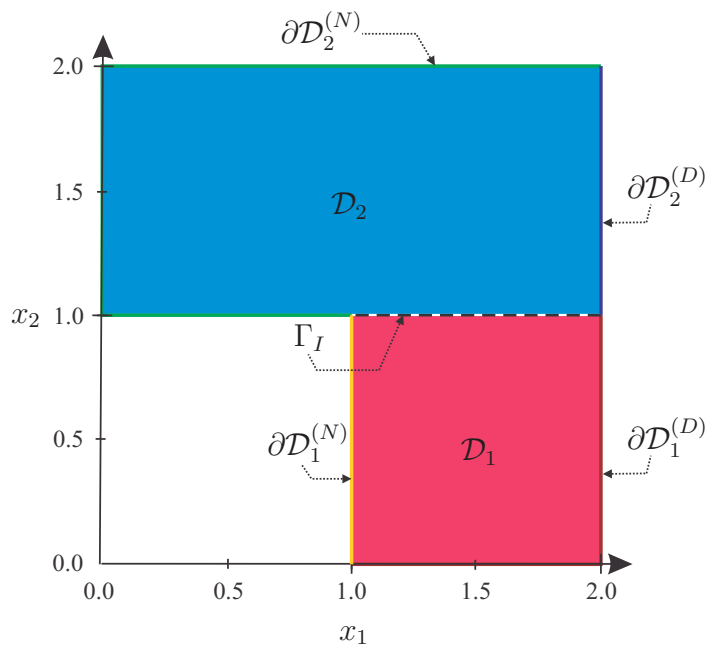}
    \put(-122,146){$\mathcal{D}_2$}
    \put(-73,53){$\mathcal{D}_1$}
    \put(-10,57){$\partial \mathcal{D}_1^{(D)}$}
    \put(-10,152){$\partial \mathcal{D}_2^{(D)}$}
    \put(-146,55){$\partial \mathcal{D}_1^{(N)}$}
    \put(-141,211){$\partial \mathcal{D}_2^{(N)}$}
    \put(-136,96){$\Gamma_I$}
    \put(-119,-15){$x_1$}
    \put(-243,105){$x_2$}
    \caption{Geometry of the 2D elliptic problem over an L-shaped domain ($\mathcal{D} = \mathcal{D}_1 \cup \mathcal{D}_2$ and $\mathcal{D}_1 \cap \mathcal{D}_2=\emptyset$).}
    \label{fig:L_shaped}
   \end{center}
\end{figure}

Consider the following elliptic PDE with random diffusion coefficient $\kappa$, 
\begin{eqnarray} 
\label{eqn:elliptic_SPDE}
&&-\nabla\cdot \left(\kappa(\bm{x},\bm{\xi}) \nabla u(\bm{x},\bm{\xi})\right) = f(\bm x)\qquad \bm{x} \in \mathcal{D}_, \nonumber \\
&&u(\bm{x},\bm{\xi})=0 \qquad  \qquad \qquad \qquad \qquad \ \ \bm{x} \in \partial \mathcal{D}^{(D)}, \nonumber \\
&&\nabla u(\bm{x},\bm{\xi}) \cdot \bm{n} =0  \qquad \qquad \qquad \qquad \ \bm{x} \in \partial \mathcal{D}^{(N)},
\end{eqnarray}
defined over an L-shaped domain $\mathcal{D}$ depicted in Fig. \ref{fig:L_shaped}. Dirichlet and Neumann boundary conditions are denoted by $\partial \mathcal{D}^{(D)}$ and $\partial \mathcal{D}^{(N)}$, respectively, and $\bm{n}$ is the unit normal vector to $\partial \mathcal{D}^{(N)}$. 

We assume $\mathcal{D}$ is composed of two non-overlapping sub-domains $\mathcal{D}_1$ and $\mathcal{D}_2$, i.e., $\mathcal{D}=\mathcal{D}_1\cup\mathcal{D}_2$ and $\mathcal{D}_1\cap\mathcal{D}_2=\emptyset$ on which $\kappa$ and $f$ take different values. Specifically, we set
\begin{equation}
\label{eqn:kappa_main}
\log\left(\kappa(\bm x,\bm \xi)-\kappa_0\right) = \left\{\begin{array}{c}G_1(\bm x,\bm \xi_1)\quad \bm x\in\mathcal{D}_1 \\ G_2(\bm x,\bm \xi_2)\quad \bm x\in\mathcal{D}_2 \end{array}\right.
\end{equation}
and
\begin{equation}
f(\bm x) = \left\{\begin{array}{c} 10 \quad \bm x\in\mathcal{D}_1 \\ 0\ \ \quad \bm x\in\mathcal{D}_2 \end{array}\right. ,\nonumber
\end{equation}
where $G_1$ and $G_2$ are statistically independent Gaussian random fields, and $\kappa_0$ is a small constant to ensure $\kappa$ is bounded away from zero. Therefore, $\kappa$ in (\ref{eqn:kappa_main}) is a (shifted) lognormal random field on sub-domains $\mathcal{D}_1$ and $\mathcal{D}_2$. We represent $G_1$ and $G_2$ by the truncated Karhunen-Lo\`eve (KL) expansions
\begin{equation}
\label{eqn:KL_Gaussian}
G_i(\bm{x},\bm{\xi}_i) = \bar{G}_i + \sum_{j=1}^{d_i} \sqrt{\tau_{i,j}} g_{i,j}(\bm{x}) {\xi}_{i,j},\qquad i=1,2,\nonumber
\end{equation}
where $\bar{G}_i $ is the mean of $G_i$, and $\left\lbrace \xi_{1,j} \right\rbrace_{j=1}^{d_1}$ and $\left\lbrace \xi_{2,j} \right\rbrace_{j=1}^{d_2}$ are i.i.d. standard Gaussian random variables. Additionally, $\left\lbrace \tau_{i,j} \right\rbrace_{j=1}^{d_i}$ and $\left\lbrace g_{i,j}(\bm{x}) \right\rbrace_{j=1}^{d_i}$ are, respectively, $d_i$ largest eigenvalues and the corresponding eigenfunctions of the Gaussian covariance kernel
\begin{equation}
\label{eqn:covariance_kernel_k}
C_i(\bm{x}_1,\bm{x}_2) = \sigma_{i}^2\mathrm{exp} {\left( -\frac{\| \bm{x}_1 - \bm{x}_2 \|_2 ^2}{l_{i}^2}\right)}, \qquad \bm{x}_1, \bm{x}_2 \in \mathcal{D}_i.
\end{equation}

Here, $\sigma_i$ controls the variability of $G_i$ and $l_i$ is the correlation length of $G_i$. To exactly compute the expectations involving $\bm K_1(\bm\xi_1)$ and $\bm K_2(\bm\xi_2)$ in (\ref{eqn:lin_sys_0_r_2}), (\ref{eqn:short_1_rg_update}), and (\ref{eqn:short_2_update}), we decompose $\kappa(\bm{x},\bm{\xi})$, separately on each sub-domain, into Hermite polynomial chaos of order $2p_i$, i.e., $\kappa(\bm{x},\bm{\xi}_i) = \sum_{\bm{j} \in \mathscr{I}_{d_i,2p_i}} \kappa_{\bm{j},i}(\bm{x}) \psi_{\bm{j}}(\bm{\xi}_i)$, $\bm x\in\mathcal{D}_i$. Here, $p_i$ is the order of the PC expansion of the stochastic functions $\{\phi_{i}^l(\bm \xi_i)\}$ in (\ref{eqn:phi_1_2_PCE}). For the case of lognormal $\kappa(\bm x,\bm \xi_i)$, the expansion coefficients $\{\kappa_{\bm{j},i}(\bm{x})\}_{\bm{j} \in \mathscr{I}_{d_i,2p_i}}$ are available analytically, see \citep{Ullmann08} or \ref{apx:Wiener_chaos_coefficients}. 

Table \ref{table:parameters_ex1} summarizes the parameter values involved in the representation of $\kappa$. 
 
\begin{table}[htb]
\caption{Assumed parameters for the PC representation of $\kappa$ in (\ref{eqn:kappa_main}).} 
\centering
\begin{tabular}{ c c c c c c c c c c}   
\hline
$d_1$ & $d_2$ & $p_1$ & $p_2$ & $l_{c,1}$ & $l_{c,2}$ & $\bar{G}$ & $\sigma_1$ & $\sigma_2$ & $\kappa_0$\\ 
\hline
4 & 6 & 3 & 3 & 2/3 & 1/3 & 1 & $1/2$  & $1/2$ & 0.28\\ 
\hline 
\end{tabular} 
\label{table:parameters_ex1} 
\end{table}

\begin{figure}[htb] 
    \centering
    \begin{tabular}{cc}
            \hspace{-0.5cm}    
      \includegraphics[width=2.8in]{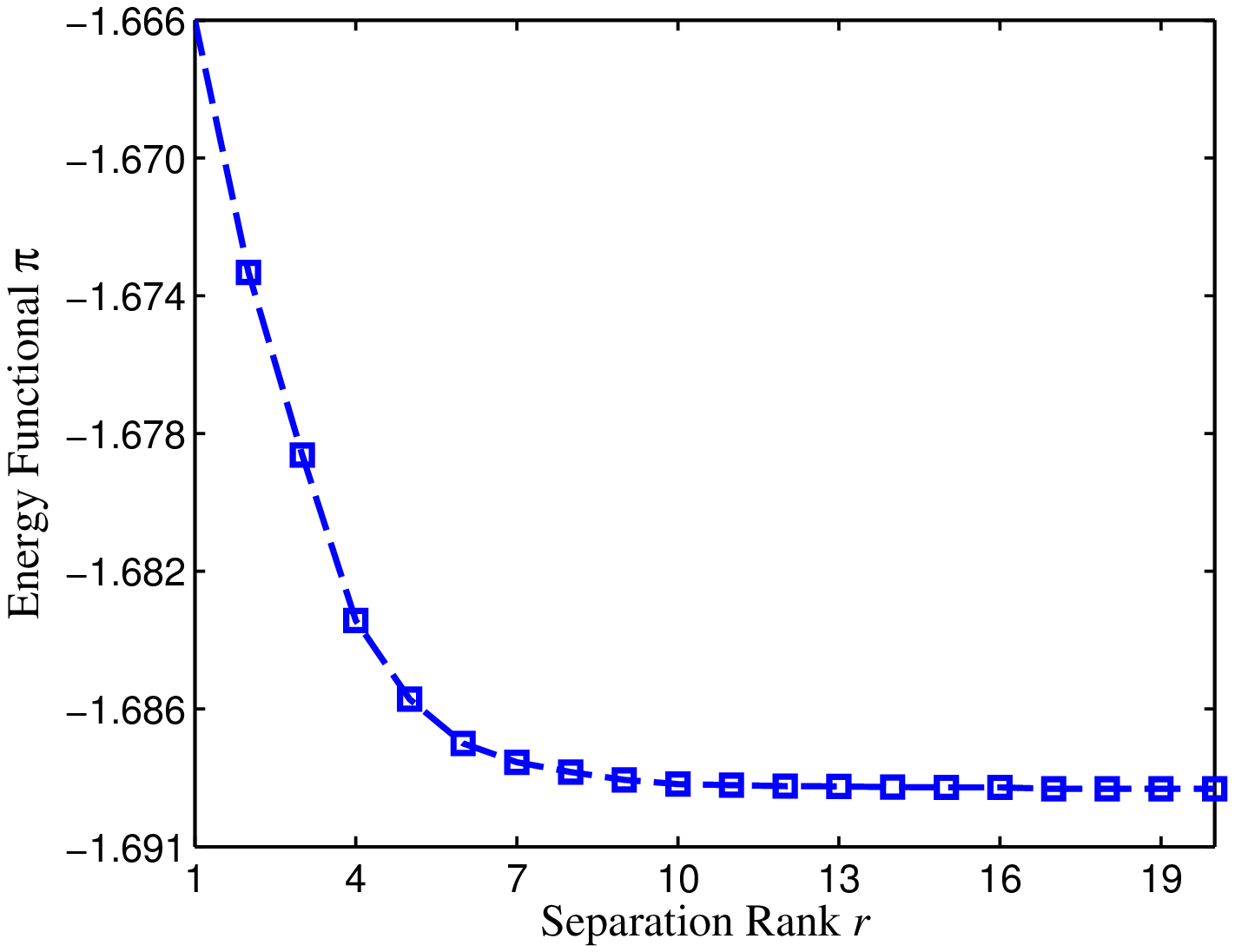}  
      &
      \includegraphics[width=2.8in]{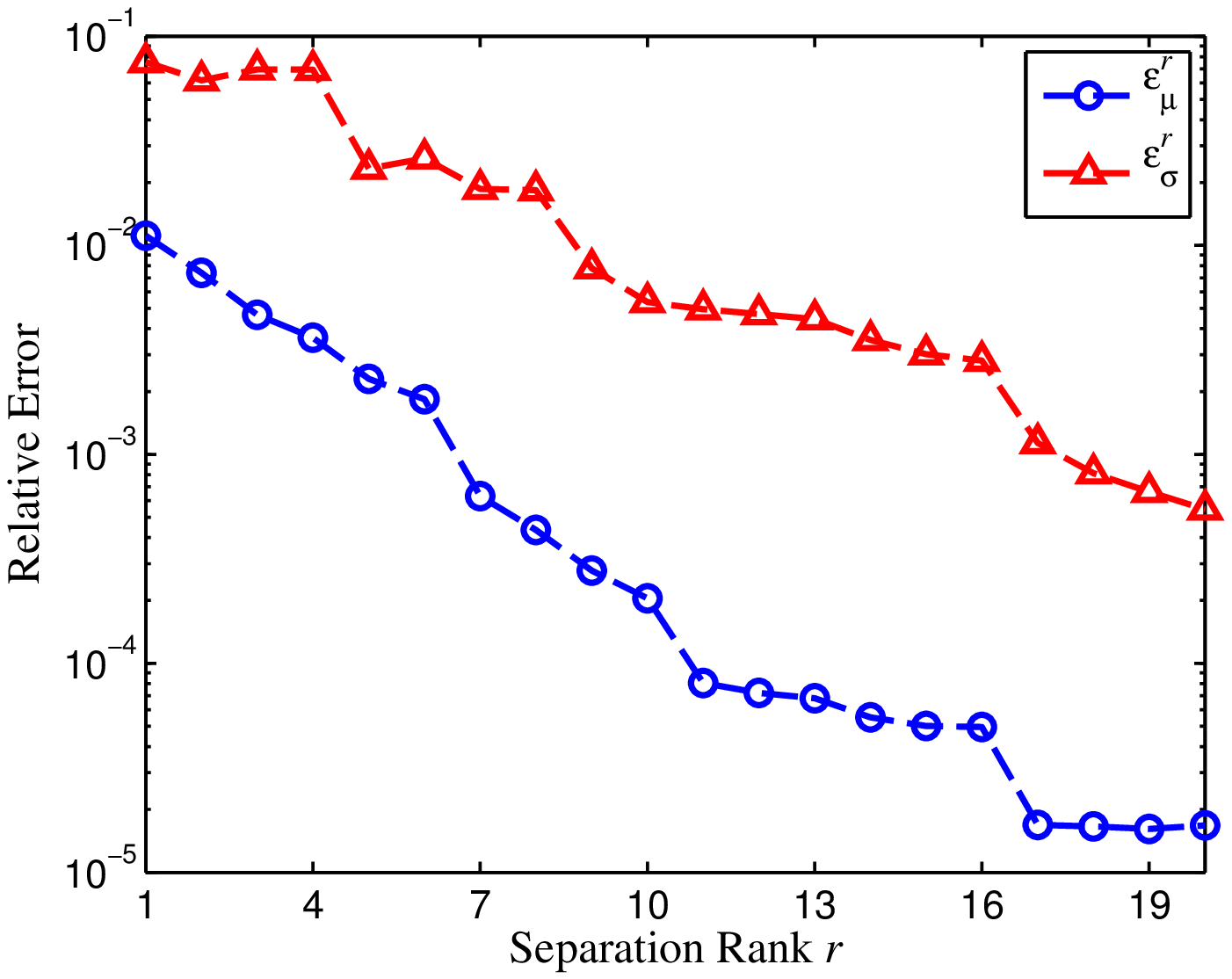} 
      \\
      (a) & (b)  
     \end{tabular}
      \caption{Energy functional $\pi$ (a) and relative errors in mean and standard deviation (b) as a function of the separation rank $r$ for the L-shaped problem. The errors are evaluated from (\ref{eqn:residual_error}).}            
\label{fig:L_eng_error}       
\end{figure}

The FE discretization of (\ref{eqn:elliptic_SPDE}) on each sub-domain is done using the FEniCS package \citep{Logg12}. In particular, we use linear triangle elements with a uniform mesh sizes $h_1 = h_2 = 1/20$ along $x_1$ and $x_2$.

\begin{figure}
    \centering
    \begin{tabular}{cc}
            \hspace{-0.5cm}    
      \includegraphics[width=2.8in]{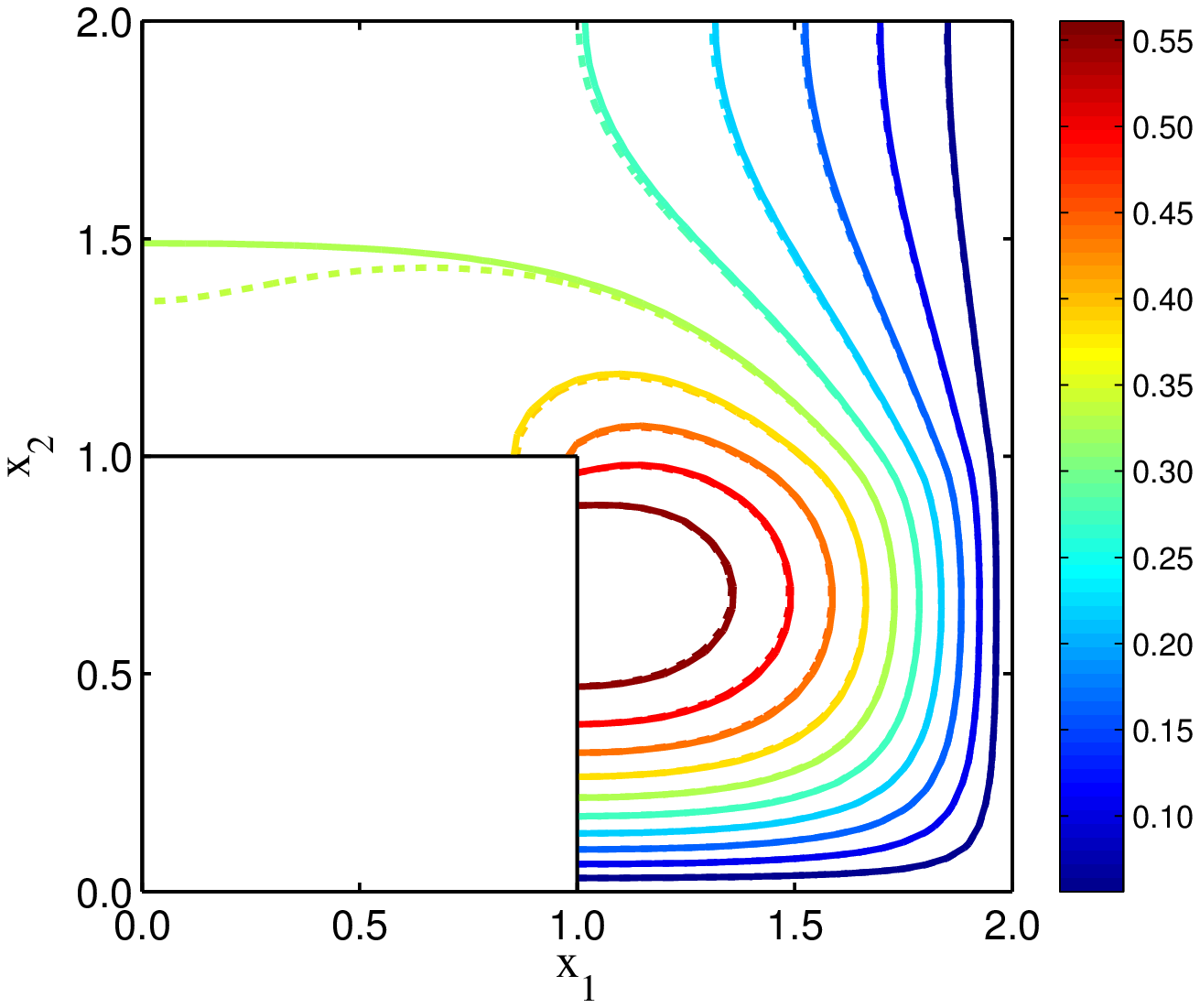}  
      &
      \includegraphics[width=2.8in]{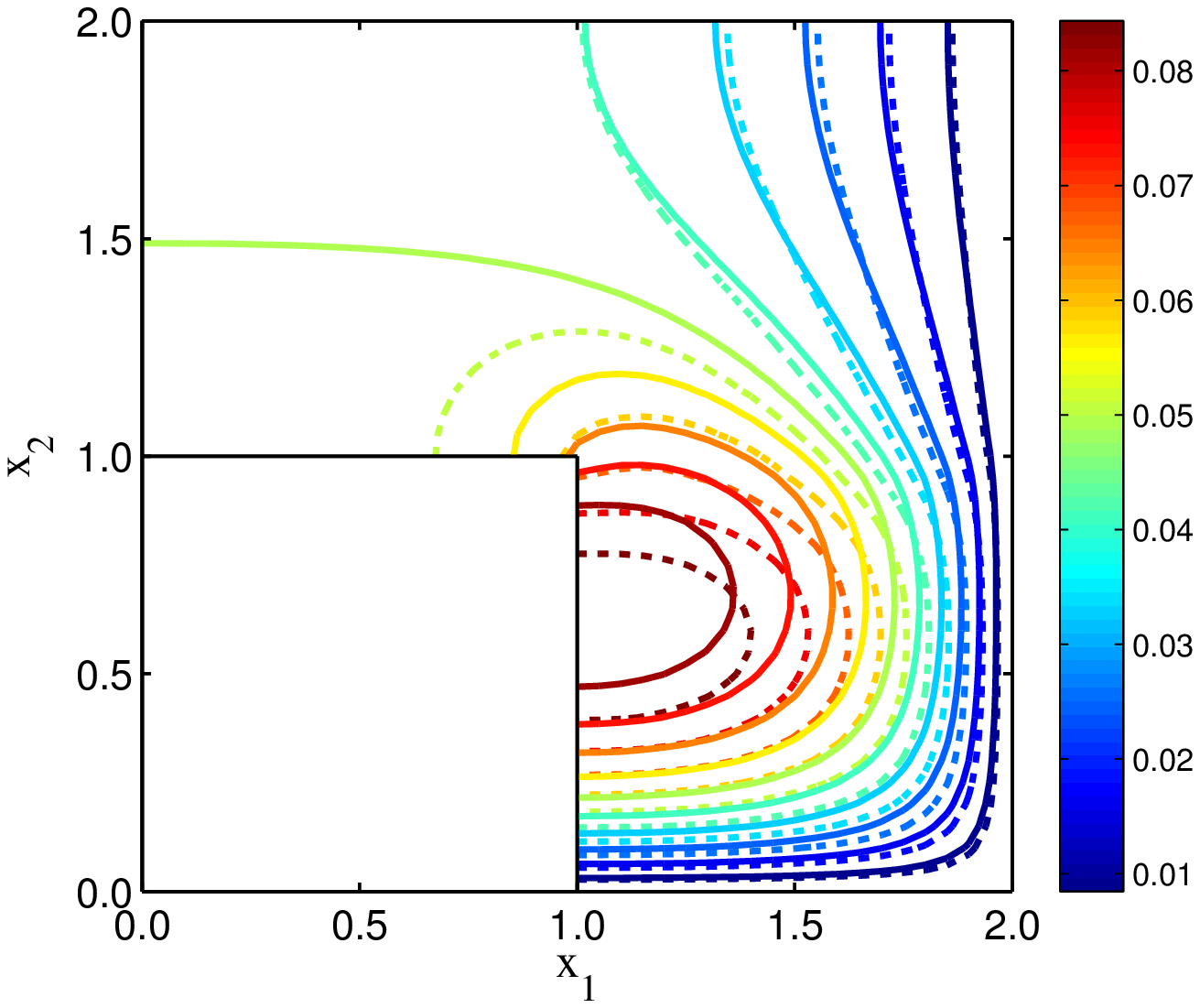} 
      \\
      (a) & (b)
      \\
      \includegraphics[width=2.8in]{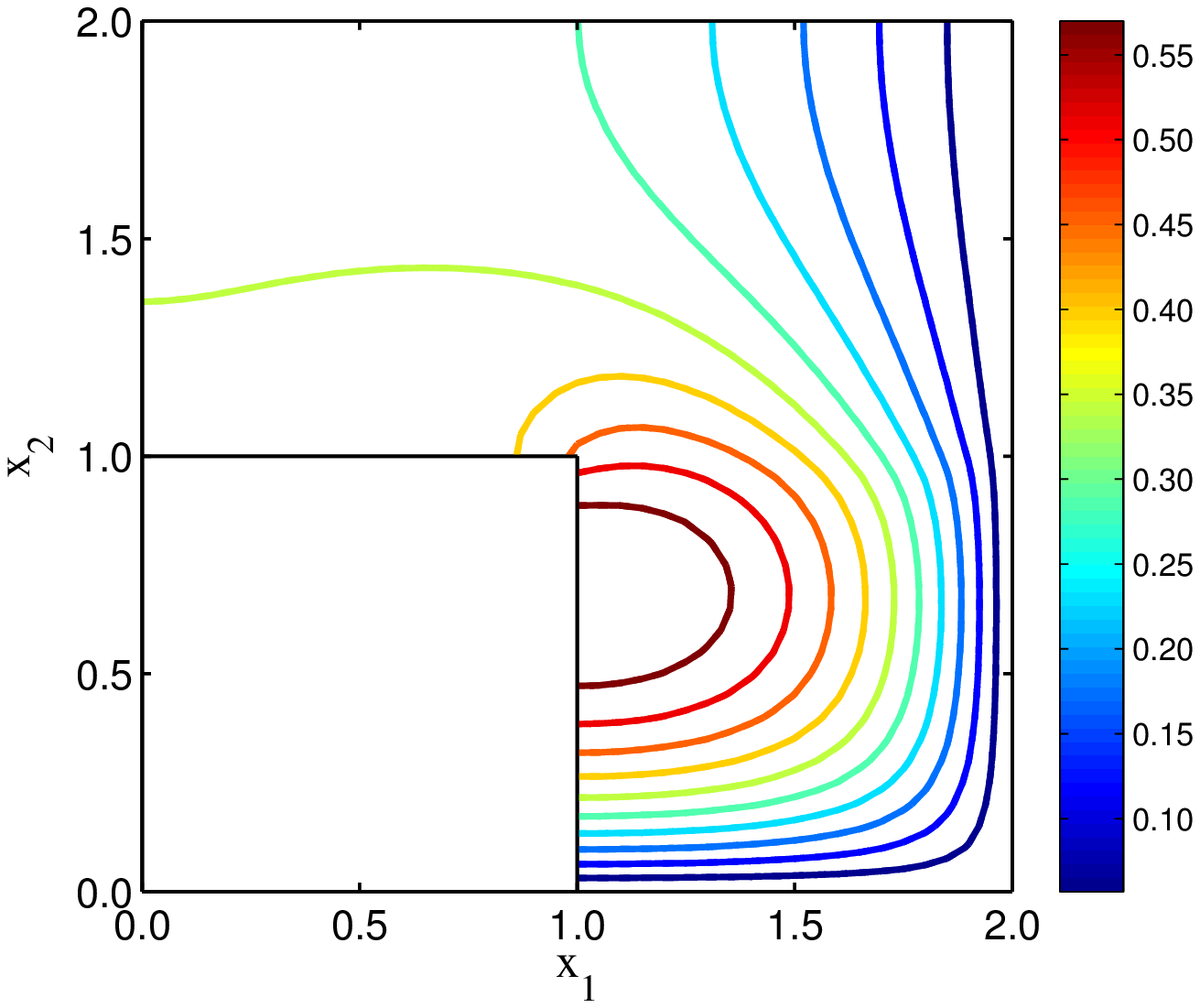}  
      &
      \includegraphics[width=2.8in]{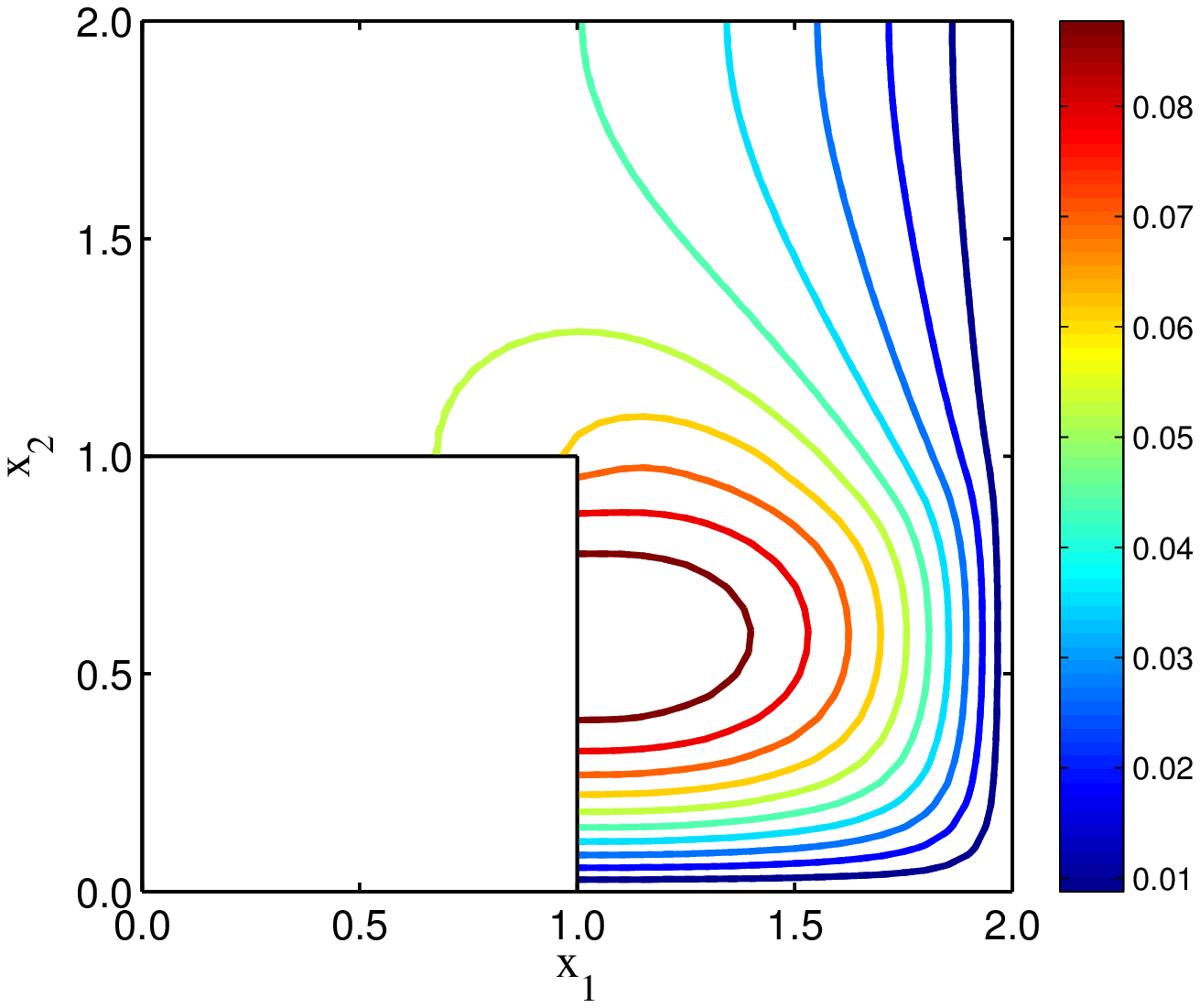}           
      \\
      (c) & (d)   
     \end{tabular}
      \caption{Contours of the solution mean and standard deviation obtained with separated representation (solid line) and the reference solution (dotted line) for the L-shaped problem. (a) Mean for $r=1$; (b) Standard deviation for $r=1$; (c) Mean for $r=20$; (d) Standard deviation for $r=20$.}            
\label{fig:L_contour}       
\end{figure}
\begin{figure}
  \begin{center}
	\includegraphics [width=10cm]{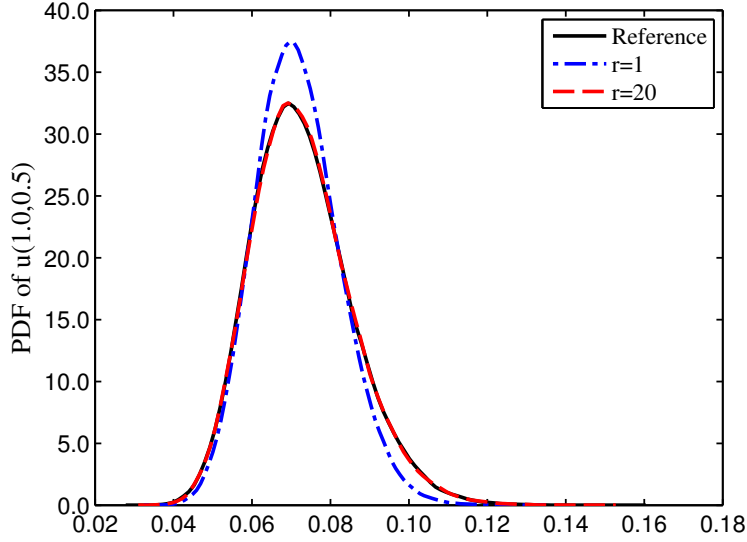}
    \caption{PDF of $u(1.0,0.5)$ for the L-shaped problem. A comparison between the rank $r$ separated representation and the reference solution.}
    \label{fig:pdf_L}
   \end{center}
\end{figure}
\begin{figure}
    \centering
    \begin{tabular}{cc}
            \hspace{-0.5cm}    
      \includegraphics[width=2.8in]{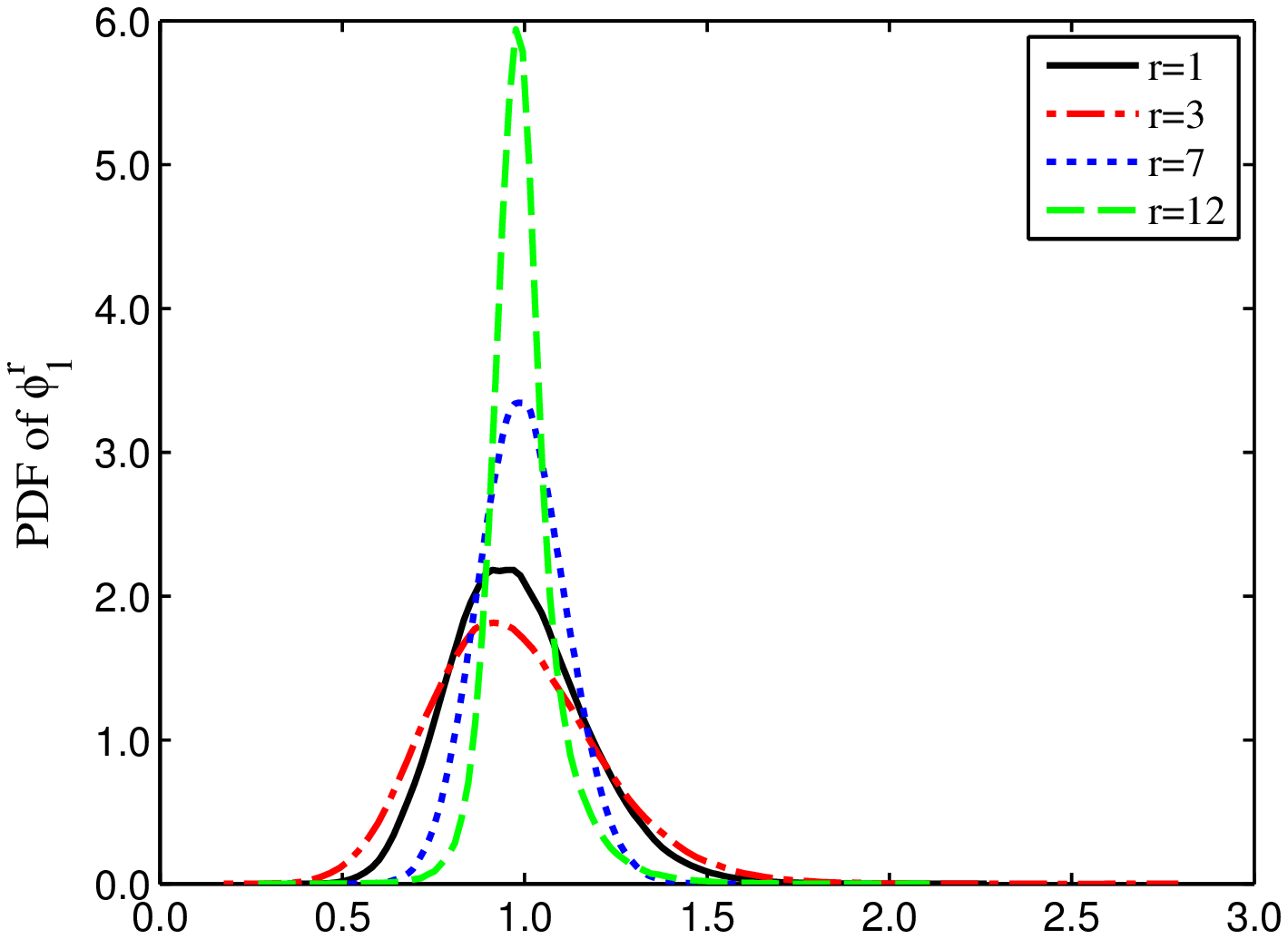}  
      &
      \includegraphics[width=2.8in]{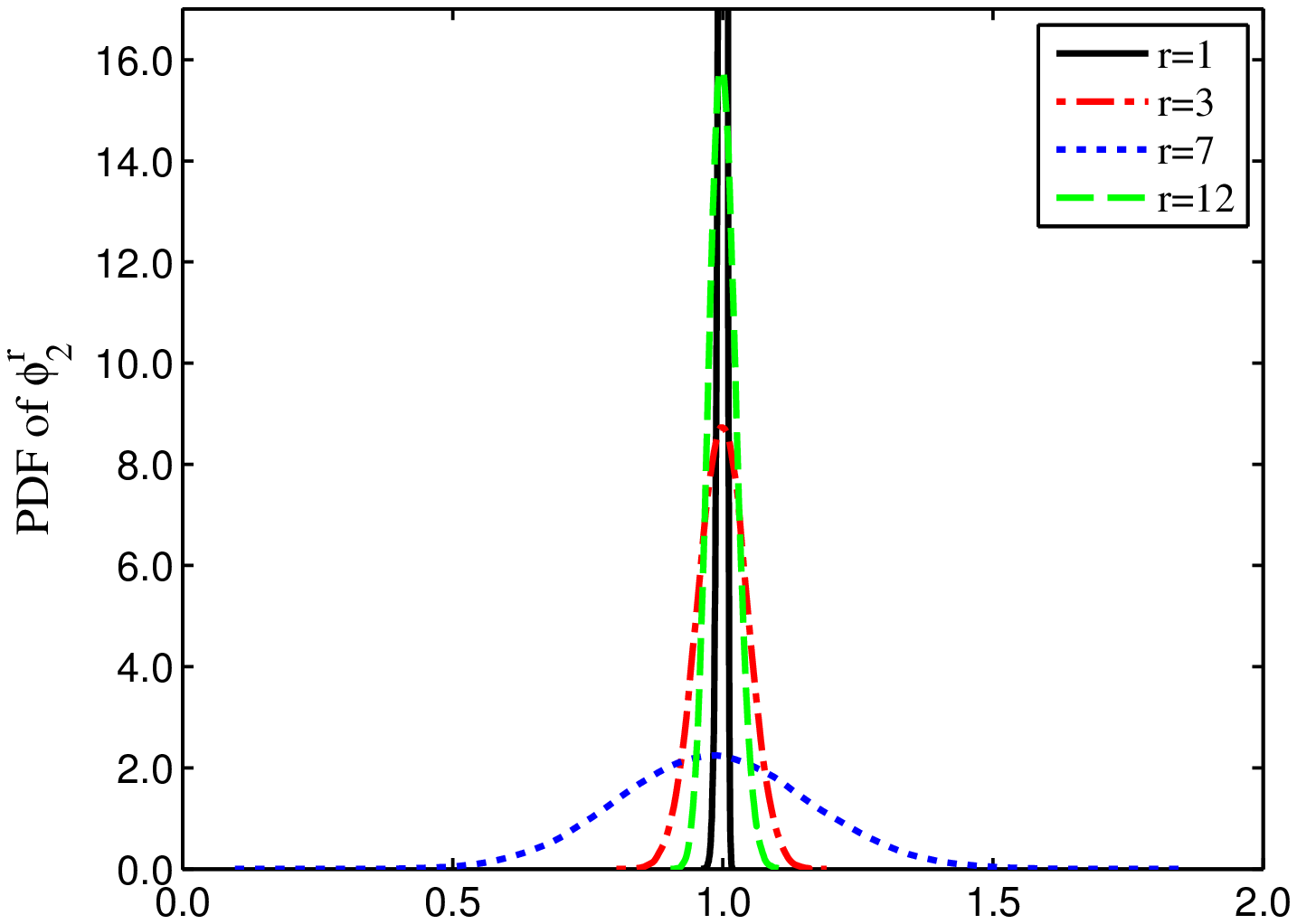} 
      \\
      (a) & (b)  
     \end{tabular}
      \caption{PDFs of the normalized $\phi_1^r$ and $\phi_2^r$ in the separated representation of the solution to the L-shaped problem. (a) PDF of the normalized $\phi_1^r$, $r=1,3,7,12$; (b) PDF of the normalized $\phi_2^r$, $r=1,3,7,12$. Each $\phi_i^r$ is normalized such that it has unit second moment. The PDF of $\phi_2^1$ is not fully shown here.}            
\label{fig:L_pdf_phi1_phi2}       
\end{figure}
\begin{figure}[htb]
  \begin{center}
	\includegraphics [width=10cm]{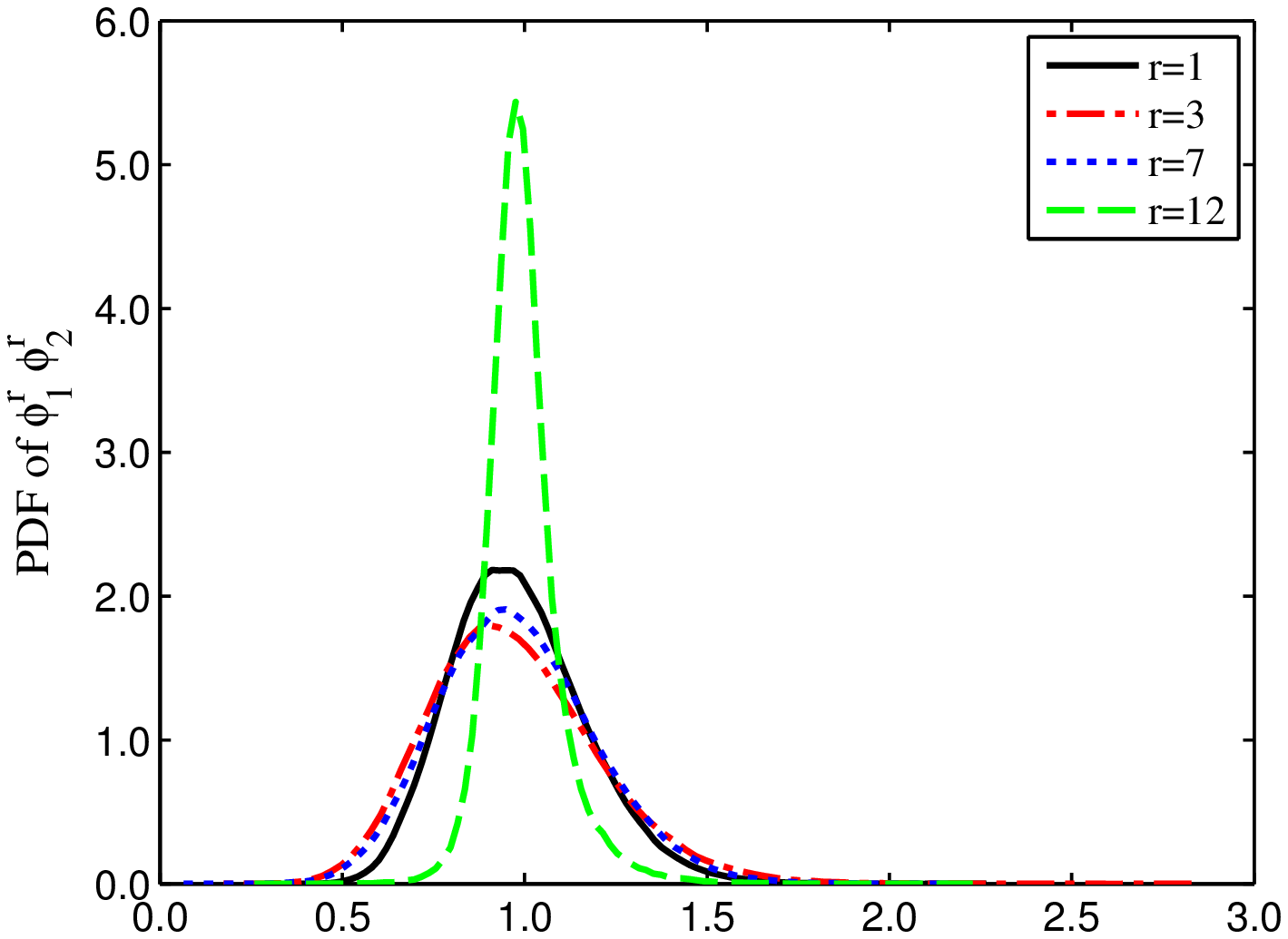}
    \caption{PDFs of the normalized stochastic basis functions $ \phi_1^{r}\phi_2^{r}$, $r=1,3,7,12$, for the L-shaped problem. Here each basis function is normalized to have unit second moment.}
    \label{fig:L_pdf_phi1phi2}
   \end{center}
\end{figure}

To demonstrate convergence of the separated representation of the solution $u$ to (\ref{eqn:elliptic_SPDE}), we increase the separation rank $r$ up to $r=20$. Alternatively, we may prescribe a target accuracy $\epsilon$ for the residual norm $\epsilon_{res}^r$ in (\ref{eqn:residual_error_general}) and identify the corresponding separation rank $r$. In the latter approach, the solution refinement is achieved by decreasing $\epsilon$. 

{\color{black}
Fig. \ref{fig:L_eng_error}(a) shows the values of the energy functional $\pi$ with respect to increments in the separation rank $r$. While this is not generally the case, here an increase in $r$ leads to a monotonic reduction in $\pi$. Fig. \ref{fig:L_eng_error}(b) presents the mean and standard deviation errors $\epsilon_{\mu}^r$ and $\epsilon_{\sigma}^r$, respectively, as defined in (\ref{eqn:residual_error}). As can be observed from this figure, there is not monotonicity in the reduction of mean and standard deviation errors as a function of $r$.} 

In Figs. \ref{fig:L_contour}(a)-(b) we display the contours of the solution mean and standard deviation obtained from the separated representation with $r=1$ as well as the reference solution. While approximation of the solution mean may be relatively accurate, the standard deviation has not yet converged with $r=1$. By increasing the separation rank to $r=20$, however, the approximation of these quantities improves considerably as can be observed from Figs. \ref{fig:L_contour}(c)-(d). 

In Fig. \ref{fig:pdf_L}, the probability density function (PDF) of the separated representation of solution at $(x_1,x_2)=(1.0,0.5)$ is compared to that of the reference solution. For the case of small separation rank, $r=1$, there is a notable difference between the two PDFs. However, this disagreement reduces considerably as the separation rank is increased to $r=20$. 

As discussed in Section \ref{sec:SR}, the stochastic functions $\lbrace \phi_1^{l} \rbrace_{l=1}^{r}$ and $\lbrace \phi_2^{l} \rbrace_{l=1}^{r}$ are not fixed \textit{a priori} and are computed by iteratively solving for the saddle point of $\pi$. Therefore, their representation is problem-dependent. Fig. \ref{fig:L_pdf_phi1_phi2} presents the PDFs of $\phi_1^{r}$ and $\phi_2^{r}$ when $r=1,3,7,12$. As it can be observed, these stochastic functions differ from each other and are not identically distributed. The actual stochastic basis functions in (\ref{eqn:rank-r-approx}) are products of these functions, i.e., $ \phi_1^{r} \phi_2^{r} $. Fig. \ref{fig:L_pdf_phi1phi2} shows the PDFs of a number of these basis functions. We note that these basis functions are generally different from the PC basis functions.

\subsection{Example II: 2D stochastic linear elasticity}
\label{sec:elasticity}
\begin{figure}[h]
  \begin{center}
	\includegraphics [width=12cm]{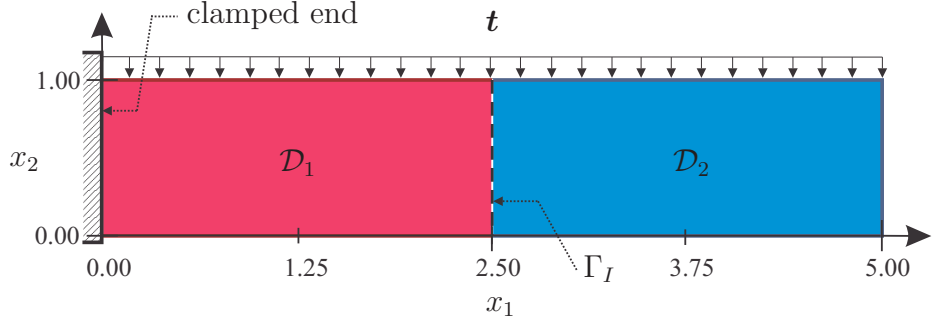}
    \put(-248,40){$\mathcal{D}_1$}
    \put(-99,40){$\mathcal{D}_2$}
    \put(-133,0){$\Gamma_I$}
    \put(-283,96){clamped end}
    \put(-170,93){$\bm t$}
    \put(-170,-13){$x_1$}
    \put(-350,42){$x_2$}
    \caption{Geometry of the 2D cantilever beam composed of two non-overlapping sub-domains $\mathcal D_1$ and $\mathcal D_2$. Young's moduli of the sub-domain materials are statistically independent random fields.}
    \label{fig:beam}
   \end{center}
\end{figure}

For the second numerical example, we consider the linear elasticity problem,
\begin{eqnarray} 
\label{eqn:elasticity}
&&-\nabla \cdot \bm\sigma \left( \bm u(\bm{x},\bm{\xi}) \right) = \bm 0 \qquad \qquad \bm{x} \in \mathcal{D}_, \nonumber \\
&&\bm u(\bm{x},\bm{\xi})=\bm 0 \quad \qquad \qquad \qquad \ \ \  \bm{x} \in \partial \mathcal{D}^{(D)}, \nonumber \\
&&\bm \sigma \left( \bm u(\bm{x},\bm{\xi}) \right) \bm{n} =\bm t(\bm x) \qquad  \qquad\, \bm{x} \in \partial \mathcal{D}^{(N)},
\end{eqnarray}
describing the deformation of the cantilever beam shown in Fig. \ref{fig:beam}. Here, $\bm u$ denotes the displacement vector field, $\bm\sigma$ is the stress tensor, $\bm n$ is the unit normal vector to the boundary, and $\bm t$ is the traction vector. We assume that $\bm t = (0,-0.1)$ on the edge corresponding to $x_2 = 1.0$ and $\bm t = (0,0)$ elsewhere. The stress tensor $\bm \sigma$ is related to the strain tensor $\bm e=\left(\nabla \bm u  + {\nabla \bm u}^T\right)/2$ via the isotropic linear elastic stress-strain relation 
\begin{equation} 
\label{eqn:stress_strain_relation}
\bm\sigma = \frac{E}{1+\nu} \left( \bm e + \frac{\nu}{1-2\nu}\ \mathrm{tr}(\bm e) \bm I \right), \nonumber
\end{equation}
where $\nu=0.3$ is the Poisson's ratio. Here, $\mathrm{tr}(\cdot)$ denotes the trace operator of a tensor and $\bm I$ is the identity tensor of order two. The Young's modulus $E(\bm x,\bm\xi)$  is the source of uncertainty in (\ref{eqn:elasticity}) and is assumed to take statistically independent values over sub-domains $\mathcal{D}_1$ and $\mathcal{D}_2$ shown in Fig. \ref{fig:beam}. Specifically,
\begin{equation}
\label{eqn:E_main}
E(\bm x,\bm \xi) = \left\{\begin{array}{c}E_1(\bm x,\bm \xi_1)\quad \bm x\in\mathcal{D}_1 \\ E_2(\bm x,\bm \xi_2)\quad \bm x\in\mathcal{D}_2 \end{array}\right.,
\end{equation}
where each $E_i(\bm x,\bm \xi_i)$, $i=1,2$, is represented by the series
\begin{equation}
\label{eqn:KL_YM}
E_i(\bm{x},\bm{\xi}_i) = \bar {E}_i + \sum_{j=1}^{d_i} \sqrt{\tau_{i,j}} g_{i,j}(\bm{x}) \xi_{i,j},\qquad \bm{x}\in\mathcal{D}_i.
\end{equation}

Here, $\left\lbrace \tau_{i,j} \right\rbrace_{j=1}^{d_i}$ and $\left\lbrace g_{i,j}(\bm{x}) \right\rbrace_{j=1}^{d_i}$ are $d_i$ largest eigenvalues and the corresponding eigenfunctions of the  Gaussian covariance kernel given in (\ref{eqn:covariance_kernel_k}). We assume that $\left\lbrace \xi_{1,j} \right\rbrace_{j=1}^{d_1}$ and $\left\lbrace \xi_{2,j} \right\rbrace_{j=1}^{d_2}$ are i.i.d. uniform random variables $U(-1,1)$. The list of parameters used in the analyses of this example is given in Table \ref{table:parameters_ex2}. These choices ensure that all realizations of $E$ are strictly positive on $\mathcal{D}$.

\begin{table}[h]
\label{table:parameters_ex2} 
\caption{Assumed parameters for the description of Young's modulus $E$ in (\ref{eqn:E_main}).} 
\centering
\begin{tabular}{ c c c c c c c c c c c }   
\hline
$d_1$ & $d_2$ & $\nu$ & $l_{c,1}$ & $l_{c,2}$ & $\bar E_1$ & $\bar E_2$ & $\sigma_1$ & $\sigma_2$  \\ 
\hline
9 & 11 & 0.3 & 2/3 & 1/3 & 100 & 100 & 35 & 35 \\ 
\hline 
\end{tabular} 
\label{table:parameters_ex2} 
\end{table}

Similar to the L-shaped problem, FE discretizations are done via the FEniCS project \citep{Logg12} and by using linear triangle elements with a uniform grid size $h_1 = h_2 = 1/10$ along $x_1$ and $x_2$ directions. Legendre PC expansions of $\lbrace \phi_1^{l} \rbrace_{l=1}^{r}$ and $\lbrace \phi_2^{l} \rbrace_{l=1}^{r}$ with degree $p=3$ were found sufficient for the solution of update equations (\ref{eqn:lin_sys_1_rg}) and (\ref{eqn:lin_sys_2_rg}), respectively. 

We note that in the present test case, $\mathcal D_2$ is a floating sub-domain with no Dirichlet boundary conditions; therefore, the PCPG solver described in Section \ref{sec:FETI} is applied to compute the vectors of Lagrange multipliers $\{\bm\lambda_0^l\}_{l=1}^r$ in (\ref{eqn:lin_sys_0_r}). We refer to Algorithm \ref{Algorithm:PCPG} for the implementation details of the PCPG solver. In our computations, we have considered $\epsilon_{PCPG}=10^{-8}$ as the stopping criterion for the PCPG solver. The condition number of $\hat{\bm{F}_I}$ determines the number of PCPG iterations and generally increases as a function of the separation rank $r$, hence asking for a preconditioner. Selection of an effective preconditioner $\bar{\bm{F}}_I^{-1}$, however, is generally a non-trivial task in the FETI approach. Among several existing choices, e.g., \citep{Farhat99,Klawonn01,Toselli99a,Rapetti01,Galvis10,Charmpis02}, we found the preconditioner of \cite{Klawonn01},
\begin{equation}
\label{eqn:feti_precond}
\bar{\bm{F}}_I^{-1} = \left(\bm{\hat{C}}^{T}_{1} \bm{\hat{C}}_{1}+\bm{\hat{C}}^{T}_{2} \bm{\hat{C}}_{2}\right)^{-1} \left( \bm{\hat{C}}^{T}_{1} \bm{\hat{K}}_{1} \bm{\hat{C}}_{1}+\bm{\hat{C}}^{T}_{2} \bm{\hat{K}}_{2} \bm{\hat{C}}_{2} \right) \left(\bm{\hat{C}}^{T}_{1} \bm{\hat{C}}_{1}+\bm{\hat{C}}^{T}_{2} \bm{\hat{C}}_{2}\right)^{-1},
\end{equation}
particularly effective for our purpose.

\begin{figure}[htb] 
    \centering
    \begin{tabular}{cc}
            \hspace{-0.5cm}    
      \includegraphics[width=2.8in]{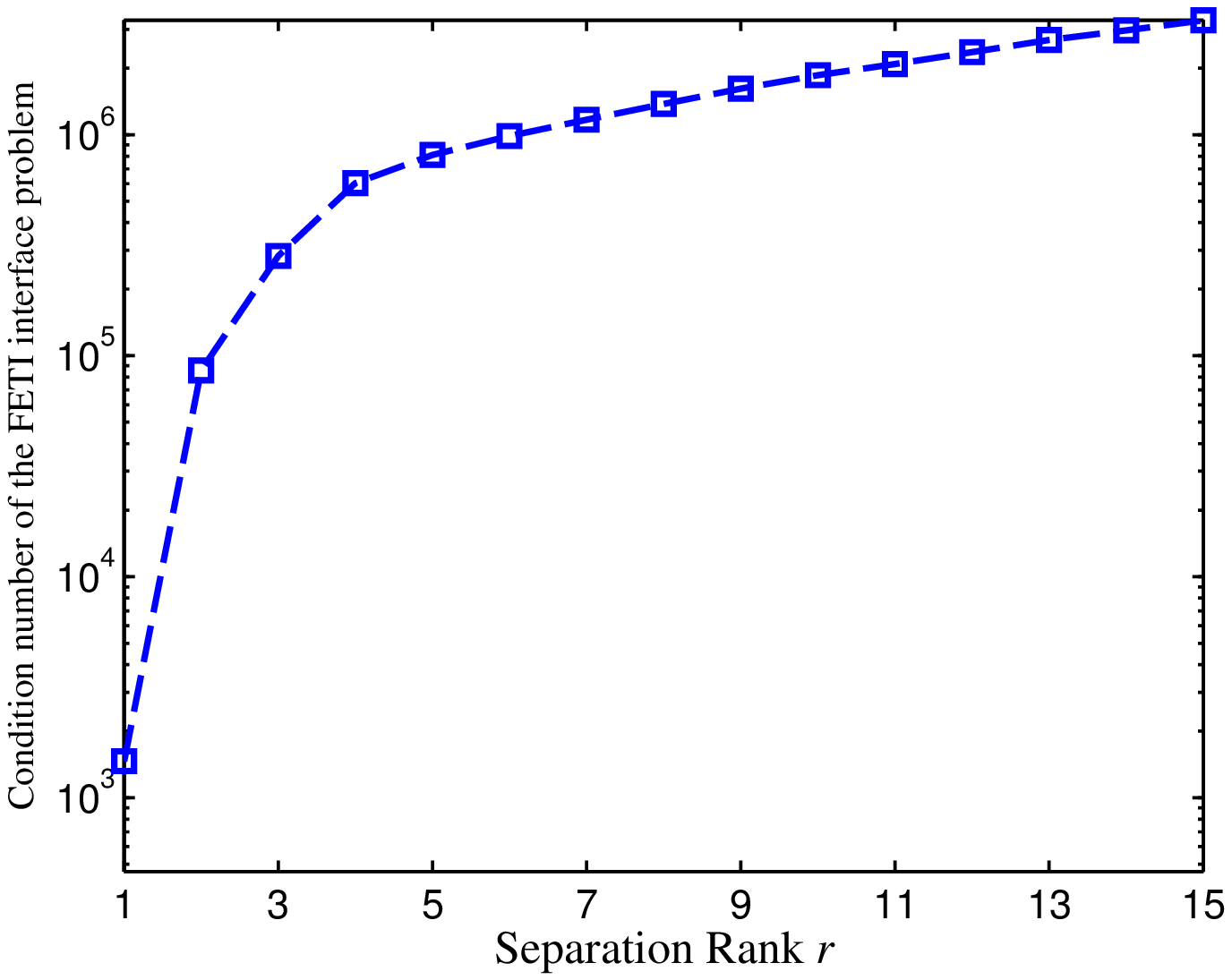}  
      &
      \includegraphics[width=2.8in]{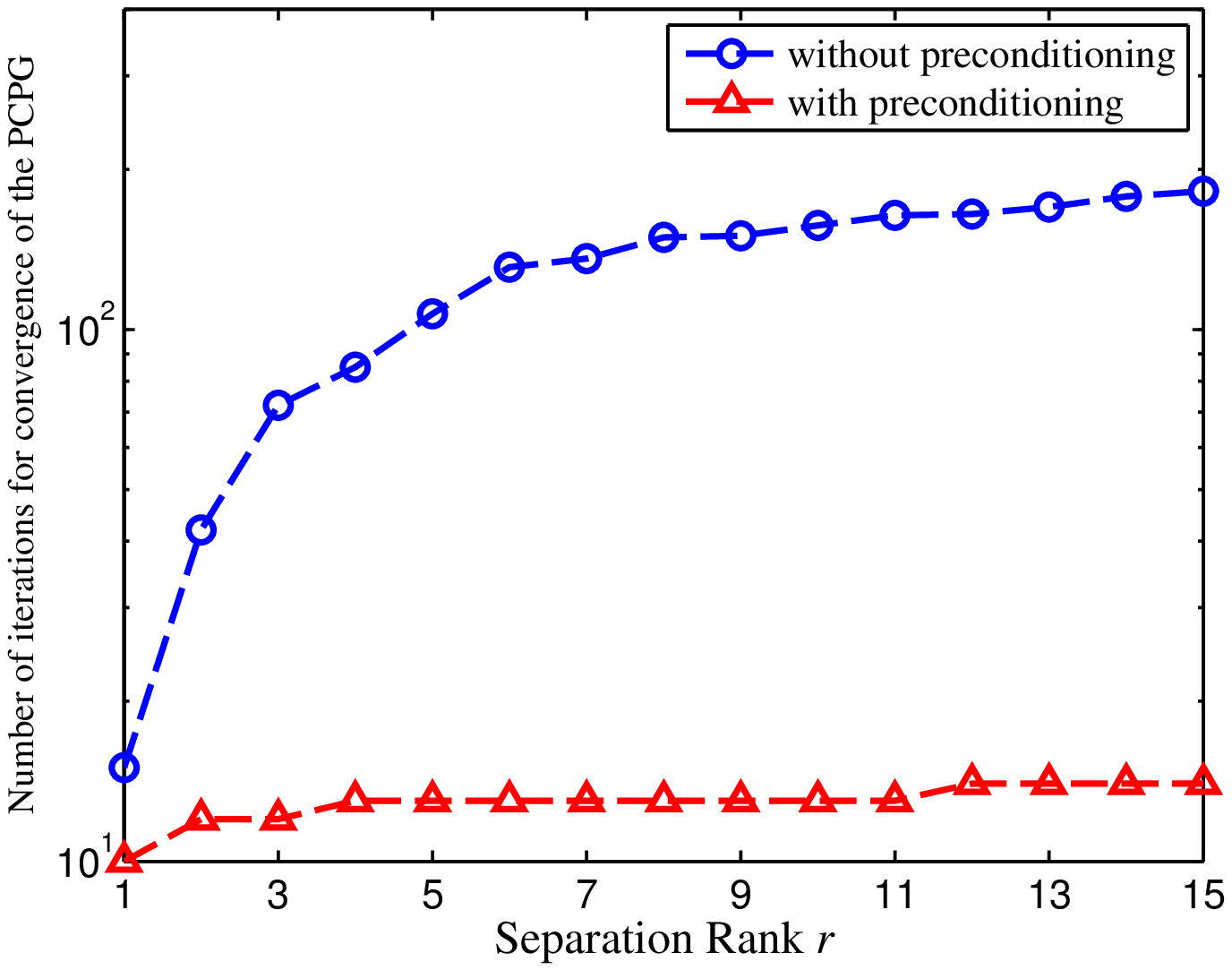} 
      \\
      (a) & (b)  
     \end{tabular}
      \caption{Convergence properties of the PCPG algorithm used for the beam problem. (a) Condition number of $\hat{\bm{F}}_I$ vs. the separation rank; (b) Number of the iterations required for the PCPG algorithm to converge with and without preconditioning vs. the separation rank.}            
\label{fig:B_FETI}       
\end{figure}

In Fig. \ref{fig:B_FETI}(a), we show the dependence of the condition number of $\hat{\bm{F}_I}$ on $r$ for the present example. A rapid increase in the condition number of $\hat{\bm{F}_I}$ can be observed when $r$ is increased. Fig. \ref{fig:B_FETI}(b) displays the number of PCPG iterations required to reach $\epsilon_{PCPG}=10^{-8}$ with the preconditioner (\ref{eqn:feti_precond}) and when no preconditioner is used, i.e., $\bar{\bm{F}}_I^{-1}=\bm I$. We observe that the choice of preconditioner (\ref{eqn:feti_precond}) makes the convergence of the PCPG solver almost independent of $r$. This is particularly crucial when one is dealing with problems in which $r$ is large. However, we note that further analysis is needed to confirm the effectiveness of (\ref{eqn:feti_precond}). 

\begin{figure}
    \centering
    \begin{tabular}{cc}
            \hspace{-0.5cm}    
      \includegraphics[width=2.8in]{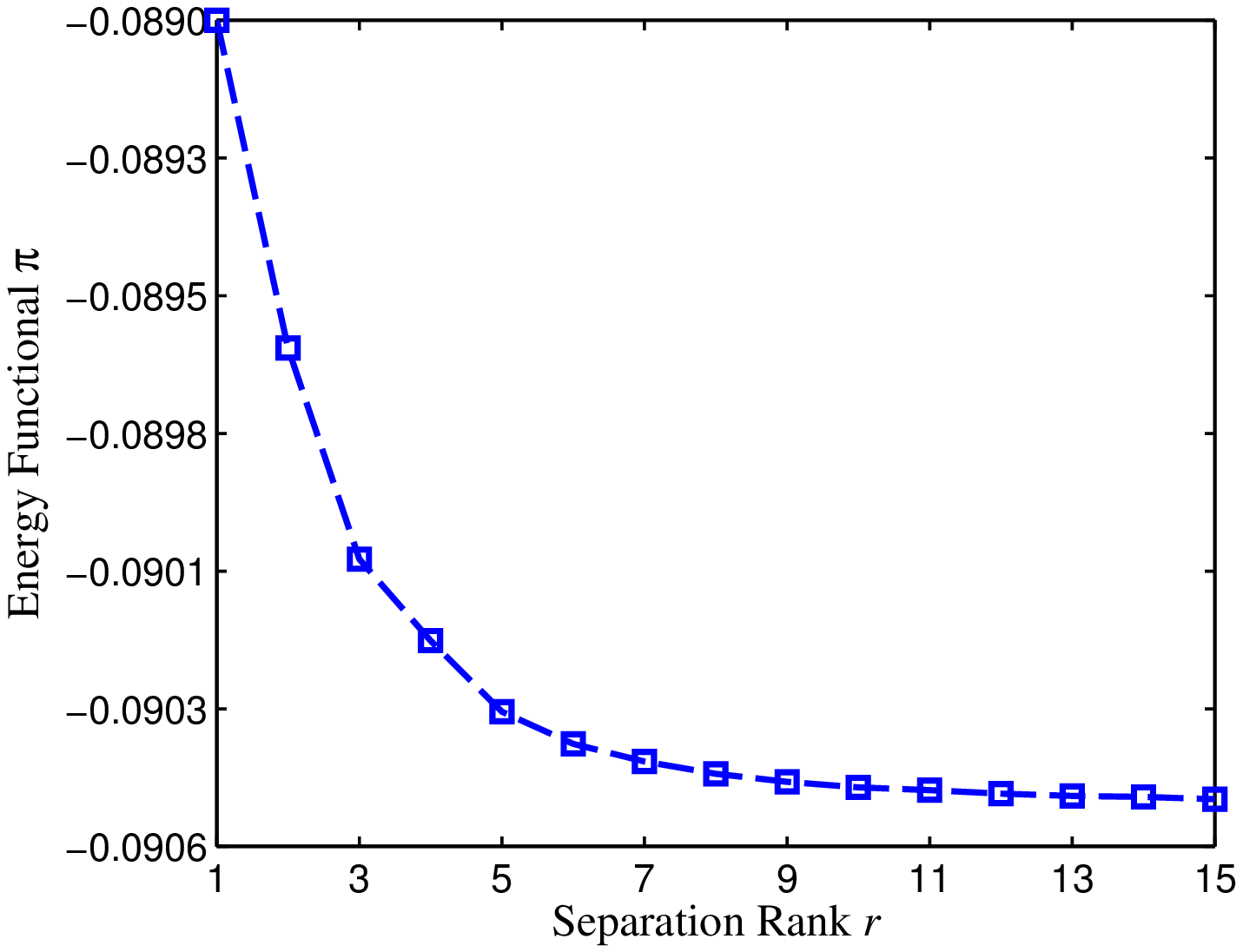}  
      &
      \includegraphics[width=2.8in]{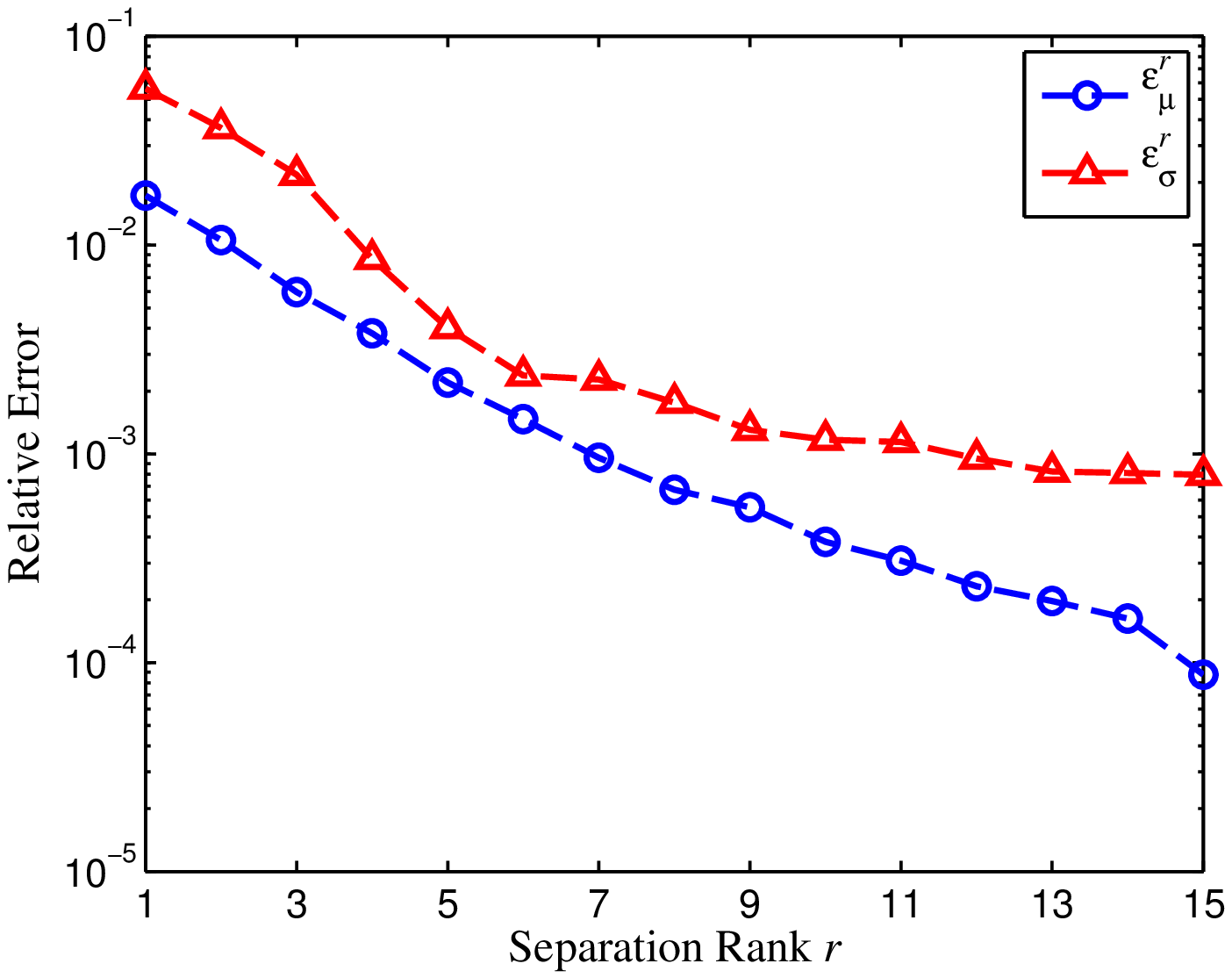} 
      \\
      (a) & (b)  
     \end{tabular}
      \caption{Energy functional $\pi$ (a) and relative errors in mean and standard deviation of the displacement (b) as a function of the separation rank $r$ for the beam problem. The errors are evaluated from (\ref{eqn:residual_error}).}            
\label{fig:B_eng_error}       
\end{figure}

Fig. \ref{fig:B_eng_error} illustrates the convergence of the separated representation. {\color{black}Similar to the L-shaped problem, increasing the separation rank $r$ results in a monotonic descries of the energy functional $\pi$},  see Fig. \ref{fig:B_eng_error}(a). The relative mean and standard deviation errors $\epsilon_{\mu}^r$ and $\epsilon_{\sigma}^r$ of the displacement are computed from (\ref{eqn:residual_error_general}) and are plotted against $r$ in Fig. \ref{fig:B_eng_error}(b).

\begin{figure}
    \centering
    \begin{tabular}{cccc} 
	 (a)   
      &
      \includegraphics[width=5in]{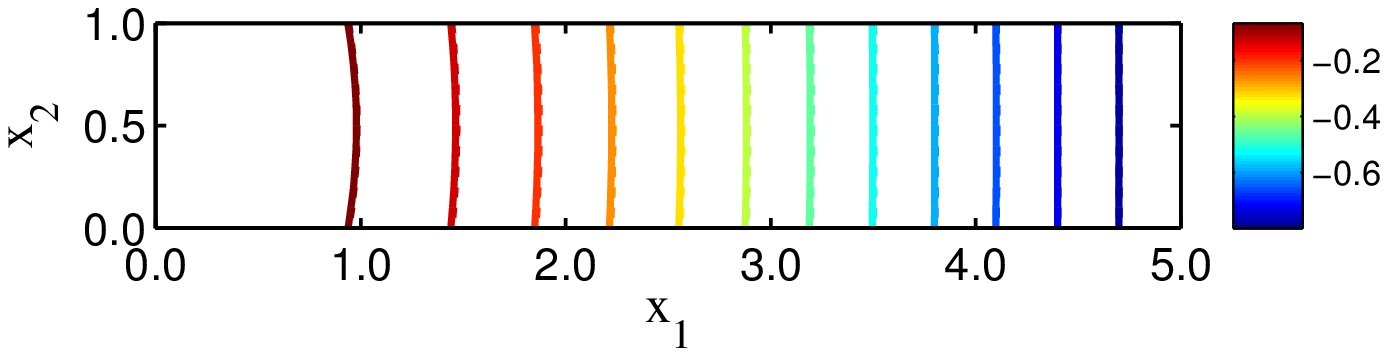}  
      \\
	 (b)  
      &
      \includegraphics[width=5in]{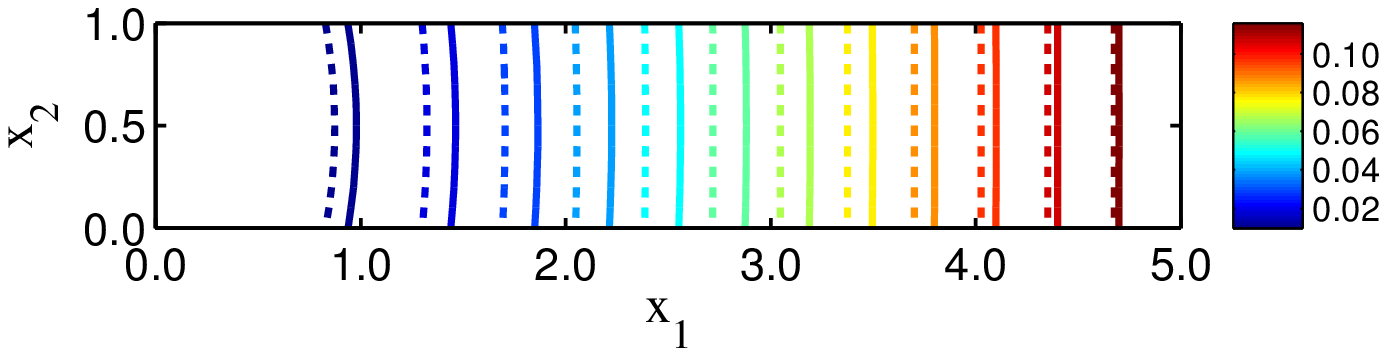} 
      \\
	 (c)  
      &
      \includegraphics[width=5in]{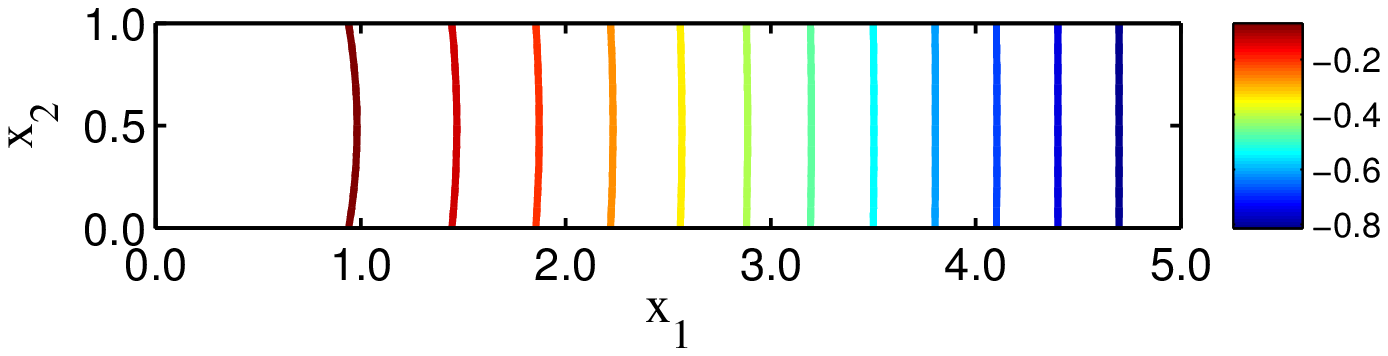}           
      \\
	 (d)  
      &
      \includegraphics[width=5in]{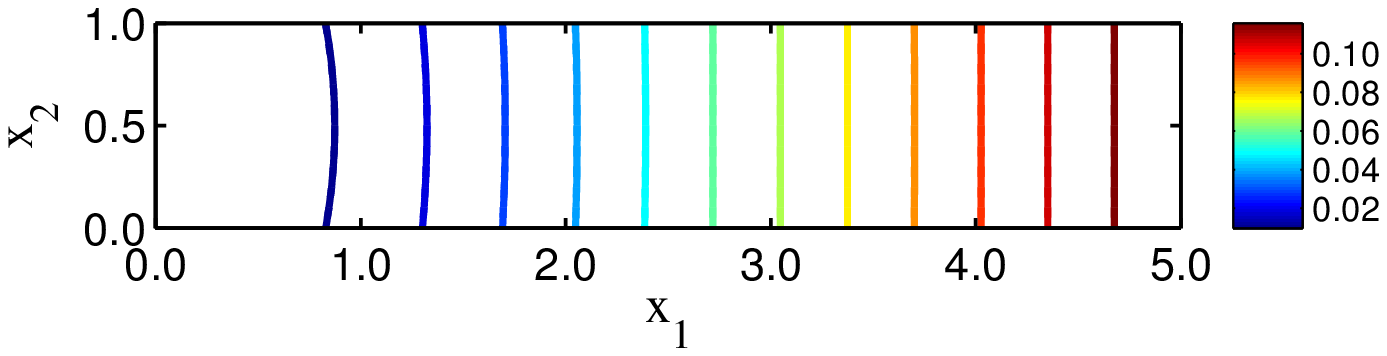}    
     \end{tabular}
      \caption{Contours of the mean and standard deviation of vertical displacement obtained with separated representation (solid line) and the reference solution (dotted line) for the beam problem. (a) Mean for $r=1$; (b) Standard deviation for $r=1$; (c) Mean for $r=15$; (d) Standard deviation for $r=15$.}            
\label{fig:B_contour}       
\end{figure}

In Fig. \ref{fig:B_contour}, we display the contours of the mean and standard deviation of the vertical displacement obtained from the separated representation, and compare them with those of the reference solution. The mean is captured fairly accurately with a rank one approximation. While for $r=1$ a considerable deviation from the standard deviation of the reference solution is observed, the rank $r=15$ approximation agrees well with the reference solution (Fig. \ref{fig:B_contour}(b) and Fig. \ref{fig:B_contour}(d)). Despite the high-dimensionality of the solution, i.e., $d=20$, we note that we only require a low separation rank, $r=15$, to accurately approximate the solution mean and standard deviation. 

\begin{figure}[htb]
  \begin{center}
	\includegraphics [width=10cm]{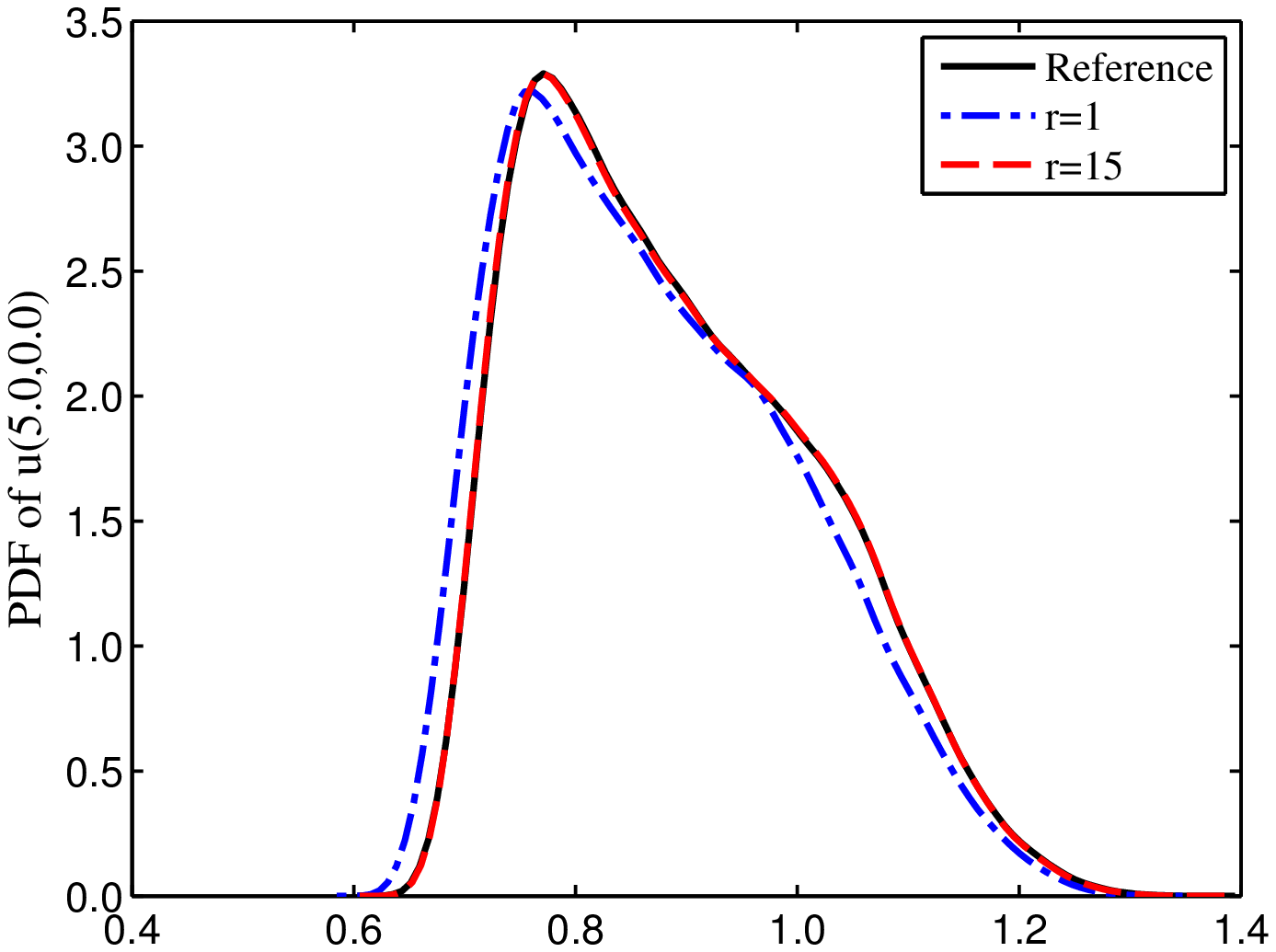}
    \caption{PDF of the total displacement at $(x_1,x_2)=(5.0,0.0)$. A comparison between the separated approximation and the reference solution.}
    \label{fig:pdf_B}
   \end{center}
\end{figure}

Fig. \ref{fig:pdf_B} compares the PDFs of the separated approximation of the total displacement at $(x_1,x_2)=(5.0,0.0)$, when $r=1$ and $r=15$, to the PDF of the corresponding reference solution. An almost identical agreement between these PDFs is observed when $r=15$.

\begin{figure}[htb] 
    \centering
    \begin{tabular}{cc}
            \hspace{-0.5cm}    
      \includegraphics[width=2.8in]{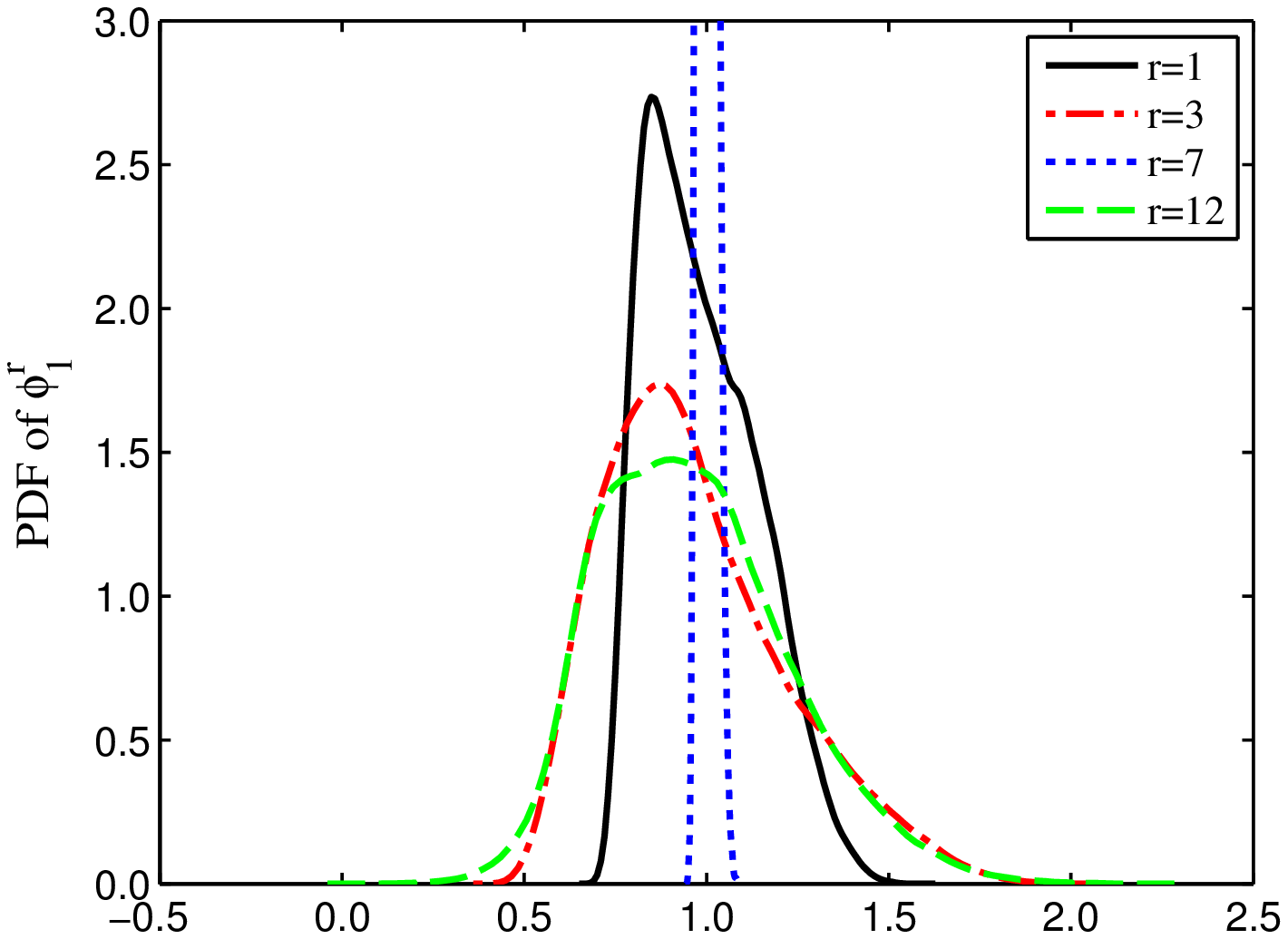}  
      &
      \includegraphics[width=2.8in]{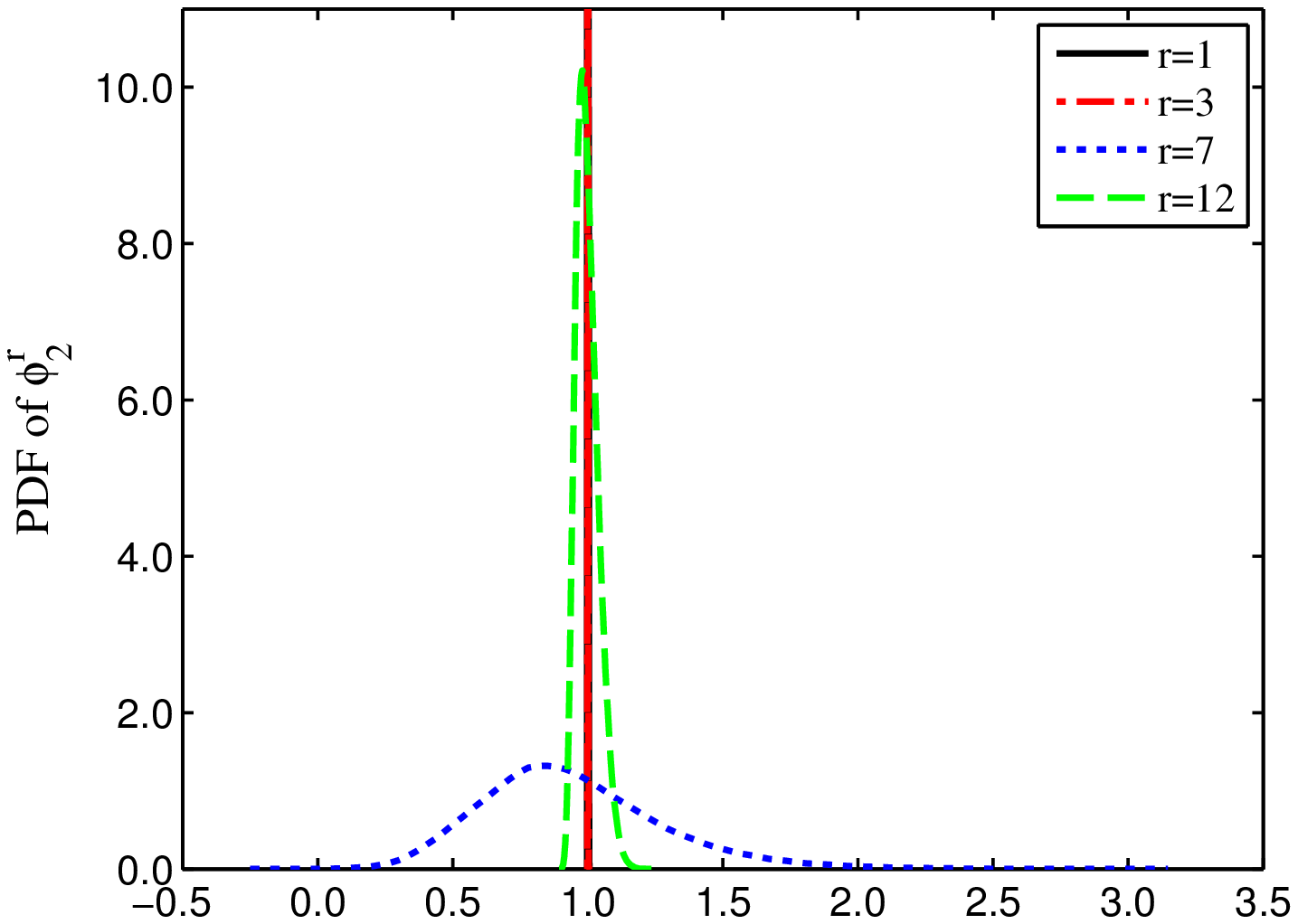} 
      \\
      (a) & (b)  
     \end{tabular}
      \caption{PDFs of the normalized stochastic functions $\phi_1^r$ and $\phi_2^r$, $r=1,3,7,12$, for the beam problem. The normalized $\phi_1^r$ and $\phi_2^r$ have unit second moment. (a) PDFs of the normalized $\phi_1^r$; (b) PDFs of the normalized $\phi_2^r$. PDFs of $\phi_1^7$, $\phi_2^1$, and $\phi_2^3$ are not shown completely.}            
\label{fig:B_pdf_phi1_phi2}       
\end{figure}
\begin{figure}[h]
  \begin{center}
	\includegraphics [width=10cm]{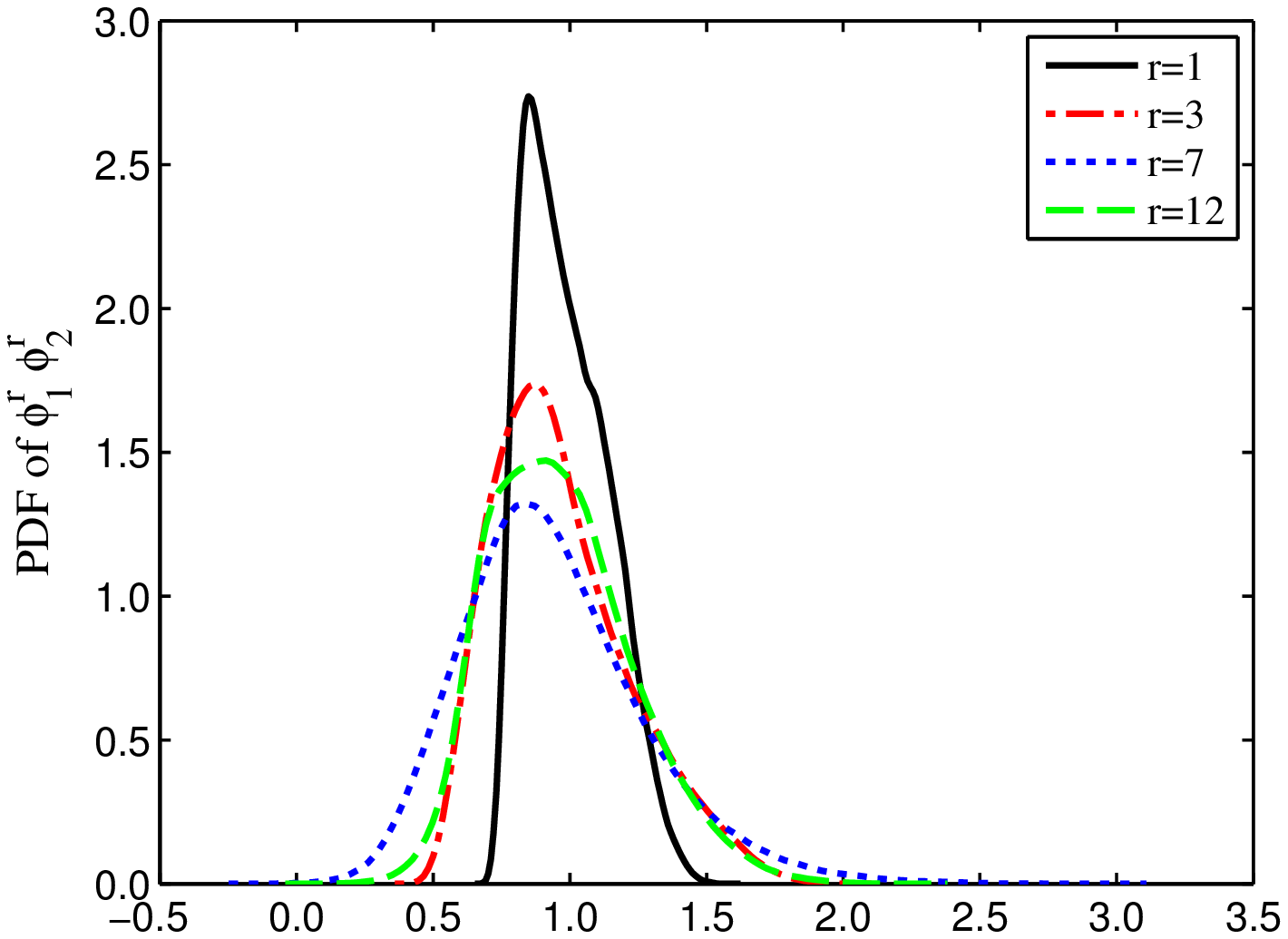}
    \caption{PDFs of the normalized stochastic basis functions $\phi_1^{r}\phi_2^{r}$, $r=1,3,7,12$, for the beam problem. Each $\phi_1^{r}\phi_2^{r}$ has unit second moment.}
    \label{fig:B_pdf_phi1phi2}
   \end{center}
\end{figure}

To demonstrate that the stochastic functions $\phi_1^{r}$ and $\phi_2^{r}$ in (\ref{eqn:rank-r-approx}) depend on the problem at hand, we present the PDFs of these quantities, for $r=1,3,7,12$, in Fig. \ref{fig:B_pdf_phi1_phi2}. As it can be seen, these PDFs are different from those of the L-shaped problem. Additionally, Fig. \ref{fig:B_pdf_phi1phi2} reports the PDF of the stochastic basis functions $\phi_1^{r}\phi_2^{r}$ for the above $r$ values. 

\section{Conclusion}
\label{sec:conclusion}

We have developed a computational framework for the propagation of uncertainty through coupled domain problems. The proposed approach hinges on the construction of a solution-adaptive stochastic basis that is of separated form with respect to the random inputs characterizing the uncertainty in each sub-domain. Such a separated construction of stochastic basis is achieved through a sequence of approximations with respect to the dimensionality, i.e., number of random inputs, of each individual sub-domain and not the combined dimensionality. This leads to a {\it partitioned} treatment of the stochastic space and, consequently, a higher scalability of the method as compared with standard uncertainty propagation approaches, such as those based on direct polynomial chaos expansions. For situations where the cardinality of the separated basis -- here referred to as the separation rank -- is small, the proposed approach provides a reduced order representation of the coupled field solution.  

The deterministic coefficients associated with each separated stochastic basis capture the spatial variability of the solution and are computed via the standard finite element tearing and interconnecting (FETI) approach. Therefore, the method achieves a high level of parallelism while requiring no intrusion in each sub-domain solver. Although our present formulation of domain coupling is based on the standard FETI approach, we foresee no major technical difficulties in employing more advanced domain coupling schemes.

The performance of the proposed framework was explored through its application to two linear elliptic PDEs with high-dimensional random inputs. Both problems were defined on physical domains consisting of two coupled sub-domains with independent sources of uncertainty. In both cases, despite the high-dimensionality of the random inputs, accurate estimations of the solution statistics were achieved with relatively low separation ranks, thus demonstrating the effectiveness of the present approach. 

The stochastic expansion of this study based on the separated representations may also be applied to other coupled problems, such as fluid structure interaction (FSI), involving uncertainty. However, different numerical strategies for the construction of the separated basis may be required. 

\section*{Acknowledgements}

The authors are indebted to the fruitful discussions they had with Prof. K.C. Park and Prof. Kurt Maute from University of Colorado, Boulder. 

AD gratefully acknowledges the financial support of the Department of Energy under Advanced Scientific Computing Research Early Career Research Award DE-SC0006402. MH's work was supported by the National Science Foundation grant CMMI-1201207. {\color{black}The work of HGM and RN has been partially supported by the Deutsche Forschungsgemeinschaft (DFG) through the SFB 880.}

\clearpage

\appendix
\appendixpage
\addappheadtotoc
\section{}
\subsection{Linear system (\ref{eqn:lin_sys_0_r}) for $\bm{u}_{0,1}^l$, $\bm{u}_{0,2}^l$, and $\bm{\lambda}_0^l$ updates}
\label{apx:deterministic_updates}

We assume all unknowns are fixed except the vector $\bm{u}_{0,1}^{l}$. Taking the derivative of $\pi$ in (\ref{eqn:energy}) with respect to $\bm{u}_{0,1}^{l}$, $l=1,\dots,r$, leads to the following linear system
\begin{equation}
\label{eqn:v1-update_r}
\sum_{l'=1}^{r} \left( \mathbb{E}\left[\phi_{1}^{l}\phi_{2}^{l} \bm{K}_{1}\phi_{1}^{l'}\phi_{2}^{l'}\right]
\bm{u}_{0,1}^{l'}-\mathbb{E}\left[
\phi_{1}^{l}\phi_{2}^{l} \phi_1^{l'} \phi_2^{l'}\right]\bm{C}_{1}\bm{\lambda}_0^{l'} \right) =\mathbb{E}\left[\phi_{1}^{l}\phi_{2}^{l}\bm{f}_{1}\right],\quad l=1,\dots,r.
\end{equation}

Similarly, the $\bm{u}_{0,2}^{l}$ updates are performed through 
\begin{equation}
\label{eqn:v2-update_r}
\sum_{l'=1}^{r} \left( \mathbb{E}\left[\phi_{1}^{l}\phi_{2}^{l} \bm{K}_{2}\phi_{1}^{l'}\phi_{2}^{l'}\right]
\bm{u}_{0,2}^{l'}+\mathbb{E}\left[
\phi_{1}^{l}\phi_{2}^{l} \phi_1^{l'} \phi_2^{l'}\right]\bm{C}_{2}\bm{\lambda}_0^{l'} \right) =\mathbb{E}\left[\phi_{1}^{l}\phi_{2}^{l}\bm{f}_{2}\right],\quad l=1,\dots,r.
\end{equation}

The derivative of $\pi$ with respect to deterministic vector $\bm{\lambda}_0^{l}$ of the Lagrange multipliers satisfies 
\begin{equation}
\label{eqn:lam0-update_r}
\sum_{l'=1}^{r} \left( \mathbb{E}\left[ \phi_1^{l} \phi_2^{l}\phi_{1}^{l'}\phi_{2}^{l'}\right]\bm{C}_{2}^{T}\bm{u}_{0,2}^{l'}-\mathbb{E}\left[ \phi_1^{l} \phi_2^{l}\phi_{1}^{l'}\phi_{2}^{l'}\right]\bm{C}_{1}^{T}\bm{u}_{0,1}^{l'} \right)= \bm{0},\quad l=1,\dots,r.
\end{equation}

Putting (\ref{eqn:v1-update_r}), (\ref{eqn:v2-update_r}), and (\ref{eqn:lam0-update_r}) together results in the  linear system of (\ref{eqn:lin_sys_0_r}) for $\bm{u}_{0,1}^l$, $\bm{u}_{0,2}^l$, and $\bm{\lambda}_0^l$ updates.

\subsection{Linear system (\ref{eqn:lin_sys_2_rg}) for $\phi_1^{l}(\bm\xi_1)$ updates}
\label{apx:stochastic_updates}

We first remind that in $\phi_1^{l}(\bm\xi_1)$ updates the vectors $\bm{u}_{0,1}^l$, $\bm{u}_{0,2}^l$, and $\bm{\lambda}_0^l$ are fixed at their current values obtained from (\ref{eqn:lin_sys_0_r}). As the solution to (\ref{eqn:lin_sys_0_r}) satisfies the constraints in (\ref{eqn:lam0-update_r}), the term $\mathbb{E} \left[ \bm{\lambda}^{T}\left(\bm{C}_{2}^T \bm{u}_2 - \bm{C}_{1}^T \bm{u}_1 \right) \right]$ in (\ref{eqn:energy}) corresponding to the the interface condition vanishes for any $\{\phi_1^{l}(\bm\xi_1)\}_{l=1}^{r}$ and $\{\phi_2^{l}(\bm\xi_2)\}_{l=1}^{r}$. Therefore, the $\phi_1^{l}(\bm\xi_1)$ updates correspond to minimization of $\pi$. More precisely, plugging the separated representation (\ref{eqn:rank-r-approx}) into condition (\ref{eqn:mpe_{phi}_{1}^{l}}), we arrive at
\begin{eqnarray}
\label{eqn:phi_1_var_rg}
&&\int\int \left(
\sum_{l'=1}^{r}{\bm{u}^{l}}^{T}_{0,1} \phi_{2}^{l}(\bm\xi_2)\bm{K}_{1}(\bm\xi_1)\bm{u}_{0,1}^{l'} \phi_{1}^{l'}(\bm\xi_1) \phi_{2}^{l'}(\bm\xi_2) \right)\delta \phi_{1}^{l}(\bm\xi_1)\rho_1(\bm\xi_1)\mathrm{d}\bm\xi_1 \rho_2(\bm\xi_2)\mathrm{d}\bm\xi_2 \nonumber \\
&&+\int\int \left(\sum_{l'=1}^{r}{\bm{u}^{l}}^{T}_{0,2} \phi_{2}^{l}(\bm\xi_2)\bm{K}_{2}(\bm\xi_2)\bm{u}_{0,2}^{l'} \phi_{1}^{l'}(\bm\xi_1) \phi_{2}^{l'}(\bm\xi_2)\right)\delta \phi_{1}^{l}(\bm\xi_1)\rho_1(\bm\xi_1)\mathrm{d}\bm\xi_1 \rho_2(\bm\xi_2)\mathrm{d}\bm\xi_2  \nonumber \\
&&= \int\int \left({\bm{u}^{l}}^{T}_{0,1}\bm{f}_{1} \phi_{2}^{l}(\bm\xi_2)+{\bm{u}^{l}}^{T}_{0,2}\bm{f}_{2} \phi_{2}^{l}(\bm\xi_2)\right)\delta \phi_{1}^{l}(\bm\xi_1)\rho_1(\bm\xi_1)\mathrm{d}\bm\xi_1 \rho_2(\bm\xi_2)\mathrm{d}\bm\xi_2 ,\quad \forall\  \delta \phi_{1}^{l}(\bm\xi_1), \nonumber \\
\end{eqnarray}
where $\rho_1(\bm\xi_1)$ and $\rho_2(\bm\xi_2)$ are the joint probability density functions of $\bm\xi_1$ and $\bm\xi_2$, respectively. One can rewrite (\ref{eqn:phi_1_var_rg}) as
\begin{eqnarray}
\label{eqn:phi_1_var_rg_rep1}
&&\int \left(\sum_{l'=1}^{r}{\bm{u}^{l}}^{T}_{0,1} \bm{K}_{1}(\bm\xi_1)\bm{u}_{0,1}^{l'}\phi_{1}^{l'}(\bm\xi_1)\mathbb{E}_{\bm\xi_2}\left[\phi_{2}^{l}
\phi_{2}^{l'}\right]\right)\delta \phi_{1}^{l}(\bm\xi_1)\rho_1(\bm\xi_1)\mathrm{d}\bm\xi_1 \nonumber \\
&&+\int \left(\sum_{l'=1}^{r}{\bm{u}^{l}}^{T}_{0,2}\mathbb{E}_{\bm\xi_2}
\left[\bm{K}_{2}\phi_{2}^{l}\phi_{2}^{l'}\right]
\bm{u}_{0,2}^{l'}\phi_{1}^{l'}(\bm\xi_1)\right)\delta \phi_{1}^{l}(\bm\xi_1)\rho_1(\bm\xi_1)\mathrm{d}\bm\xi_1 \nonumber \\
&&= \int \left({\bm{u}^{l}}^{T}_{0,1}\bm{f}_{1}
\mathbb{E}_{\bm\xi_2}\left[\phi_{2}^{l}\right]+
{\bm{u}^{l}}^{T}_{0,2}\bm{f}_{2}\mathbb{E}_{\bm\xi_2}
\left[\phi_{2}^{l}\right]\right)\delta \phi_{1}^{l}(\bm\xi_1)\rho_1(\bm\xi_1)\mathrm{d}\bm\xi_1 ,\quad \forall\  \delta \phi_{1}^{l}(\bm\xi_1),\nonumber
\end{eqnarray}
which, {\color{black}if the integrands are continuous}, is equivalent to
\begin{eqnarray}
\label{eqn:phi_1_var_rg_rep2}
&&\sum_{l'=1}^{r}{\bm{u}^{l}}^{T}_{0,1}\bm{K}_{1}(\bm\xi_1)\bm{u}_{0,1}^{l'}\phi_{1}^{l'}(\bm\xi_1)\mathbb{E}_{\bm\xi_2}\left[\phi_{2}^{l}
\phi_{2}^{l'}\right]+\sum_{l'=1}^{r}{\bm{u}^{l}}^{T}_{0,2}\mathbb{E}_{\bm\xi_2}\left[
\bm{K}_{2}\phi_{2}^{l}\phi_{2}^{l'}\right]
\bm{u}_{0,2}^{l'}\phi_{1}^{l'}(\bm\xi_1)\nonumber \\
&&= {\bm{u}^{l}}^{T}_{0,1}\bm{f}_{1}\mathbb{E}_{\bm\xi_2}
\left[\phi_{2}^{l}\right] + {\bm{u}^{l}}^{T}_{0,2}\bm{f}_{2}\mathbb{E}_{\bm\xi_2}
\left[\phi_{2}^{l}\right], \qquad l=1,\dots,r.
\end{eqnarray}

It is straightforward to check that (\ref{eqn:lin_sys_1_rg}) -- with its components given in (\ref{eqn:short_1_rg_update}) --  is the matrix representation of  (\ref{eqn:phi_1_var_rg_rep2}). A similar approach can be used to derive (\ref{eqn:lin_sys_2_rg}).

\section{}
\label{apx:Wiener_chaos_coefficients}
Let $\kappa(\bm x,\bm \xi)=\exp(G(\bm x,\bm \xi))$ be a lognormal random field, where $G(\bm x,\bm\xi)$ is a Gaussian random field given by the Karhunen-Lo\`eve expansion $G = \bar{G}+\sum_{j=1}^{d}\sqrt{\tau_j}g_j(\bm x)\xi_j$. Here, $\bar{G}$ is the mean of $G$, $\left\lbrace \tau_{j} \right\rbrace_{j=1}^{d_i}$ and $\left\lbrace g_{j}(\bm{x}) \right\rbrace_{j=1}^{d_i}$ are, respectively, $d$ largest eigenvalues and the corresponding eigenfunctions of the covariance function of $G$, and $\{\xi_j\}_{j=1}^d$ are i.i.d. normal Gaussian random variables. 

Following \citep{Ullmann08}, the coefficients ${\kappa}_{\bm{i}}(\bm{x})$ in the Hermite polynomial chaos expansion $\kappa(\bm{x},\bm{\xi}) \approx \sum_{\bm{j} \in \mathscr{I}_{d,p}} \kappa_{\bm{j}}(\bm{x}) \psi_{\bm{j}}(\bm{\xi})$ can be computed by
\begin{equation}
\label{eqn:WPC_coeff_formula}
{\kappa}_{\bm{i}}(\bm{x}) = \frac{\bar{\kappa}}{\sqrt{\bm{i}!}} \prod_{j=1}^d \left[ \sqrt{\tau_j} g_j(\bm{x}) \right]^{i_j},
\end{equation}
where $\bar{\kappa} = \mathrm{exp}\left[ \bar{G} + \frac{ \mathrm{var}[G] }{2}\right]$, and $\sqrt{\bm{i}!} = \prod_{j=1}^d i_j !$.

\section*{References}

\bibliographystyle{model2-names}
\bibliography{Bib_V1}

\begin{thebibliography}{67}
\expandafter\ifx\csname natexlab\endcsname\relax\def\natexlab#1{#1}\fi
\expandafter\ifx\csname url\endcsname\relax
  \def\url#1{\texttt{#1}}\fi
\expandafter\ifx\csname urlprefix\endcsname\relax\def\urlprefix{URL }\fi
\providecommand{\eprint}[2][]{\url{#2}}
\providecommand{\bibinfo}[2]{#2}
\ifx\xfnm\relax \def\xfnm[#1]{\unskip,\space#1}\fi
\bibitem[{A.~Falc{\'o} and Nouy(2012)}]{Falco11}
\bibinfo{author}{A.~Falc{\'o}, A.}, \bibinfo{author}{Nouy, A.},
  \bibinfo{year}{2012}.
\newblock \bibinfo{title}{Proper generalized decomposition for nonlinear convex
  problems in tensor banach spaces}.
\newblock \bibinfo{journal}{Numerische Mathematik} \bibinfo{volume}{121},
  \bibinfo{pages}{503--530}.
\bibitem[{A.~Toselli(2005)}]{Toselli05}
\bibinfo{author}{A.~Toselli, O.B.W.}, \bibinfo{year}{2005}.
\newblock \bibinfo{title}{Domain Decomposition Methods, Algorithms and Theory}.
\newblock \bibinfo{publisher}{Springer, Berlin}.
\bibitem[{Arnst et~al.(2012)Arnst, Ghanem, Phipps and Red-Horse}]{Arnst12}
\bibinfo{author}{Arnst, M.}, \bibinfo{author}{Ghanem, R.},
  \bibinfo{author}{Phipps, E.}, \bibinfo{author}{Red-Horse, J.},
  \bibinfo{year}{2012}.
\newblock \bibinfo{title}{Reduced chaos expansions with random coefficients in
  reduced-dimensional stochastic modeling of coupled problems}.
\newblock \bibinfo{journal}{arXiv preprint arXiv:1207.0910} .
\bibitem[{Babu{\v{s}}ka et~al.(2007)Babu{\v{s}}ka, Nobile and
  Tempone}]{Babuska07a}
\bibinfo{author}{Babu{\v{s}}ka, I.}, \bibinfo{author}{Nobile, F.},
  \bibinfo{author}{Tempone, R.}, \bibinfo{year}{2007}.
\newblock \bibinfo{title}{A stochastic collocation method for elliptic partial
  differential equations with random input data}.
\newblock \bibinfo{journal}{SIAM Journal on Numerical Analysis}
  \bibinfo{volume}{45}, \bibinfo{pages}{1005--1034}.
\bibitem[{Beylkin and Mohlenkamp(2002)}]{Beylkin02}
\bibinfo{author}{Beylkin, G.}, \bibinfo{author}{Mohlenkamp, M.},
  \bibinfo{year}{2002}.
\newblock \bibinfo{title}{{Numerical operator calculus in higher dimensions}}.
\newblock \bibinfo{journal}{Proceedings of the National Academy of Science}
  \bibinfo{volume}{99}, \bibinfo{pages}{10246--10251}.
\bibitem[{Bieri(2011)}]{Bieri11}
\bibinfo{author}{Bieri, M.}, \bibinfo{year}{2011}.
\newblock \bibinfo{title}{A sparse composite collocation finite element method
  for elliptic {SPDEs}.}
\newblock \bibinfo{journal}{SIAM Journal on Numerical Analysis}
  \bibinfo{volume}{49}, \bibinfo{pages}{2277--2301}.
\bibitem[{Bieri et~al.(2009)Bieri, Andreev and Schwab}]{Bieri09c}
\bibinfo{author}{Bieri, M.}, \bibinfo{author}{Andreev, R.},
  \bibinfo{author}{Schwab, C.}, \bibinfo{year}{2009}.
\newblock \bibinfo{title}{Sparse tensor discretization of elliptic {SPDEs}}.
\newblock \bibinfo{journal}{SIAM Journal on Scientific Computing}
  \bibinfo{volume}{31}, \bibinfo{pages}{4281--4304}.
\bibitem[{Bieri and Schwab(2009)}]{Bieri09a}
\bibinfo{author}{Bieri, M.}, \bibinfo{author}{Schwab, C.},
  \bibinfo{year}{2009}.
\newblock \bibinfo{title}{Sparse high order {FEM} for elliptic {sPDEs}}.
\newblock \bibinfo{journal}{Computer Methods in Applied Mechanics and
  Engineering} \bibinfo{volume}{198}, \bibinfo{pages}{1149--1170}.
\bibitem[{Cai(1993)}]{Cai93}
\bibinfo{author}{Cai, X.C.}, \bibinfo{year}{1993}.
\newblock \bibinfo{title}{Some nonoverlapping domain decomposition methods}.
\newblock \bibinfo{journal}{SIAM Journal on Scientific Computing}
  \bibinfo{volume}{14}, \bibinfo{pages}{239--247}.
\bibitem[{Cameron and Martin(1947)}]{Cameron47}
\bibinfo{author}{Cameron, R.}, \bibinfo{author}{Martin, W.},
  \bibinfo{year}{1947}.
\newblock \bibinfo{title}{The orthogonal development of non-linear functionals
  in series of fourier-hermite functionals}.
\newblock \bibinfo{journal}{The Annals of Mathematics} \bibinfo{volume}{48},
  \bibinfo{pages}{385--392}.
\bibitem[{Chan and Mathew(1990)}]{Chan90}
\bibinfo{author}{Chan, T.F.}, \bibinfo{author}{Mathew, T.P.},
  \bibinfo{year}{1990}.
\newblock \bibinfo{title}{Domain decomposition algorithms}.
\newblock \bibinfo{journal}{Third International Symposium on Domain
  Decomposition Methods for Partial Differential Equations, T. F. Chan, R.
  Glowinski, J. Periaux, and O. B. Widlund, eds., SIAM, Philadelphia} .
\bibitem[{Charmpis and Papadrakakis(2002)}]{Charmpis02}
\bibinfo{author}{Charmpis, D.C.}, \bibinfo{author}{Papadrakakis, M.},
  \bibinfo{year}{2002}.
\newblock \bibinfo{title}{Enhancing the performance of the {FETI} method with
  preconditioning techniques implemented on clusters of networked computers}.
\newblock \bibinfo{journal}{Computational Mechanics} \bibinfo{volume}{30},
  \bibinfo{pages}{12--28}.
\bibitem[{Chevreuil et~al.(2013)Chevreuil, Nouy and Safatly}]{Chevreuil12}
\bibinfo{author}{Chevreuil, M.}, \bibinfo{author}{Nouy, A.},
  \bibinfo{author}{Safatly, E.}, \bibinfo{year}{2013}.
\newblock \bibinfo{title}{A multiscale method with patch for the solution of
  stochastic partial differential equations with localized uncertainties}.
\newblock \bibinfo{journal}{Computer Methods in Applied Mechanics and
  Engineering} \bibinfo{volume}{255}, \bibinfo{pages}{255--274}.
\bibitem[{Chinesta et~al.(2011)Chinesta, Ammar, Leygue and
  Keunings}]{Chinesta11}
\bibinfo{author}{Chinesta, F.}, \bibinfo{author}{Ammar, A.},
  \bibinfo{author}{Leygue, A.}, \bibinfo{author}{Keunings, R.},
  \bibinfo{year}{2011}.
\newblock \bibinfo{title}{An overview of the proper generalized decomposition
  with applications in computational rheology}.
\newblock \bibinfo{journal}{Journal of NonNewtonian Fluid Mechanics}
  \bibinfo{volume}{166}, \bibinfo{pages}{578--592}.
\bibitem[{Cottereau et~al.(2011)Cottereau, Clouteau, Dhia and
  Zaccardi}]{Cottereau11}
\bibinfo{author}{Cottereau, R.}, \bibinfo{author}{Clouteau, D.},
  \bibinfo{author}{Dhia, H.B.}, \bibinfo{author}{Zaccardi, C.},
  \bibinfo{year}{2011}.
\newblock \bibinfo{title}{A stochastic-deterministic coupling method for
  continuum mechanics}.
\newblock \bibinfo{journal}{Computer Methods in Applied Mechanics and
  Engineering} \bibinfo{volume}{200}, \bibinfo{pages}{3280--3288}.
\bibitem[{Dohrmann(2003)}]{Dohrmann03}
\bibinfo{author}{Dohrmann, C.R.}, \bibinfo{year}{2003}.
\newblock \bibinfo{title}{A preconditioner for substructuring based on
  constrained energy minimization}.
\newblock \bibinfo{journal}{SIAM Journal on Scientific Computing}
  \bibinfo{volume}{25}, \bibinfo{pages}{246--258}.
\bibitem[{Doostan et~al.(2007)Doostan, Ghanem and Red-Horse}]{Doostan07}
\bibinfo{author}{Doostan, A.}, \bibinfo{author}{Ghanem, R.},
  \bibinfo{author}{Red-Horse, J.}, \bibinfo{year}{2007}.
\newblock \bibinfo{title}{Stochastic model reduction for chaos
  representations}.
\newblock \bibinfo{journal}{Computer Methods in Applied Mechanics and
  Engineering} \bibinfo{volume}{196}, \bibinfo{pages}{3951--3966}.
\bibitem[{Doostan and Iaccarino(2009)}]{Doostan09}
\bibinfo{author}{Doostan, A.}, \bibinfo{author}{Iaccarino, G.},
  \bibinfo{year}{2009}.
\newblock \bibinfo{title}{A least-squares approximation of partial differential
  equations with high-dimensional random inputs}.
\newblock \bibinfo{journal}{Journal of Computational Physics}
  \bibinfo{volume}{228}, \bibinfo{pages}{4332--4345}.
\bibitem[{Doostan and Owhadi(2011)}]{Doostan11a}
\bibinfo{author}{Doostan, A.}, \bibinfo{author}{Owhadi, H.},
  \bibinfo{year}{2011}.
\newblock \bibinfo{title}{A non-adapted sparse approximation of {PDEs} with
  stochastic inputs}.
\newblock \bibinfo{journal}{Journal of Computational Physics}
  \bibinfo{volume}{230}, \bibinfo{pages}{3015--3034}.
\bibitem[{Doostan et~al.(2013)Doostan, Validi and Iaccarino}]{Doostan13}
\bibinfo{author}{Doostan, A.}, \bibinfo{author}{Validi, A.},
  \bibinfo{author}{Iaccarino, G.}, \bibinfo{year}{2013}.
\newblock \bibinfo{title}{Non-intrusive low-rank separated approximation of
  high-dimensional stochastic models}.
\newblock \bibinfo{journal}{Computer Methods in Applied Mechanics and
  Engineering} \bibinfo{note}{In press}.
\bibitem[{Farhat et~al.(1995)Farhat, Chen and Mandel}]{Farhat95}
\bibinfo{author}{Farhat, C.}, \bibinfo{author}{Chen, P.S.},
  \bibinfo{author}{Mandel, J.}, \bibinfo{year}{1995}.
\newblock \bibinfo{title}{A scalable {L}agrange multiplier based domain
  decomposition method for time-dependent problems}.
\newblock \bibinfo{journal}{International Journal for Numerical Methods in
  Engineering} \bibinfo{volume}{38}, \bibinfo{pages}{3831--3853}.
\bibitem[{Farhat et~al.(2001)Farhat, Lesoinne, LeTallec, Pierson and
  Rixen}]{Farhat01}
\bibinfo{author}{Farhat, C.}, \bibinfo{author}{Lesoinne, M.},
  \bibinfo{author}{LeTallec, P.}, \bibinfo{author}{Pierson, K.},
  \bibinfo{author}{Rixen, D.}, \bibinfo{year}{2001}.
\newblock \bibinfo{title}{{FETI-DP}: a dual-primal unified {FETI} method--part
  i: a faster alternative to the two-level {FETI} method}.
\newblock \bibinfo{journal}{International Journal for Numerical Methods in
  Engineering} \bibinfo{volume}{50}, \bibinfo{pages}{1523--1544}.
\bibitem[{Farhat and Mandel(1998)}]{Farhat98}
\bibinfo{author}{Farhat, C.}, \bibinfo{author}{Mandel, J.},
  \bibinfo{year}{1998}.
\newblock \bibinfo{title}{The two-level {FETI} method for static and dynamic
  plate problems—part i: An optimal iterative solver for biharmonic systems}.
\newblock \bibinfo{journal}{Computer Methods in Applied Mechanics and
  Engineering} \bibinfo{volume}{155}, \bibinfo{pages}{129--152}.
\bibitem[{Farhat and Roux(1991)}]{Farhat91}
\bibinfo{author}{Farhat, C.}, \bibinfo{author}{Roux, F.}, \bibinfo{year}{1991}.
\newblock \bibinfo{title}{A method of finite element tearing and
  interconnecting and its parallel solution algorithm}.
\newblock \bibinfo{journal}{International Journal for Numerical Methods in
  Engineering} \bibinfo{volume}{32}, \bibinfo{pages}{1205--1227}.
\bibitem[{Foo and Karniadakis(2010)}]{Foo10}
\bibinfo{author}{Foo, J.}, \bibinfo{author}{Karniadakis, G.},
  \bibinfo{year}{2010}.
\newblock \bibinfo{title}{Multi-element probabilistic collocation method in
  high dimensions}.
\newblock \bibinfo{journal}{Journal of Computational Physics}
  \bibinfo{volume}{229}, \bibinfo{pages}{1536--1557}.
\bibitem[{Franca and Macedo(1998)}]{Franca98}
\bibinfo{author}{Franca, L.}, \bibinfo{author}{Macedo, A.},
  \bibinfo{year}{1998}.
\newblock \bibinfo{title}{A two-level finite element method and its application
  to the {H}elmholtz equation}.
\newblock \bibinfo{journal}{International Journal for Numerical Methods in
  Engineering} \bibinfo{volume}{43}, \bibinfo{pages}{23--32}.
\bibitem[{Galvis and Sarkis(2010)}]{Galvis10}
\bibinfo{author}{Galvis, J.}, \bibinfo{author}{Sarkis, M.},
  \bibinfo{year}{2010}.
\newblock \bibinfo{title}{{FETI} and {BDD} preconditioners for
  {Stokes-Mortar-Darcy} systems}.
\newblock \bibinfo{journal}{Commun. Appl. Math. Comput. Sci}
  \bibinfo{volume}{5}, \bibinfo{pages}{1--30}.
\bibitem[{Gao and Hesthaven(2010)}]{Gao10b}
\bibinfo{author}{Gao, Z.}, \bibinfo{author}{Hesthaven, J.},
  \bibinfo{year}{2010}.
\newblock \bibinfo{title}{On anova expansions and strategies for choosing the
  anchor point}.
\newblock \bibinfo{journal}{Applied Mathematics and Computation}
  \bibinfo{volume}{217}, \bibinfo{pages}{3274--3285}.
\bibitem[{Gao and Hesthaven(2011)}]{Gao11a}
\bibinfo{author}{Gao, Z.}, \bibinfo{author}{Hesthaven, J.},
  \bibinfo{year}{2011}.
\newblock \bibinfo{title}{Efficient solution of ordinary differential equations
  with high-dimensional parametrized uncertainty}.
\newblock \bibinfo{journal}{Communications in Computational Physics}
  \bibinfo{volume}{10}, \bibinfo{pages}{253}.
\bibitem[{Geradin et~al.(1997)Geradin, Coulon and Delsemme}]{Geradin97}
\bibinfo{author}{Geradin, M.}, \bibinfo{author}{Coulon, D.},
  \bibinfo{author}{Delsemme, J.}, \bibinfo{year}{1997}.
\newblock \bibinfo{title}{Parallelization of the {SAMCEF} finite element
  software through domain decomposition and {FETI} algorithm}.
\newblock \bibinfo{journal}{International Journal of Supercomputer
  Applications} \bibinfo{volume}{11}, \bibinfo{pages}{286--298}.
\bibitem[{Ghanem(1999)}]{Ghanem99c}
\bibinfo{author}{Ghanem, R.}, \bibinfo{year}{1999}.
\newblock \bibinfo{title}{Higher order sensitivity of heat conduction problems
  to random data using the spectral stochastic finite element method}.
\newblock \bibinfo{journal}{ASME Journal of Heat Transfer}
  \bibinfo{volume}{121}, \bibinfo{pages}{290--299}.
\bibitem[{Ghanem and Spanos()}]{Ghanem91a}
\bibinfo{author}{Ghanem, R.}, \bibinfo{author}{Spanos, P.}, .
\newblock \bibinfo{title}{Stochastic Finite Elements: A Spectral Approach}.
\newblock \bibinfo{publisher}{Springer Verlag}, \bibinfo{address}{Berlin}.
\bibitem[{Ghosh et~al.(2009)Ghosh, Avery and Farhat}]{Ghosh09}
\bibinfo{author}{Ghosh, D.}, \bibinfo{author}{Avery, P.},
  \bibinfo{author}{Farhat, C.}, \bibinfo{year}{2009}.
\newblock \bibinfo{title}{A {FETI}-preconditioned congugate gradient method for
  large-scale stochastic finite element problems}.
\newblock \bibinfo{journal}{International Journal For Numerical Methods In
  Engineering} \bibinfo{volume}{80}, \bibinfo{pages}{914--931}.
\bibitem[{Gosselet and Rey(2006)}]{Gosselet06}
\bibinfo{author}{Gosselet, P.}, \bibinfo{author}{Rey, C.},
  \bibinfo{year}{2006}.
\newblock \bibinfo{title}{Non-overlapping domain decomposition methods in
  structural mechanics}.
\newblock \bibinfo{journal}{Archives of Computational Methods in Engineering}
  \bibinfo{volume}{13}, \bibinfo{pages}{515--572}.
\bibitem[{Hackbusch(2012)}]{Hackbusch12}
\bibinfo{author}{Hackbusch, W.}, \bibinfo{year}{2012}.
\newblock \bibinfo{title}{Tensor spaces and numerical tensor calculus}.
  volume~\bibinfo{volume}{42}.
\newblock \bibinfo{publisher}{Springer}.
\bibitem[{Khoromskij and Schwab(2011)}]{Khoromskij11}
\bibinfo{author}{Khoromskij, B.}, \bibinfo{author}{Schwab, C.},
  \bibinfo{year}{2011}.
\newblock \bibinfo{title}{Tensor-structured {G}alerkin approximation of
  parametric and stochastic elliptic {PDEs}}.
\newblock \bibinfo{journal}{SIAM Journal on Scientific Computing}
  \bibinfo{volume}{33}, \bibinfo{pages}{364--385}.
\bibitem[{Klawonn and Widlund(2001)}]{Klawonn01}
\bibinfo{author}{Klawonn, A.}, \bibinfo{author}{Widlund, O.},
  \bibinfo{year}{2001}.
\newblock \bibinfo{title}{{FETI} and {N}eumann-{N}eumann iterative
  substructuring methods: Connections and new results}.
\newblock \bibinfo{journal}{Communications on Pure and Applied Mathematics}
  \bibinfo{volume}{54}, \bibinfo{pages}{57--90}.
\bibitem[{LeTallec et~al.(1991)LeTallec, Roeck and Vidrascu}]{Tallec91}
\bibinfo{author}{LeTallec, P.}, \bibinfo{author}{Roeck, T.D.},
  \bibinfo{author}{Vidrascu, M.}, \bibinfo{year}{1991}.
\newblock \bibinfo{title}{Domain decomposition methods for large linearly
  elliptic three-dimensional problems}.
\newblock \bibinfo{journal}{Journal of Computational and Applied Mathematics}
  \bibinfo{volume}{34}, \bibinfo{pages}{93--117}.
\bibitem[{Logg et~al.(2012)Logg, Mardal and Wells}]{Logg12}
\bibinfo{author}{Logg, A.}, \bibinfo{author}{Mardal, K.A.},
  \bibinfo{author}{Wells, G.}, \bibinfo{year}{2012}.
\newblock \bibinfo{title}{Automated Solution of Differential Equations by the
  Finite Element Method}.
\newblock \bibinfo{publisher}{Springer}.
\bibitem[{Ma and Zabaras(2009)}]{Ma09a}
\bibinfo{author}{Ma, X.}, \bibinfo{author}{Zabaras, N.}, \bibinfo{year}{2009}.
\newblock \bibinfo{title}{An adaptive hierarchical sparse grid collocation
  algorithm for the solution of stochastic differential equations}.
\newblock \bibinfo{journal}{Journal of Computational Physics}
  \bibinfo{volume}{228}, \bibinfo{pages}{3084--3113}.
\bibitem[{Maitre and Knio(2010)}]{LeMaitre10}
\bibinfo{author}{Maitre, O.L.}, \bibinfo{author}{Knio, O.},
  \bibinfo{year}{2010}.
\newblock \bibinfo{title}{Spectral Methods for Uncertainty Quantification with
  Applications to Computational Fluid Dynamics}.
\newblock \bibinfo{publisher}{Springer}.
\bibitem[{Maitre et~al.(2004)Maitre, Knio, Najm and Ghanem}]{LeMaitre04}
\bibinfo{author}{Maitre, O.L.}, \bibinfo{author}{Knio, O.},
  \bibinfo{author}{Najm, H.}, \bibinfo{author}{Ghanem, R.},
  \bibinfo{year}{2004}.
\newblock \bibinfo{title}{Uncertainty propagation using {Wiener-Haar}
  expansions}.
\newblock \bibinfo{journal}{Journal of Computational Physics}
  \bibinfo{volume}{197}, \bibinfo{pages}{28--57}.
\bibitem[{Mathelin and Hussaini(2003)}]{Mathelin03}
\bibinfo{author}{Mathelin, L.}, \bibinfo{author}{Hussaini, M.},
  \bibinfo{year}{2003}.
\newblock \bibinfo{title}{A Stochastic Collocation Algorithm for Uncertainty
  Analysis}.
\newblock \bibinfo{type}{Technical Report} \bibinfo{number}{NAS 1.26:212153;
  NASA/CR-2003-212153}. NASA Langley Research Center.
\bibitem[{Matthies(2008)}]{HGMatth08}
\bibinfo{author}{Matthies, H.G.}, \bibinfo{year}{2008}.
\newblock \bibinfo{title}{Stochastic finite elements: Computational approaches
  to stochastic partial differential equations}.
\newblock \bibinfo{journal}{Z. Angew. Math. Mech.} \bibinfo{volume}{88},
  \bibinfo{pages}{849--873}.
\bibitem[{Najm(2009)}]{Najm09}
\bibinfo{author}{Najm, H.}, \bibinfo{year}{2009}.
\newblock \bibinfo{title}{Uncertainty quantification and polynomial chaos
  techniques in computational fluid dynamics}.
\newblock \bibinfo{journal}{Annual Reviews} \bibinfo{volume}{41},
  \bibinfo{pages}{35--52}.
\bibitem[{Nobile et~al.(2008)Nobile, Tempone and Webster}]{Nobile08b}
\bibinfo{author}{Nobile, F.}, \bibinfo{author}{Tempone, R.},
  \bibinfo{author}{Webster, C.}, \bibinfo{year}{2008}.
\newblock \bibinfo{title}{An anisotropic sparse grid stochastic collocation
  method for partial differential equations with random input data}.
\newblock \bibinfo{journal}{SIAM Journal on Numerical Analysis}
  \bibinfo{volume}{46}, \bibinfo{pages}{2411--2442}.
\bibitem[{Nouy(2007)}]{Nouy07}
\bibinfo{author}{Nouy, A.}, \bibinfo{year}{2007}.
\newblock \bibinfo{title}{{A generalized spectral decomposition technique to
  solve a class of linear stochastic partial differential equations}}.
\newblock \bibinfo{journal}{Computer Methods in Applied Mechanics and
  Engineering} \bibinfo{volume}{196}, \bibinfo{pages}{4521--4537}.
\bibitem[{Nouy(2008)}]{Nouy08}
\bibinfo{author}{Nouy, A.}, \bibinfo{year}{2008}.
\newblock \bibinfo{title}{{Generalized spectral decomposition method for
  solving stochastic finite element equations: Invariant subspace problem and
  dedicated algorithms}}.
\newblock \bibinfo{journal}{Computer Methods in Applied Mechanics and
  Engineering} \bibinfo{volume}{197}, \bibinfo{pages}{4718--4736}.
\bibitem[{Nouy(2010)}]{Nouy_2010}
\bibinfo{author}{Nouy, A.}, \bibinfo{year}{2010}.
\newblock \bibinfo{title}{Proper generalized decompositions and separated
  representations for the numerical solution of high dimensional stochastic
  problems}.
\newblock \bibinfo{journal}{Archives of Computational Methods in Engineering}
  \bibinfo{volume}{17}, \bibinfo{pages}{403--434}.
\bibitem[{Park and Felippa(1998)}]{Park98}
\bibinfo{author}{Park, K.C.}, \bibinfo{author}{Felippa, C.A.},
  \bibinfo{year}{1998}.
\newblock \bibinfo{title}{A variational framework for solution method
  developments in structural mechanics}.
\newblock \bibinfo{journal}{Journal of Applied Mechanics}
  \bibinfo{volume}{56/1}, \bibinfo{pages}{242--249}.
\bibitem[{Quarteroni and Valli(1999)}]{Quarteroni99}
\bibinfo{author}{Quarteroni, A.}, \bibinfo{author}{Valli, A.},
  \bibinfo{year}{1999}.
\newblock \bibinfo{title}{Domain decomposition methods for partial differential
  equations}. volume~\bibinfo{volume}{10}.
\newblock \bibinfo{publisher}{Clarendon Press Oxford}.
\bibitem[{Rapetti and Toselli(2001)}]{Rapetti01}
\bibinfo{author}{Rapetti, F.}, \bibinfo{author}{Toselli, A.},
  \bibinfo{year}{2001}.
\newblock \bibinfo{title}{A {FETI} preconditioner for two dimensional edge
  element approximations of {Maxwell's} equations on nonmatching grids}.
\newblock \bibinfo{journal}{SIAM Journal on Scientific Computing}
  \bibinfo{volume}{23}, \bibinfo{pages}{92--108}.
\bibitem[{Rixen et~al.(1999)Rixen, Farhat, Tezaur and Mandel}]{Farhat99}
\bibinfo{author}{Rixen, D.J.}, \bibinfo{author}{Farhat, C.},
  \bibinfo{author}{Tezaur, R.}, \bibinfo{author}{Mandel, J.},
  \bibinfo{year}{1999}.
\newblock \bibinfo{title}{Theoretical comparison of the {FETI} and
  algebraically partitioned {FETI} methods, and performance comparisons with a
  direct sparse solver}.
\newblock \bibinfo{journal}{International Journal for Numerical Methods in
  Engineering} \bibinfo{volume}{46}, \bibinfo{pages}{501--533}.
\bibitem[{Smith et~al.(2004)Smith, Bjorstad and Gropp}]{Smith04}
\bibinfo{author}{Smith, B.}, \bibinfo{author}{Bjorstad, P.},
  \bibinfo{author}{Gropp, W.}, \bibinfo{year}{2004}.
\newblock \bibinfo{title}{Domain decomposition}.
\newblock \bibinfo{publisher}{Cambridge University Press}.
\bibitem[{Soize and Ghanem(2009)}]{Soize09}
\bibinfo{author}{Soize, C.}, \bibinfo{author}{Ghanem, R.G.},
  \bibinfo{year}{2009}.
\newblock \bibinfo{title}{Reduced chaos decomposition with random coefficients
  of vector-valued random variables and random fields}.
\newblock \bibinfo{journal}{Computer Methods in Applied Mechanics and
  Engineering} \bibinfo{volume}{198}, \bibinfo{pages}{1926--1934}.
\bibitem[{Subber and Sarkar(2012)}]{Subber12}
\bibinfo{author}{Subber, W.}, \bibinfo{author}{Sarkar, A.},
  \bibinfo{year}{2012}.
\newblock \bibinfo{title}{Domain decomposition method of stochastic pdes: a
  two-level scalable preconditioner}.
\newblock \bibinfo{journal}{Journal of Physics: Conference Series}
  \bibinfo{volume}{341}, \bibinfo{pages}{012033}.
\bibitem[{Todor and Schwab(2007)}]{Todor07a}
\bibinfo{author}{Todor, R.A.}, \bibinfo{author}{Schwab, C.},
  \bibinfo{year}{2007}.
\newblock \bibinfo{title}{Convergence rates for sparse chaos approximations of
  elliptic problems with stochastic coefficients}.
\newblock \bibinfo{journal}{IMA Journal of Numerical Analysis}
  \bibinfo{volume}{27}, \bibinfo{pages}{232--261}.
\bibitem[{Toselli and Klawonn(1999)}]{Toselli99a}
\bibinfo{author}{Toselli, A.}, \bibinfo{author}{Klawonn, A.},
  \bibinfo{year}{1999}.
\newblock \bibinfo{title}{A FETI Domain Decomposition Method For {M}axwell's
  Equations With Discontinuous Coefficients In Two Dimensions}.
\newblock \bibinfo{type}{Technical Report}. Courant Institute, New York
  University.
\bibitem[{Ullmann(2008)}]{Ullmann08}
\bibinfo{author}{Ullmann, E.}, \bibinfo{year}{2008}.
\newblock \bibinfo{title}{Solution Strategies for Stochastic Finite Element
  Discretizations}.
\newblock Ph.D. thesis. Technische Universit{\"a}t Bergakademie Freiberg.
\bibitem[{Wiener(1938)}]{Wiener38}
\bibinfo{author}{Wiener, N.}, \bibinfo{year}{1938}.
\newblock \bibinfo{title}{The homogeneous chaos}.
\newblock \bibinfo{journal}{Amer. J. Math} \bibinfo{volume}{60},
  \bibinfo{pages}{897--936}.
\bibitem[{Xiu(2009)}]{Xiu09a}
\bibinfo{author}{Xiu, D.}, \bibinfo{year}{2009}.
\newblock \bibinfo{title}{Fast numerical methods for stochastic computations: A
  review}.
\newblock \bibinfo{journal}{Communications in Computational Physics}
  \bibinfo{volume}{5}, \bibinfo{pages}{242--272}.
\bibitem[{Xiu(2010)}]{Xiu10a}
\bibinfo{author}{Xiu, D.}, \bibinfo{year}{2010}.
\newblock \bibinfo{title}{Numerical Methods for Stochastic Computations: A
  Spectral Method Approach}.
\newblock \bibinfo{publisher}{Princeton University Press}.
\bibitem[{Xiu and Hesthaven(2005)}]{Xiu05a}
\bibinfo{author}{Xiu, D.}, \bibinfo{author}{Hesthaven, J.},
  \bibinfo{year}{2005}.
\newblock \bibinfo{title}{High-order collocation methods for differential
  equations with random inputs}.
\newblock \bibinfo{journal}{SIAM Journal on Scientific Computing}
  \bibinfo{volume}{27}, \bibinfo{pages}{1118--1139}.
\bibitem[{Xiu and Karniadakis(2002)}]{Xiu02}
\bibinfo{author}{Xiu, D.}, \bibinfo{author}{Karniadakis, G.},
  \bibinfo{year}{2002}.
\newblock \bibinfo{title}{The {W}iener-{A}skey polynomial chaos for stochastic
  differential equations}.
\newblock \bibinfo{journal}{SIAM Joural on Scientific Computing}
  \bibinfo{volume}{24}, \bibinfo{pages}{619--644}.
\bibitem[{Xiu and Karniadakis(2003)}]{Xiu03}
\bibinfo{author}{Xiu, D.}, \bibinfo{author}{Karniadakis, G.},
  \bibinfo{year}{2003}.
\newblock \bibinfo{title}{Modeling uncertainty in flow simulations via
  generalized polynomial chaos}.
\newblock \bibinfo{journal}{Journal of Computational Physics}
  \bibinfo{volume}{187}, \bibinfo{pages}{137--167}.
\bibitem[{Xu and Zou(1998)}]{Xu98}
\bibinfo{author}{Xu, J.}, \bibinfo{author}{Zou, J.}, \bibinfo{year}{1998}.
\newblock \bibinfo{title}{Some nonoverlapping domain decomposition methods}.
\newblock \bibinfo{journal}{SIAM Review.} \bibinfo{volume}{40},
  \bibinfo{pages}{857--914}.
\bibitem[{Zhang et~al.(2009)Zhang, Choi and Karniadakis}]{Zhang11}
\bibinfo{author}{Zhang, Z.}, \bibinfo{author}{Choi, M.},
  \bibinfo{author}{Karniadakis, G.E.}, \bibinfo{year}{2009}.
\newblock \bibinfo{title}{Anchor points matter in anova decomposition}, in:
  \bibinfo{booktitle}{Spectral and Higher Order Methods for Partial
  Differential Equations, Lecture Notes in Computational Science and
  Engineering}, pp. \bibinfo{pages}{347--355}.

\end{thebibliography}

\end{document}